\documentclass[12pt,leqno]{amsbook}
\usepackage{latexsym,amssymb,amscd,fancyhdr,eqname,smartref,graphicx,amsopn}
\usepackage[Sonny]{fncychap}

\def\B'c{{\mathcal{B'}}}
\def\U'c{{\mathcal{U'}}}

\def\opn#1#2{\def#1{\operatorname{#2}}} 
\def\opn#1#2{\def#1{\operatorname{#2}}}
\opn\chara{char}
\opn\length{\ell}
\opn\projdim{proj\,dim}
\opn\injdim{inj\,dim}
\opn\ini{in}
\opn\rank{rank}
\opn\depth{depth}
\opn\height{ht}
\opn\embdim{emb\,dim}
\opn\codim{codim}

\opn\Tr{Tr}
\opn\bigrank{big\,rank}
\opn\superheight{superheight}\opn\lcm{lcm}
\opn\trdeg{tr\,deg}%
\opn\reg{reg}
\opn\lreg{lreg}
\opn\set{set}
\opn\supp{Supp}
\opn\shad{Shad}
\opn\del{del}
\opn\div{div}
\opn\Div{Div}
\opn\cl{cl}
\opn\Cl{Cl}
\opn\Spec{Spec}
\opn\Supp{Supp}
\opn\supp{supp}
\opn\Sing{Sing}
\opn\Ass{Ass}
\opn\Ann{Ann}
\opn\Rad{Rad}
\opn\Soc{Soc}
\opn\depth{depth}
\opn\Ker{Ker}
\opn\Coker{Coker}
\opn\Im{Im}
\opn\Hom{Hom}
\opn\Tor{Tor}
\opn\Ext{Ext}
\opn\End{End}
\opn\Aut{Aut}
\opn\id{id}

\opn\nat{nat}
\opn\GL{GL}
\opn\SL{SL}
\opn\mod{mod}
\opn\ord{ord}
\opn\id{id}
\opn\chr{char}
\opn\aff{aff}
\opn\con{conv}
\opn\relint{relint}
\opn\st{st}
\opn\lk{lk}
\opn\cn{cn}
\opn\core{core}
\opn\vol{vol}
\opn\gr{gr}

\def\pot#1#2{#1[\kern-0.28ex[#2]\kern-0.28ex]}

\opn\dirlim{\underrightarrow{\lim}}
\opn\invlim{\underleftarrow{\lim}}
\let\union=\cup

\let\dirsum=\oplus

\def\pnt{{\raise0.5mm\hbox{\large\bf.}}}

\def\twoline#1#2{\aoverb{\scriptstyle {#1}}{\scriptstyle {#2}}}
\newcommand{\aoverb}[2]{{\genfrac{}{}{0pt}{1}{#1}{#2}}}

\def\Implies{\ifmmode\Longrightarrow \else
     \unskip${}\Longrightarrow{}$\ignorespaces\fi}
\def\implies{\ifmmode\Rightarrow \else
     \unskip${}\Rightarrow{}$\ignorespaces\fi}
\def\iff{\ifmmode\Longleftrightarrow \else
     \unskip${}\Longleftrightarrow{}$\ignorespaces\fi}

\let\:=\colon

\newtheorem{Theorem1}{\bfseries Theorem}[section]
\newtheorem{Lemma1}[Theorem1]{\bfseries Lemma}
\newtheorem{Corollary1}[Theorem1]{\bfseries Corollary}
\newtheorem{Proposition1}[Theorem1]{\bfseries Proposition}
\newtheorem{Remark1}[Theorem1]{\bfseries Remark}

\newtheorem{Example1}[Theorem1]{\bfseries Example}
\newtheorem{Examples1}[Theorem1]{\bfseries Examples}
\newtheorem{Definition1}[Theorem1]{\bfseries Definition}

\let\epsilon=\varepsilon
\let\phi=\varphi
\let\kappa=\varkappa
  \renewcommand{\headrulewidth}{0pt} 
  \renewcommand{\footrulewidth}{0pt} 

\numberwithin{section}{chapter}
\numberwithin{equation}{chapter}
\numberwithin{figure}{chapter}
\title{Classes of monomial ideals}
\author{}
\date{}

\begin{document}
\begin{titlepage}

\pagestyle{empty}
\vspace{0in}
\begin{center}
\large\textbf{     UNIVERSITATEA "OVIDIUS" CONSTAN\c TA} 
	
\noindent\large\textbf{FACULTATEA DE MATEMATIC\u A \c SI INFORMATIC\u A } 	

\large\textbf{\c SCOALA DOCTORAL\u A}

\vspace{1.5in}

\huge{\textbf{TEZ\u A DE DOCTORAT}}
\end{center}
\vspace{2in}
\begin{flushleft}
	\large{\textbf{$\ \ $CONDUC\u ATOR \c STIIN\c TIFIC}
	
		\textbf{PROF. UNIV. DR. MIRELA \c STEF\u ANESCU}	}
\end{flushleft}

\vspace{0.5in}
\begin{flushright}
	\[
\]
\large{\textbf{DOCTORAND$\ \ $\qquad\qquad}

	\textbf{ANDA-GEORGIANA OLTEANU}}
\end{flushright}

\vspace{0.5in}
\vspace{0.5in}
\vspace{0.5in}
\begin{center}
\large{\textbf{CONSTAN\c TA 2008}}
\end{center}
\end{titlepage}
\newpage
\large\textbf{    } 
	
\noindent\large\textbf{ } 	

\large\textbf{ }

\vspace{2in}
\begin{center}
\LARGE\textbf{{CLASSES OF MONOMIAL IDEALS}}
\end{center}
\vspace{2in}
\begin{flushleft}
	\large{\textbf{$\ \ $}
	
		\textbf{}	}
\end{flushleft}

\vspace{0.5in}
\begin{flushright}
\large{\textbf{}

	\textbf{}}
\end{flushright}
\normalsize

\pagestyle{empty}
\chapter*{Acknowledgments}
\thispagestyle{empty}
I am very grateful to my advisor Professor Mirela \c Stef\u anescu for her 
support and guidance provided during my PhD study and research. 

I would like to thank Professor
J\"urgen Herzog and Professor Ezra Miller for valuable suggestions and discussions during the preparation of this thesis.

I am also deeply grateful to Professor Viviana Ene for the invaluable support, permanent encouragement, and good advices.

Not in the last place, I would like to thank my family for helping and sustaining me all the time. 
 
During the preparation of my thesis, I was supported 
by the CEEX Program of the Romanian Ministry of Education and Research, contract CEX 05-D11-11/2005 and by the CNCSIS contract TD 507/2007. 

\nopagebreak
\pagestyle{empty}
\tableofcontents

\frontmatter
\addcontentsline{toc}{chapter}{Preface} 
 
\chapter*{Preface}
\pagestyle{plain}
Monomial ideals are at the intersection of commutative algebra and combinatorial algebra. Many important problems in polynomial rings can be reduced to the study of monomial ideals. It became a standard method in commutative algebra to study invariants of arbitrary graded ideals in a polynomial ring by passing to their initial ideals with respect to some term order. By this process, one obtains a monomial ideal associated to a given graded ideal, which shares many basic invariants with the original ideal, but it is also accesible to combinatorial techniques and faster computational methods. A prominent example in this context is provided by the study of Hilbert functions for standard graded algebras over a field $k$. Lexsegment ideals play a key role in the discussion. The combinatorial characterization of Hilbert functions of graded ideals is due to F.S. Macaulay \cite{Ma}. Its square-free analogue -- the Kruskal-Katona theorem -- describes the possible $f$--vectors of the simplicial complexes on a given vertex set. Based on the work of S. Eliahou and M. Kervaire \cite{EK} on stable monomial ideals and using basic techniques on generic initial ideals, in $1993$, A. Bigatti \cite{B} and H. Hulett \cite{Hul} independently proved that, among all the graded ideals with a given Hilbert function, the lexsegment ideal posseses the maximal graded Betti numbers, provided the base field is of characteristic zero. This theorem has a square-free analogue due to A. Aramova, J. Herzog and T. Hibi \cite{AHH}.

Combinatorial algebra is closely related to square-free monomial ideals. Using the polarization process, to any monomial ideal, a square-free monomial ideal in a polynomial ring with more variables can be associated. The homological behaviour of a monomial ideal is preserved under passing to its polarization which is a square-free monomial ideal. To square-free monomial ideals, one may attach simplicial complexes and viceversa. On one hand, this allows the study of simplicial complexes by algebraic methods. This fact was impressively demonstrated by R. Stanley \cite{St1} who proved by this method the upper bound conjecture. On the other hand, one may apply combinatorial techniques to simplicial complexes in order to study the square-free monomial ideals.

The Alexander duality is one of the most efficient tools to study square-free monomial ideals nowadays. J.A. Eagon and V. Reiner \cite{EaRe} showed that the Stanley--Reisner ring of a simplicial complex $\Delta$ is Cohen--Macaulay if and only if the Stanley--Reisner ideal of the Alexander dual $I_{\Delta^{\vee}}$ has a linear resolution. This theorem has many interesting applications.

In this thesis, we focus on the study of some classes of monomial ideals, namely lexsegment ideals and monomial ideals with linear quotients. We define a new class of monomial ideals which, in the square-free case, is related to pure constructible simplicial complexes.

The thesis is structured in four chapters. Chapter 1 starts with regular sequences and their connections with the depth of a module. Next, we look at some properties of Cohen--Macaulay ideals. We recall some basic definitions and known facts about monomial ideals with linear quotients, stable ideals and lexsegment ideals. Next, we pay attention to the Koszul complex, to the Eliahou-Kervaire resolution for stable ideals, and to the construction of a resolution by mapping cones in the special context of monomial ideals with linear quotients. Next we give a brief description of simplicial complexes. Finally, we remind some properties of Coxeter groups and of subword complexes in Coxeter groups. These properties will be needed in the last chapter of the thesis.

Chapter $2$ is devoted to the study of lexsegment ideals. Initial lexsegment ideals are useful in extremal combinatorics and in the theory of Hilbert functions. Arbitrary lexsegment ideals were defined by H. Hulett and H.M. Martin \cite{HM}. A. Aramova, J. Herzog and E. De Negri also studied these ideals \cite{ADH} and \cite{DH}. They characterized all lexsegment ideals which have a linear resolution.

Let $S=k[x_1,\ldots, x_n]$ be the polynomial ring in $n$ variables over a field $k$. We assume that the monomials of $S$ are lexicographically ordered such that $x_1>x_2>\ldots>x_n$. Let $d\geq2$ be an integer and denote by $\mathcal{M}_d$ the set of all the monomials of degree $d$ of $S$. For two monomials $u,v$ of $\mathcal{M}_d$, with $u\geq_{lex}v$, the set $$\mathcal L(u,v)=\{w\in\mathcal{M}_d\ |\ u\geq_{lex}w\geq_{lex}v\}$$ is called the lexsegment determined by $u$ and $v$. A lexsegment ideal in $S$ is a monomial ideal of $S$ which is generated by a lexsegment. 

We treat separately the cases of a completely lexsegment ideal and of a non-completely lexsegment ideal in order to show that all lexsegment ideals with a linear resolution have linear quotients.

For the case of completely lexsegment ideals, we define a total order on the set of monomials in $\mathcal{M}_d$, denoted $\prec$, as follows: for two monomials $w,w'$ of degree $d$ in $S$, $w=x_1^{\alpha_1}\cdots x_n^{\alpha_n}$ and $w'=x_1^{\beta_1}\cdots x_n^{\beta_n},$ we set $w\prec w'$ if $\alpha_1<\beta_1$ or $\alpha_1=\beta_1$ and $w>_{lex}w'$.

Let $u,v$ be monomials in $\mathcal{M}_d$ such that $I=(\mathcal{L}(u,v))$ is a completely lexsegment ideal with a linear resolution. If $\mathcal{L}(u,v)=\{w_1,\ldots,w_r\}$ where $w_1\prec w_2\prec\ldots\prec w_r$, we show that $I$ has linear quotients with respect to this order of the generators. Next, we look to the decomposition function. We determine the decomposition function with respect to the ordering $\prec$ and we prove that it is regular. Therefore, one may apply the procedure developed in \cite{HeTa} to get the explicit resolutions for this class of ideals.   

Next, we consider the class of non-completely lexsegment ideals with a linear resolution. We prove that these ideals have also linear quotients, but with respect to a different order of the monomials. Let $u,\ v\in\mathcal{M}_d$ such that $I=(\mathcal{L}(u,v))$ is a non-completely lexsegment ideal with a linear resolution. We write $I$ as a sum of two ideals $J$ and $K$, where $J$ is generated by all the monomials from $\mathcal{L}(u,v)$ which are not divisible by $x_1$ and $K$ is generated by all the monomials from $\mathcal{L}(u,v)$ which are divisible by $x_1$. Let $G(J)=\{g_1,\ldots,g_m\}$, $g_1>_{lex}\ldots>_{lex} g_m$ and $G(K)=\{h_1,\ldots,h_r\}$,  $h_1>_{\overline{lex}}\ldots>_{\overline{lex}} h_r$ where we denote by $>_{\overline{lex}}$ the lexicographical order corresponding to $x_n>x_{n-1}>\ldots>x_1$. We prove that $I$ has linear quotients with respect to the sequence of monomials $g_1,\ldots,g_m,h_1,\ldots,h_r$.

Next, we compute the Krull dimension and the depth for arbitrary lexsegment ideals. We obtain that these invariants can be determined just looking at the monomials $u$ and $v$ that define the lexsegment. As a consequence, we characterize all the lexsegment ideals which are Cohen--Macaulay. The results of this chapter are contained in the joint paper \cite{EOS} with V. Ene and L. Sorrenti.

In Chapter $3$, we define a new class of monomial ideals, namely \textit{constructible ideals}\rm. This new concept has at the base the notion of constructible simplicial complex. Using the relation between Cohen--Macaulay simplicial complexes, shellable simplicial complexes and their Stanley--Reisner ideals of the Alexander dual, we have the following diagram:

\[\begin{array}{ccccc}
	\Delta\ \mbox{is shellable}&\Longrightarrow&\Delta\ \mbox{is constructible}&\Longrightarrow&\Delta\ \mbox{is Cohen--Macaulay}\\
	\Updownarrow&&&&\Updownarrow\\
	I_{\Delta^{\vee}}\ \mbox{is an ideal with}&&&&I_{\Delta^{\vee}}\ \mbox{is an ideal with a}\\
	\mbox{linear quotients}&&&&\mbox{linear resolution} 
\end{array}
\]

It naturally appears the problem of determining the Stanley--Reisner ideal of the Alexander dual associated to a constructible simplicial complex.

We define the notion of constructible ideal and we prove that a pure simplicial complex is constructible if and only if the Stanley--Reisner ideal of the Alexander dual is square-free constructible. Therefore, in the case of square-free monomial ideals generated in one degree, we get the following implications: 
\[\begin{array}{ccccc}
	
\mbox{ square-free monomial}&&\mbox{square-free }&&\mbox{ square-free monomial}\\
	\mbox{ideals with linear}&\Longrightarrow&\mbox{}&\Longrightarrow&\mbox{ideals with a linear }\\
	\mbox{quotients}&&\mbox{constructible ideals}&&\mbox{resolution} 
\end{array}.
\]

We show that the same implications hold for monomial ideals, not necessarily square-free. More precisely, we prove that all constructible ideals have linear resolutions. This allows us to obtain a recursive formula for the Betti numbers of a constructible ideal. Next, we pay attention to the behaviour of the "constructibility" in the process of polarization. We get that the property of being constructible is preserved during the polarization process. This helps us, in the last section of the chapter, to find an example of a constructible ideal which is not square-free and does not have linear quotients.

We prove next that all the monomial ideals with linear quotients generated in one degree are constructible ideals. The Betti numbers of monomial ideals with linear quotients generated in one degree are known \cite{HeTa}. We obtain the same Betti numbers using our recursive formula for constructible ideals.

In the end of the chapter, we analyse some examples. Starting with simplicial complexes, we get examples of square-free constructible ideals which do not have linear quotients and of square-free ideals with a linear resolution which are not constructible. We end with an example of a constructible ideal which is not square-free and does not have linear quotients. The results of this chapter are contained in our paper \cite{O1}.

In Chapter $4$ we study a particular class of simplicial complexes, that of the subword complexes in Coxeter groups. Subword complexes were defined by E. Miller and A. Knutson \cite{KM} for the study of Schubert polynomials and combinatorics of determinantal ideals. They proved that these simplicial complexes are vertex-decomposable \cite{KnMi}. Therefore, subword complexes in Coxeter groups are shellable.

We prove directly that subword complexes in Coxeter groups are shellable using the Alexander duality. More precisely, we show that the Stanley--Reisner ideal of the Alexander dual associated to a subword complex has linear quotients with respect to the lexicographical order of the minimal monomial generators. As a consequence, we get a shelling on the facets of the subword complex. For the Stanley--Reisner ideal of the Alexander dual, we get an upper bound for the projective dimension, and, in consequence, we obtain an upper bound for the regularity of the Stanley--Reisner ideal. 

In the last section of this chapter, we study a special class of subword complexes. Let $(W,S)$ be a Coxeter system, $Q$ a word of size $n$ and $\pi$ an element in $W$. Denote by $\Delta$ the subword complex, $G(I_{\Delta^{\vee}})=\{u_1,\ldots,u_r\}$ the minimal monomial generating set of $I_{\Delta^{\vee}}$, with $u_1>_{lex}\ldots>_{lex}u_r$, $r\leq n-\ell(\pi)+1$ and $d_r=r-1$, where $d_r=|\set(u_r)|$. For this class, we get that the Stanley--Reisner ideal of the Alexander dual is isomorphic with a monomial prime ideal, $J$. Therefore, the Koszul complex associated to the sequence $G(J)$ is isomorphic to a minimal graded free resolution for the Stanley--Reisner ring of the Alexander dual. Also, the Stanley--Reisner ring of $\Delta$ is a complete intersection ring. 

In the end, we determine all the subword complexes in this special class which are simplicial spheres or Cohen--Macaulay. The results in Chapter $4$ are contained in our paper \cite{O2}.

We acknowledge the support provided by the Computer Algebra Systems \textsc{CoCoA} \cite{Co} and \textsc{Singular} \cite{GPS} for the extensive experiments which helped us to obtain some of the results in this thesis.

\mainmatter
\pagestyle{fancy}
\chapter{Background}
In this chapter we recall some basic notions and some results that will be used later.
\section{Regular sequences and depth}
The notion of regular sequence appeared first in the papers of M. Auslander and D.A. Buchsbaum \cite{AB}, and D. Rees \cite{Re}. Using regular sequences, one can define the depth of a module. For details, see W. Bruns and J. Herzog \cite{BH}, R.H. Villareal \cite{Vi}.

Through this section, $R$ is a Noetherian ring and we assume that all the modules are finitely generated.

\begin{Definition1}\rm Let $M$ be a module over $R$. An element $x\in R$ is an \textit{$M$--regular element} if $xz=0$ for $z\in M$ implies $z=0$.
\end{Definition1}

\begin{Examples1}\rm 1. Let $R=k[x]$, where $k$ is a field and $M=R$. Then $x$ is regular on $R$.

2. Let $k$ be a field, $R=k[x,y,z]$ be the polynomial ring over the field $k$ and $M=R/(xy)$. Let $z^2\in R$. Then $z^2$ is an $M$-regular element. Indeed, let $m\in R$ be a monomial such that $mz^2\in(xy)$, that is $xy\mid mz^2$. Since $\gcd(xy,z^2)=1$, we must have $xy\mid m$.
\end{Examples1}
Using the above definition, one may consider the concept of a regular sequence on a module.

\begin{Definition1}\rm\label{reg} Let $M$ be a module over $R$. A sequence $\mathbf{x}=x_1,\ldots,x_r$ is a \textit{regular sequence on $M$} (or simply an \textit{$M$-sequence}) if the following conditions are satisfied:
\begin{itemize}
	\item[(i)] $x_i$ is an $M/(x_1,\ldots,x_{i-1})M$-regular element for $i=1,\ldots,r$
	\item[(ii)] $M/(\mathbf{x})M\neq0$. 
\end{itemize}
\end{Definition1}

 \begin{Remark1}\rm An $R$-regular sequence is called \textit{a regular sequence}.
 \end{Remark1}

\begin{Examples1}\rm 1. Let $R=k[x_1,\ldots,x_n]$ be the polynomial ring over a field $k$. Then the sequence of variables $x_1,\ldots,x_n$ is a regular sequence.

2. Let $k$ be a field and $R=k[x,y,z]$ be the polynomial ring over the field $k$. Consider the sequence $\textit{\textbf{f}}=xy,z^2$. Then $\textit{\textbf{f}}$ is a regular sequence.
\end{Examples1}

\begin{Proposition1}[\cite{BH}] Let $M$ be an $R$-module and $\mathbf{x}$ an $M$-regular sequence. Then an exact sequence
	\[N_2\stackrel{\varphi_2}{\longrightarrow}N_1\stackrel{\varphi_1}{\longrightarrow}N_0\stackrel{\varphi_0}{\longrightarrow}M\longrightarrow0
\]
of $R$-modules induces an exact sequence
\[N_2/(\mathbf{x})N_2\longrightarrow N_1/(\mathbf{x})N_1\longrightarrow N_0/(\mathbf{x})N_0\longrightarrow M/(\mathbf{x})M\longrightarrow0.
\]
\end{Proposition1}
The fact that $R$ is a Noetherian local ring is essential in the following proposition.
\begin{Proposition1}[\cite{BH}] Let $R$ be a local Noetherian ring, $\frak{m}$ its maximal ideal, $M$ be $R-$module, and $\mathbf{x}=x_1,\ldots,x_r$ an $M$-sequence. Then every permutation of $\mathbf{x}$ is an $M$-sequence.
\end{Proposition1}

 If $R$ is not a local ring, one may find examples of a regular sequence $\mathbf{x}$ such that a permutation of $\mathbf{x}$ is not a regular sequence.

\begin{Example1}\rm Let $R=k[x,y,x]$ be the polynomial ring over a field $k$. Then the sequence $\textit{\textbf{f}}=x,y(1-x),z(1-x)$ is a regular sequence, but $\textit{\textbf{f}}\ ^\prime=y(1-x),z(1-x),x$ is not a regular sequence. The fact that $\textit{\textbf{f}}$ is a regular sequence can be easily verified. For the second sequence, one may note that $z(1-x)$ is not regular on $R/(y(1-x))$. Indeed, $z(1-x)\cdot y\in(y(1-x))$ but $y\notin(y(1-x))$.  
\end{Example1}

Next we describe from the homological point of view the fact that an ideal $I$ contains regular elements.

\begin{Proposition1}[\cite{BH}] Let $M,N$ be $R$-modules and set $I=\Ann(N)$. Then $I$ contains an $M$-regular element if and only if $\Hom_R(N,M)=(0)$. 
\end{Proposition1}
 
 Let $M$ be an $R$-module and $\mathbf{x}=x_1,\ldots,x_n,\ldots$ be an $M$-sequence. Then we have the strictly ascending sequence of ideals:
	\[(x_1)\subset (x_1,x_2)\subset\ldots\subset(x_1,x_2,\ldots,x_n)\subset\ldots\ .
\]
Since $R$ is a Noetherian ring, this sequence must terminates. Thus, any $M$-sequence is finite.

\begin{Definition1}\rm Let $M$ be an $R$-module and $I$ be an ideal of $R$ such that $IM\neq M$. A \textit{maximal $M$-sequence} in $I$ is one which cannot be extended to a longer $M$-sequence in $I$. 
\end{Definition1}

\begin{Proposition1}[\cite{Vi}] Let $M$ be an $R$-module and $I$ be an ideal of $R$ such that $IM\neq M$. If $\mathbf{x}=x_1,\ldots,x_r$ is an $M$-sequence in $I$, then $\mathbf{x}$ can be extended to a maximal $M$-sequence in $I$.
\end{Proposition1}

The next result describes the relation between the length of a regular sequence on a module and the dimension of the module.

\begin{Lemma1}[\cite{Vi}]\label{rdim} Let $R$ be a Noetherian local ring, $\frak{m}$ be its maximal ideal, and $M$ be an $R$-module. If $\mathbf{x}=x_1,\ldots,x_r$ is an $M$-sequence in $\frak{m}$, then $r\leq\dim(M)$.
\end{Lemma1}

The following result is very important.
\begin{Theorem1}[Rees \cite{BH}]\label{Rees} Let $M$ be an $R$-module and $I$ an ideal of $R$ such that $IM\neq M$. Then all the maximal $M$-sequences in $I$ have the same length $n$ given by
	\[n=\min\{i:\Ext_R^i(R/I,M)\neq(0)\}.
\]
\end{Theorem1}

If $R$ is a Noetherian local ring, $\frak{m}$ is its maximal ideal, and $M$ is an $R$-module, it can be defined the concept of depth of the $R$-module $M$.

\begin{Definition1}\rm Let $R$ be a Noetherian local ring, $\frak{m}$ be its maximal ideal, and $M$ be an $R$-module. The \textit{depth of $M$}, denoted by $\depth(M)$, is the length of any maximal regular sequence on $M$ which is contained in $\frak{m}$.
\end{Definition1}

By Theorem \ref{Rees}, one has
\[\depth(M)=\min\{i:\Ext_R^i(R/\frak{m},M)\neq(0)\}.
\]
In general, by Lemma \ref{rdim}, we have $\depth(M)\leq\dim(M)$.
\section{Cohen--Macaulay ideals}
The origins of the theory of Cohen--Macaulay rings are the unmixedness theorems of F.S. Macaulay \cite{Ma} and I.S. Cohen \cite{co}. The foundations of the current shape of the theory were laid by M. Auslander, D.A. Buchasbaum, M. Nagata and D. Rees. A comprehensive monograph on Cohen--Macaulay rings is the book of W. Bruns and J. Herzog \cite{BH}. 

Through this section, $R$ is a Noetherian local ring and $\frak{m}$ its maximal ideal. We assume that all the modules are finitely generated.
 
We saw in the previous section that, for an $R$-module $M$, $\depth(M)\leq\dim(M)$.
 
\begin{Definition1}\rm An $R$-module $M$ is called \textit{Cohen--Macaulay} if $\depth(M)=\dim(M)$. A local ring $(R,\frak{m})$ is a \textit{Cohen--Macaulay ring} if $R$ is Cohen--Macaulay as an $R$-module.
\end{Definition1}

\begin{Example1}\rm Let $k$ be a field. Then $k[x_1,\ldots,x_n]$ and $k[[x_1,\ldots,x_n]]$ are Cohen--Macaulay rings.  
\end{Example1}

\begin{Definition1}\rm Let $I$ be an ideal of $R$. $I$ is a \textit{Cohen--Macaulay ideal} if $R/I$ is a Cohen--Macaulay $R$-module.
\end{Definition1}

\begin{Example1}\rm Let $R=k[x_1,\ldots,x_5]$ be the polynomial ring over a field $k$ and $$I=(x_1x_2, x_2x_3,x_3x_4,x_1x_4,x_3x_5,x_4x_5)$$ an ideal of $R$. One may note that $\height(I)=3$, hence $\dim(S/I)=2$ and $\depth(S/I)\leq2$. It can be verified that $x_1+x_2$ is regular on $R/I$ and $x_2+x_3+x_4+x_5$ is regular on $R/(I,x_1+x_2)$. Thus $\depth(R/I)=2$ and $R/I$ is a Cohen--Macaulay module.
\end{Example1}
The behaviour of the depth along an exact sequence is described in the next proposition.

\begin{Proposition1}[Depth lemma \cite{BH}] If
	\[0\longrightarrow N\longrightarrow M\longrightarrow L\longrightarrow 0
\]
is a short exact sequence of modules over a local ring $R$, then
\begin{itemize}
	\item[(a)] If $\depth(M)<\depth(L)$, then $\depth(N)=\depth(M)$;
	\item[(b)] If $\depth(M)=\depth(L)$, then $\depth(N)\geq\depth(M)$; 
	\item[(c)] If $\depth(M)>\depth(L)$, then $\depth(N)=\depth(L)+1$.
\end{itemize}
\end{Proposition1}

The dimension and the depth of a module $M$ are closely related to $M$-regular elements.

\begin{Proposition1}[\cite{BH}] If $M$ is a module over $R$ and $z$ is a regular element of $M$, then
\begin{itemize}
	\item[(a)] $\depth(M/(z)M)=\depth(M)-1$;
	\item[(b)] $\dim(M/(z)M)=\dim(M)-1$. 
\end{itemize}
\end{Proposition1}

We will use the following notion in the last chapter of this thesis.

\begin{Definition1}\rm Let $I$ be an ideal of $R$. If $I$ is generated by a regular sequence, $I$ is \textit{a complete intersection ideal} and $R/I$ is \textit{a complete intersection ring}.
\end{Definition1}

\begin{Proposition1}[\cite{BH}] A complete intersection ring is Cohen--Macaulay.
\end{Proposition1}

The following theorem is an instrument for computing the depth of a module.

\begin{Theorem1}[Auslander--Buchsbaum \cite{BH}] Let $(R,\frak{m})$ be a Noetherian local ring, and $M\neq 0$ a finite $R$-module. If $\projdim(M)<\infty$, then
	\[\projdim(M)+\depth(M)=\depth(R).
\]
\end{Theorem1}
\section{Classes of monomial ideals}
\subsection{Monomial ideals with linear quotients}
In \cite{HeTa}, J. Herzog and Y. Takayama defined the notion of ideal with linear quotients for describing the comparison maps in the construction of a free resolution as iterated mapping cones.

In this section we consider only monomial ideals with linear quotients.

Henceforth, we denote by $S$ the polynomial ring in $n$ variables over a field $k$, that is $S=k[x_1,\ldots,x_n]$. For a monomial ideal $I$ of $S$, we denote by $G(I)$ the minimal monomial generating set.

\begin{Definition1}\rm Let $I$ be a monomial ideal of $S$. $I$ \textit{has linear quotients} (or it is an \textit{ideal with linear quotients}) if there exists an ordering $u_1,\ldots,u_m$ of its minimal monomial generators such that, for all $1\leq i\leq m$, the colon ideals $(u_1,\ldots,u_{i-1}):u_i$ are generated by variables.
\end{Definition1}

The following result essentially uses the fact that $I$ is a monomial ideal. (See, for instance, J. Herzog \cite{HH}.)

\begin{Lemma1}[\cite{H}] The monomial ideal $I$ of $S$ has linear quotients with respect to the sequence of the minimal monomial generators $u_1,\ldots, u_m$ if and only if for all $i$ and for all $j<i$ there exist an integer $k<i$ and an integer $l$ such that
	\[\frac{u_k}{\gcd(u_k,u_i)}=x_l\ \mbox{and}\ x_l\ \mbox{divides}\ \frac{u_j}{\gcd(u_j,u_i)}.
\]
 \end{Lemma1}
 
 If $I$ is a square-free monomial ideal of $S$, we get the following consequence:

 \begin{Corollary1}[\cite{H}]\label{sqfreelinquot} Let $I$ be a square-free monomial ideal with $G(I)=\{u_1,\ldots,$ $u_m\}$ and let $F_i=\supp(u_i)$, for all $i=1,\ldots,m$, where we denote $\supp(u_i)=\{j\in[n]\ :\ x_j|u_i\}$. Then $I$ has linear quotients if and only if for all $i$ and for all $j<i$ there exist an integer $l\in F_j\setminus F_i$ and an integer $k<i$ such that $F_k\setminus F_i=\{l\}$.
 \end{Corollary1}
 
Let $d>0$ be an integer. Recall that an ideal $I$ of $S$ \textit{has a $d$--linear resolution} if the minimal graded free resolution of $I$ is of the form
	\[\ldots\longrightarrow S(-d-2)^{\beta_2}\longrightarrow S(-d-1)^{\beta_1}\longrightarrow S(-d)^{\beta_0}\longrightarrow I\longrightarrow 0.
\]

Equivalently, $I$ has a $d$--linear resolution if and only if $\Tor_i^S(I;k)_{i+j}=0$ for all $j\neq d$.

\begin{Proposition1}[\cite{H}]\label{linquotlinres} Let $I$ be a monomial ideal of $S$ generated in degree $d$ and assume that $I$ has linear quotients. Then $I$ has a $d$--linear resolution.
\end{Proposition1}

As a consequence, for a monomial ideal with linear quotients generated in one degree, we can compute the Betti numbers. Firstly, we fix some notations. Let $I$ be a monomial ideal of $S$ with $G(I)=\{u_1,\ldots,u_m\}$ and assume that $I$ has linear quotients with respect to the sequence $u_1,\ldots,u_m$. We denote $L_k=(u_1,\ldots,u_{k-1}):u_k$ and $r_k=|G(L_k)|$.

\begin{Corollary1}[\cite{HeTa}]\label{Bettilinquot} Let $I\subset S$ be a monomial ideal with linear quotients generated in one degree. Then, with the above notations, one has
\[\beta_i(I)=\sum_{k=1}^{m}\left(\twoline{r_k}{i} \right),
\]
for all $i$. In particular, it follows that $$\projdim(I)=\max\{r_1,\ldots,r_m\}.$$ 
\end{Corollary1}

\begin{Example1}\rm Let $S=k[x_1,\ldots,x_4]$ and $I\subset S$ be the monomial ideal with $G(I)=\{x_1x_2,\ x_1x_3^2x_4,\ x_2x_3\}$. One may easily check that $I$ has linear quotients with respect to this order of the minimal monomial generators.
\end{Example1}

On the other hand, we note that the ideal $J=(x_1^3,\ x_1x_2^2)\subset k[x_1,x_2]$ does not have linear quotients since $(x_1^3):( x_1x_2^2)=(x_1^2)$ and $(x_1x_2^2):(x_1^3)=(x_2^2)$.

Let $I$ be a monomial ideal of $S$ with $G(I)=\{u_1,\ldots,u_m\}$ and assume that $I$ has linear quotients with respect to the sequence $u_1,\ldots,u_m$. One usually denotes
	\[\set(u_i)=\{j\in[n]\ :\ x_j\in(u_1,\ldots,u_{i-1}):u_i\},
\]
for all $i$.

\begin{Example1}\label{ideal}\rm Let $S=k[x_1,\ldots,x_4]$ and $I\subset S$ be the monomial ideal with $G(I)=\{x_1^2x_2,\ x_1x_2^2,\ x_1x_2x_3,\ x_2x_3^2\}$. Let us denote $u_1=x_1^2x_2,\ u_2=x_1x_2^2,\ u_3=x_1x_2x_3,\ u_4=x_2x_3^2$. One may easily check that $I$ has linear quotients with respect to this order of the minimal monomial generators. We have $\set(u_1)=\emptyset,\ \set(u_2)=\{1\},\ \set(u_3)=\{1,\ 2\},$ and $ \set(u_4)=\{1\}$.
\end{Example1}

For a monomial ideal with linear quotients, J. Herzog and Y. Takayama \cite{HeTa} defined the decomposition function. Let us assume that $I$ is a monomial ideal of $S$ with $G(I)=\{u_1,\ldots,u_m\}$ and $I$ has linear quotients with respect to the sequence $u_1,\ldots,u_m$. Denote $I_j=(u_1,\ldots,u_j)$ for all $1\leq j\leq m$ and by $M(I)$ the set of all the monomials in $I$.

\begin{Definition1}\rm The map $g:M(I)\rightarrow G(I)$ defined by $g(u)=u_j$ if $j$ is the smallest number such that $u\in I_j$ is called the \textit{decomposition function} of $I$ with respect to the sequence $u_1,\ldots,u_m$.
\end{Definition1}

\begin{Example1}\rm Let $S=k[x_1,\ldots,x_4]$ and $I\subset S$ be the monomial ideal from Example \ref{ideal}. Let $g:M(I)\rightarrow G(I)$ be the decomposition function and $u=x_1^2x_2^2x_3^2$ be a monomial in $I$. Then $g(x_1^2x_2^2x_3^2)=x_1^2x_2$.
\end{Example1}

\begin{Definition1}\rm Let $I$ be a monomial ideal of $S$ with linear quotients. The decomposition function $g:M(I)\rightarrow G(I)$ is called \textit{regular} if $\set(g(x_su))\subseteq\set(u)$ for all $s\in\set(u)$ and $u\in G(I)$.
\end{Definition1}
	\[
\]
\subsection{Stable ideals}
S. Eliahou and M. Kervaire \cite{EK} introduced and studied the stable ideals. In \cite{AH}, A. Aramova and J. Herzog computed the resolution of a stable ideal using Koszul homology. Here we recall some notations and results given in S. Eliahou and M. Kervaire \cite{EK}.

Firstly, let us fix some notations. Let $S=k[x_1,\ldots,x_n]$ be the polynomial ring in $n$ variables over a field $k$. For a monomial $u=x_1^{\alpha_1}\cdots x_n^{\alpha_n}$ we denote by
	\[\max(u)=\max\{j\in[n]\ :\ x_j\mid u\}=\max \{j\in[n]\ :\ \alpha_j>0\}
\]
and by
\[\min(u)=\min\{j\in[n]\ :\ x_j\mid u\}=\min \{j\in[n]\ :\ \alpha_j>0\}.
\] 

\begin{Definition1}[\cite{EK}]\rm A monomial ideal $I$ of $S$ is called \textit{stable} if for every monomial $w\in I$ and for all positive integers $i$ with $i<\max(w)$ it holds
	\[\frac{x_iw}{x_{\max(w)}}\in I.
\]
\end{Definition1}
As usual, we denote by $G(I)$ the minimal monomial generating set of the monomial ideal $I$ of $S$.

The following lemma shows that one may check the stability of a monomial ideal using a finite number of tests.
 
\begin{Lemma1}[\cite{EK}]\label{stable} Let $I$ be a monomial ideal of $S$. The following conditions are equivalent:
\begin{itemize}
	\item[(i)] $I$ is stable.
	\item[(ii)] If $u\in G(I)$ and $i<\max(u)$, then $x_iu/x_{\max(u)}\in I$. 
\end{itemize}
\end{Lemma1}

\begin{Example1}\rm Let $S=k[x_1,x_2,x_3]$ and $I$ be the monomial ideal with the minimal monomial generating system $G(I)=\{x_1^2,\ x_1x_2^2,\ x_1x_2x_3\}$. One may easily check, using Lemma \ref{stable}, that $I$ is a stable ideal.
\end{Example1}

For a monomial $u=x_1^{\alpha_1}\cdots x_n^{\alpha_n}$ in $S$, we denote by $\nu_i(u)$ the exponent of the variable $x_i$ in the monomial $u$, that is $\nu_i(u)=\alpha_i$. Recall that, if $u,v$ are monomials in $S$, $u>_{lex}v$ means that $\deg(u)>\deg(v)$ or $\deg(u)=\deg(v)$ and there exists an integer $l$ such that, for all $i<l$, we have that $\nu_i(u)=\nu_i(v)$ and $\nu_l(u)>\nu_l(v)$.

\begin{Proposition1}[\cite{S}]\label{stablelq} Let $I$ be a stable monomial ideal of $S$ with $G(I)=\{u_1,\ldots,u_m\}$ where $u_1>_{lex}\ldots>_{lex}u_m$ with $x_1>x_2>\ldots>x_n$. Then $I$ has linear quotients with respect to this order of the minimal monomial generators.
\end{Proposition1}

An imediat consequence is the following:

\begin{Corollary1} If $I$ is a stable ideal of $S$ generated in one degree, then $I$ has a linear resolution. 
\end{Corollary1}

\begin{Lemma1}[\cite{EK}] Let $I\subset S$ be a stable monomial ideal with the minimal monomial generating system $G(I)$. For every monomial $w\in I$, there exists a unique decomposition $w=uy$ with $u\in G(I)$ and $\max(u)\leq\min(y)$. For convenience, we define $\min(1)=\infty$.
\end{Lemma1}

The unique decomposition of a monomial $w\in I$, $w=uy$, with $u\in G(I)$ and $\max(u)\leq\min(y)$, is called the \textit{canonical decomposition} of $w$.

\begin{Example1}\rm Let $S=k[x_1,x_2,x_3]$ and $I$ be the monomial ideal with the minimal monomial generating system $G(I)=\{x_1^2,\ x_1x_2^2,\ x_1x_2x_3\}$. Let $u=x_1^3x_3x_4$ be a monomial in $I$. Then $u=x_1^2\cdot x_1x_3x_4$, where $x_1^2\in G(I)$ and $\max(x_1^2)=\min(x_1x_3x_4)$.
\end{Example1}

Let us denote by $M(I)$ the set of all the monomials in $I$, where $I$ is a stable ideal. The map
	\[g:M(I)\rightarrow G(I)
\]
defined by $g(w)=u$ if $w=uy$ is the canonical decomposition of $w$, is called the \textit{decomposition function} of the stable ideal $I$.

One may note that the decomposition function for stable ideals is a particular case of the decomposition function defined for monomial ideals with linear quotients with respect to the lexicographical order of the minimal monomial generators.

	\[
\]

\subsection{Lexsegment ideals}\label{Lexsegment ideals}

Let $S=k[x_1,\ldots,x_n]$ be the polynomial ring in $n$ variables over a field $k$ and let us assume that all the monomials of $S$ are ordered by the lexicographical order with $x_1>x_2>\ldots>x_n$. Let $d>0$ be an integer and denote by $\mathcal{M}_d$ the set of all the monomials of $S$ of degree $d$. Let $u$ be a monomial in $\mathcal{M}_d$. The set
	\[\mathcal{L}^i(u)=\{w\in\mathcal{M}_d\ :\ w\geq_{lex} u\}
\]
is called an \textit{initial lexsegment} and the set
\[\mathcal{L}^f(u)=\{w\in\mathcal{M}_d\ :\ u\geq_{lex} w\}
\]
is a \textit{final lexsegment}. An \textit{initial lexsegment ideal} is a monomial ideal generated by an initial lexsegment. Similarly, a \textit{final lexsegment ideal} is a monomial ideal generated by a final lexsegment. Initial and final lexsegment ideals have been well-studied, see \cite{Ma}, \cite{Hul}, \cite{B}, \cite{De}.

Let $u,v\in\mathcal{M}_d$ and $u\geq_{lex}v$. The lexsegment defined by the monomials $u$ and $v$ is the set
\[\mathcal{L}(u,v)=\{w\in\mathcal{M}_d\ :\ u\geq_{lex} w\geq_{lex} v\}.
\]
A monomial ideal generated by a lexsegment is called a \textit{lexsegment ideal}. 

\begin{Example1}\rm Let $S=k[x_1,x_2,x_3,x_4]$. The set of all the monomials of degree $3$ is
	\[\mathcal {M}_3=\{x_1^3,\ x_1^2x_2,\ x_1^2x_3,\ x_1^2x_4,\ x_1x_2^2,\ x_1x_2x_3,\ x_1x_2x_4,\ x_1x_3^2,\ x_1x_3x_4,\ x_1x_4^2,\ x_2^3,\ x_2^2x_3,
\]
	\[x_2^2x_4, \ x_2x_3^2,\ x_2x_3x_4,\ x_2x_4^2,\ x_3^3,\ x_3^2x_4,\ x_3x_4^2,\ x_4^3\}\ \ \ \ \ \ \ \ \ \ \ \ \ \ \ \ \ \ \ \ \ \ \ \ \ \ \ 
\]
Let $u=x_1x_2x_3$ and $v=x_2x_3^2$, $u>_{lex} v,$ and let $I=(\mathcal{L}(u,v))$. Then
	\[G(I)=\mathcal{L}(u,v)=\{x_1x_2x_3,\ x_1x_2x_4,\ x_1x_3^2,\ x_1x_3x_4,\ x_1x_4^2,\ x_2^3,\ x_2^2x_3,\ x_2^2x_4, \ x_2x_3^2\}.\]
	
\end{Example1} 
Arbitrary lexsegment ideals have been defined by H. Hulett and H.M. Martin \cite{HM}. A. Aramova, E. De Negri and J. Herzog \cite{ADH} also studied these ideals and they characterized all the lexsegment ideals with a linear resolution.
One may note that initial lexsegment ideals are stable in the sense of Eliahou and Kervaire. Final lexsegment ideals are also stable in the sense of Eliahou and Kervaire, but with respect to the order of the variables $x_n>x_{n-1}>\ldots>x_1$. Thus, initial and final lexsegment ideals have linear quotients with respect to the lexicographical order, Proposition \ref{stablelq}.  

A useful tool in the study of lexsegment ideals is the shadow of a lexsegment. Recall that if $\mathcal{L}$ is a lexsegment, the \textit{shadow of $\mathcal{L}$} is the set
	\[\shad(\mathcal{L})=\{x_iw\ :\ w\in\mathcal{L},\ i=1,\ldots, n\}.
\]The $i$\textit{-th shadow} is recursively defined as $\shad^i(\mathcal{L})=\shad(\shad^{i-1}(\mathcal{L}))$. Any initial lexsegment has the property that its shadow is again an initial lexsegment, a fact which is not true for arbitrary lexsegments.

\begin{Example1}\rm Let $u=x_1x_3^2$, $v=x_2^3$ be monomials in $k[x_1,x_2,x_3]$. The lexsegment defined by the monomials $u$ and $v$ is $\mathcal{L}(u,v)=\{x_1x_3^2,\ x_2^3\}$. The set $\shad{\mathcal{L}(u,v)}$ is 
	\[\shad(\mathcal{L}(u,v))=\{x_1^2x_3^2,\ x_1x_2^3,\ x_1x_2x_3^2,\ x_1x_3^3,\ x_2^4,\ x_2^3x_3\}.
\]
One may easily check that $\shad{\mathcal{L}(u,v)}$ is not a lexsegment.
\end{Example1}

\begin{Definition1}[\cite{HM}]\rm$\ $ A lexsegment $\mathcal{L}$ is called a \textit{completely lexsegment} if all the iterated shadows of $\mathcal{L}$ are lexsegments. 
\end{Definition1}
 
A monomial ideal generated by a completely lexsegment is called a \textit{completely lexsegment ideal}. One may note that initial lexsegments are completely lexsegments.

The following theorem gives us a method to check whether a lexsegment ideal is a completely lexsegment ideal.

\begin{Theorem1}[Persistence \cite{DH}] Let $I$ be a lexsegment ideal generated in degree $d$ and suppose that $I_{d+1}$ is a lexsegment. Then $I$ is a completely lexsegment ideal. 
\end{Theorem1}

The following result gives a characterization of the final completely lexsegment ideals.

\begin{Proposition1}[\cite{DH}] Let $u$ be a monomial in $\mathcal{M}_d$ and let $K\subset S$ be the final lexsegment ideal generated by $\mathcal{L}^f(u)$. Then the following statements are equivalent:
\begin{itemize}
	\item[(a)] $K$ is a completely lexsegment ideal. 
	\item[(b)] $K_{d+1}$ is a lexsegment.
	\item[(c)] $u\geq_{lex}x_2^d$.  
\end{itemize}
\end{Proposition1}

The next result characterizes completely lexsegment ideals which are not final lexsegments.

\begin{Theorem1}[\cite{DH}]\label{condcompletely} Let $u=x_1^{a_1}\cdots x_n^{a_n}$ and $v=x_1^{b_1}\cdots x_n^{b_n}$ be monomials of degree $d$ in $S$, $v\neq x_n^d$, $u\geq_{lex}v$ and let $I$ be the ideal generated by $\mathcal{L}(u,v)$. Then $I$ is a completely lexsegment ideal if and only if $a_1\neq0$ and one of the following conditions holds:
\begin{itemize}
	\item[(a)] $u=x_1^px_2^{d-p}$ and $v=x_1^px_{n}^{d-p}$ for some integer $0<p\leq d$.
	\item[(b)] $a_1\neq b_1$ and for every $w<v$ there exists an integer $i>1$ such that $x_i|w$ and $x_1w/x_i\leq_{lex}u$. 
\end{itemize}
\end{Theorem1}

\begin{Example1}\rm Let $u=x_1x_3^2$, $v=x_2x_4^2$ be monomials in $k[x_1,\ldots, x_4]$. Then $I=(\mathcal{L}(u,v))\subset k[x_1,\ldots, x_4]$ is a completely lexsegment ideal.
\end{Example1}

\begin{Example1}\rm Let $u=x_1x_3x_4$, $v=x_2x_4^2$ be monomials in $k[x_1,\ldots, x_4]$. One may note that, for the monomial $x_3^3<_{lex}x_2x_4^2$, the condition (b) in Theorem \ref{condcompletely} does not hold. Thus, the ideal $I=(\mathcal{L}(u,v))\subset k[x_1,\ldots, x_4]$ is not a completely lexsegment ideal.
\end{Example1}
A. Aramova, E. De Negri and J. Herzog \cite{ADH} classified all lexsegment ideals with a linear resolution.

\begin{Theorem1}[\cite{ADH}]\label{completelylex} Let $u=x_1^{a_1}\cdots x_n^{a_n}$ and $v=x_1^{b_1}\cdots x_n^{b_n}$ be monomials of degree $d$ in $S$, $u\geq_{lex}v$ such that $I=(\mathcal{L}(u,v))$ is a completely lexsegment ideal. Then $I$ has a linear resolution if and only if one of the following conditions holds:
\begin{itemize}
	\item[(a)]  $u=x_1^px_2^{d-p}$ and $v=x_1^px_{n}^{d-p}$ for some integer $0<p\leq d$.
	\item[(b)] $b_1<a_1-1$.
	\item[(c)] $b_1=a_1-1$ and for the greatest $w<v$, $w\in\mathcal{M}_d$, one has $x_1w/x_{\max(w)}\leq_{lex} u$.
\end{itemize}
\end{Theorem1}
If $I$ is a non-completely lexsegment ideal, we have the following result:
\begin{Theorem1}[\cite{ADH}]\label{noncompletelylex} Let $u=x_1^{a_1}\cdots x_n^{a_n}$ and $v=x_1^{b_1}\cdots x_n^{b_n}$ be monomials of degree $d$ in $S$, $a_1\neq0$, $u\geq_{lex}v$. Suppose that $I=(\mathcal{L}(u,v))$ is not a completely lexsegment ideal. Then $I$ has a linear resolution if and only if $u$ and $v$ are of the form
	\[u=x_1x_{l+1}^{a_{l+1}}\cdots x_n^{a_n},\ v=x_lx_n^{d-1},
\]
for some $l$, $2\leq l<n$.
\end{Theorem1}
\section{The Koszul complex}

We will recall the Koszul complex following W. Bruns and J. Herzog \cite[Chapter 1, Section 1.6]{BH}.
  
Let $R$ be a commutative ring with unity and $\textit{\textbf{f}}=f_1,\ldots,f_r$ be a sequence of elements of $R$. Let $F$ be a free $R$--module with basis $\{e_1,\ldots,e_r\}$ and $K_j(\textit{\textbf{f}}\,;R)=\bigwedge\limits^j F$ for $0\leq j\leq r$. Note that $K_j(\textit{\textbf{f}}\,;R)$ is a free $R$-module with basis $$\{e_{{\sigma}}=e_{i_1}\wedge e_{i_2}\wedge\ldots \wedge e_{i_j}\ :\ \sigma=(i_1,\ldots,i_j),\ 1\leq i_1<i_2<\ldots<i_j\leq n\}.$$ In particular, $\rank(K_j(\textit{\textbf{f}}\,;R))=\left(\twoline{r}{j}\right)$.

The \textit{Koszul complex associated to the sequence $\textbf{f}$} is denoted by $K_{\bullet}(\textit{\textbf{f}}\,;R)$ and it is defined as:
	\[K_{\bullet}(\textit{\textbf{f}}\,;R):\ 0\longrightarrow K_r(\textit{\textbf{f}}\,;R)\stackrel{\partial_r}{\longrightarrow}K_{r-1}(\textit{\textbf{f}}\,;R)\longrightarrow\ldots\longrightarrow  K_1(\textit{\textbf{f}}\,;R)\stackrel{\partial_1}{\longrightarrow}K_0(\textit{\textbf{f}}\,;R)\longrightarrow0
\]
where the differential $\partial_j:K_j(\textit{\textbf{f}}\,;R)\longrightarrow K_{j-1}(\textit{\textbf{f}}\,;R)$ is given by 
	\[\partial_j(e_{\sigma})=\sum_{k\in\sigma}(-1)^{\alpha(\sigma,k)} f_{k} e_{\sigma\setminus k},
\]
with $\alpha(\sigma,k)=|\{i\in\sigma\ :\ i<k\}|$, $\sigma\subset\{1,\ldots,r\}$, $|\sigma|=j$.
One may easily check that $\partial_j\circ\partial_{j+1}=0$ for all $j$, so $K_{\bullet}(\textit{\textbf{f}}\,;R)$ is indeed a complex.

Let us denote $\textit{\textbf{f}}\,^{\prime}=f_1,\ldots,f_{r-1}$. Then, the Koszul complex associated to the sequence $\textit{\textbf{f}}=f_1,\ldots,f_r$ can be obtained as follows:
	\[K_{\bullet}(\textit{\textbf{f}}\,;R)=K_{\bullet}(\textit{\textbf{f}}\,^{\prime}\,;R)\otimes K_{\bullet}(f_r;R).
\]

\begin{Example1}\rm Let $S=k[x_1,x_2,x_3]$ and $\textit{\textbf{f}}=x_1^2,\ x_1x_2x_3,\ x_3^3$ be a sequence of monomials in $S$. The Koszul complex associated to the sequence $\textit{\textbf{f}}$ is:
	\[\  0\longrightarrow S(-8)\stackrel{\partial_3}{\longrightarrow} S(-5)^2\oplus S(-6)\stackrel{\partial_2}{\longrightarrow} S(-2)\oplus S(-3)^2\stackrel{\partial_1}{\longrightarrow} S\longrightarrow0
\]
where the differentials are:

$\partial_1(e_i)=f_i$ for all $i\in\{1,2,3\}$, that is
\[\partial_1=
\left(\begin{array}{ccc}
	x_1^2&x_1x_2x_3&x_3^3
\end{array}\right).
\]

$\partial_2(e_1\wedge e_2)=x_1^2\,e_2-x_1x_2x_3\,e_1$,

$\partial_2(e_1\wedge e_3)=x_1^2\,e_3-x_3^3\,e_1$,

$\partial_2(e_2\wedge e_3)=x_1x_2x_3\,e_3-x_3^3\,e_2$,
that is
 
 \[\partial_2=
\left(\begin{array}{ccc}
	-x_1x_2x_3&-x_3^3&0\\
	x_1^2&0&-x_3^3\\
	0&x_1^2&x_1x_2x_3
\end{array}\right).
\]

$\partial_3(e_1\wedge e_2\wedge e_3)=x_1^2\,e_2\wedge e_3-x_1x_2x_3\, e_1\wedge e_3+x_3^3\,e_1\wedge e_2$, that is

\[\partial_3=
\left(\begin{array}{c}
	x_3^3\\
	-x_1x_2x_3\\
	x_1^2
\end{array}\right).
\]
\end{Example1}
Let $M$ be an arbitrary $R$--module. The complex
	\[K_{\bullet}(\textit{\textbf{f}}\,;M)=K_{\bullet}(\textit{\textbf{f}}\,;R)\otimes M
\]
is the \textit{Koszul complex associated to the sequence $\textbf{f}$ with coefficients in $M$}. Therefore, $K_{\bullet}(\textit{\textbf{f}}\,;M)$ is the complex
\[K_{\bullet}(\textit{\textbf{f}}\,;M):\ 0\rightarrow K_r(\textit{\textbf{f}}\,;R)\otimes M\rightarrow \ldots\rightarrow  K_1(\textit{\textbf{f}}\,;R)\otimes M\rightarrow K_0(\textit{\textbf{f}}\,;R)\otimes M\rightarrow0
\]
Denote $K_j(\textit{\textbf{f}}\,;M):=K_j(\textit{\textbf{f}}\,;R)\otimes M$, for all $0\leq j\leq r$. Then, the differentials are $\partial_j\otimes \id: K_j(\textit{\textbf{f}}\,;M)\rightarrow K_{j-1}(\textit{\textbf{f}}\,;M)$, defined by
	\[(\partial_j\otimes\id)(e_{\sigma}\otimes m)=\sum_{k\in\sigma}(-1)^{\alpha(\sigma,k)} f_{k} e_{\sigma\setminus k}\otimes m,
\]
for all $\sigma\subset\{1,\ldots,r\}$ with $|\sigma|=j$ and for all $m\in M$.

Denote by $H_i(\textit{\textbf{f}}\,;R)$ the homology modules of the Koszul complex $K_{\bullet}(\textit{\textbf{f}}\,;R)$ and by $H_i(\textit{\textbf{f}}\,;M)$ the homology modules of $K_{\bullet}(\textit{\textbf{f}}\,;M)$.

\begin{Proposition1}[\cite{BH}] Let $I$ be the ideal of $R$ generated by the sequence $\textit{\textbf{f}}=f_1,\ldots,f_r$ of elements of $R$. Then the following statements hold:
\begin{itemize}
	\item[(a)] $H_0(\textit{\textbf{f}}\,;R)=R/I$, $H_0(\textit{\textbf{f}}\,;M)=M/IM$.
	\item[(b)] $H_r(\textit{\textbf{f}}\,;R)=\Ann_R(I)$, $H_r(\textit{\textbf{f}}\,;M)=\Ann_M(I)$, where $$\Ann_M(I)=\{m\in M\ :\ f_i\, m=0,\ 1\leq i\leq r\}.$$
	\item[(c)] The Koszul complex is an exact functor i.e. if 
	\[0\longrightarrow U\longrightarrow M\longrightarrow N\longrightarrow0
\]
is an exact sequence of $R$--modules, then
\[
0\longrightarrow K_{\bullet}(\textit{\textbf{f}}\,;U)\longrightarrow K_{\bullet}(\textit{\textbf{f}}\,;M)\longrightarrow K_{\bullet}(\textit{\textbf{f}}\,;N)\longrightarrow0
\]
is an exact sequence of complexes. In particular, there exists a long exact sequence
	\[\ldots\longrightarrow H_i(\textit{\textbf{f}}\,;U)\longrightarrow H_i(\textit{\textbf{f}}\,;M)\longrightarrow H_i(\textit{\textbf{f}}\,;N)\longrightarrow H_{i-1}(\textit{\textbf{f}}\,;U)\longrightarrow\ldots
\]
of homology modules.
\end{itemize}
\end{Proposition1}
The next result describes the relation between regular sequences and the Koszul complex. 
\begin{Theorem1}[\cite{BH}]\label{Koszulres} Let $\textit{\textbf{f}}=f_1,\ldots,f_m$ be a sequence of elements of $R$ and $M$ an $R$--module. 
\begin{itemize}
	\item[(a)] If $\textit{\textbf{f}}$ is an $M$--sequence, then $K_{\bullet}(\textit{\textbf{f}}\,;M)$ is acyclic.
	\item[(b)] If $\textit{\textbf{f}}$ is an $R$--sequence, then $K_{\bullet}(\textit{\textbf{f}}\,;R)$ is a free resolution of $R/(\textit{\textbf{f}}\,)$.
\end{itemize}
\end{Theorem1}

 \section{The Eliahou-Kervaire resolution}
 
 We follow A. Aramova and J. Herzog \cite{AH} for the construction of the Eliahou--Kervaire resolution.
 
 Let $S=k[x_1,\ldots,x_n]$ be the polynomial ring and denote $\underline{\textbf{x}}=x_1,\ldots,x_n$ the sequence of variables. Let $M$ be a finitely generated graded $S$--module and $I$ a stable ideal of $S$. Henceforth, we denote by $H_i(\underline{\textbf{x}})$ the Koszul homology module $H_i(\underline{\textbf{x}}; S/I)$. Let $\varepsilon:S\rightarrow S/I$ be the canonical surjection. As usual, we denote $\frak{m}=(x_1,\ldots,x_n)$ the maximal ideal generated by the set of all the variables. 
 
For all $i\geq0$, there exist isomorphisms 
 
	\[H_i(\underline{\mathbf{x}})\cong\Tor_i^S(k,S/I).
\]
In particular, this implies that $\beta_i(S/I)=\dim_k H_i(\underline{\textbf{x}})$, for all $i\geq0$.
Hence, the free $S$--resolution of $S/I$ may be written in the form
	\[\ldots\longrightarrow S\otimes H_2(\underline{\mathbf{x}})\stackrel{\nu_2}{\longrightarrow}S\otimes H_1(\underline{\mathbf{x}})\stackrel{\nu_1}{\longrightarrow}S\otimes H_0(\underline{\mathbf{x}})\stackrel{\nu_0}{\longrightarrow}S/I\longrightarrow 0, 
\]
where the maps $\nu_i,\ i\geq0$ are defined below.
 
 \begin{Theorem1}[\cite{EK}, \cite{AH}]\label{base} For all $j=1,\ldots,n$ and $i>0$, the Koszul homology $H_i(x_j,\ldots,x_n)$ is annihilated by $\frak{m}$. A basis of $H_i(x_j,\ldots,x_n)$ is given by the homology classes of the cycles
	\[\varepsilon(u')e_{\sigma}\wedge e_{\max(u)},\ u\in G(I),\ |\sigma|=i-1,\ j\leq\min(\sigma),\ \max(\sigma)<\max(u)
\]
where $u'=u/x_{\max(u)}$.
 \end{Theorem1}
 
 An important consequence of the Theorem \ref{base} is the following corollary:

 \begin{Corollary1}[Eliahou-Kervaire \cite{EK}] Let $I\subset S$ be a stable ideal. Then:
\begin{itemize}
	\item[(a)] $\beta_i(I)=\sum\limits_{u\in G(I)}\left(\twoline{\max(u)-1}{i}\right)$, for all $i\geq 0$;
	\item[(b)] $\projdim_S(I)=\max\{\max(u)-1\ :\ u\in G(I)\}$;
	\item[(c)] $\reg(I)=\max\{\deg(u)\ :\ u\in G(I)\}$. 
\end{itemize}
\end{Corollary1}

\begin{Corollary1}[\cite{EK}] Let $I\subset S$ be a stable ideal. Then the Hilbert series of $I$ is
	\[H(I,t)=\sum_{u\in G(I)}\frac{t^{\deg(u)}}{(1-t)^{n-\max(u)+1}}.
\]
\end{Corollary1}

Let $\mathbb{F}_{\bullet}$ be the minimal free resolution of $S/I$ over $S$
	\[\mathbb{F}_{\bullet}:\ \ldots\rightarrow S\otimes H_2(\underline{\mathbf{x}})\stackrel{\nu_2}{\rightarrow}S\otimes H_1(\underline{\mathbf{x}})\stackrel{\nu_1}{\rightarrow}S\otimes H_0(\underline{\mathbf{x}})\stackrel{\varepsilon}{\rightarrow}S/I\rightarrow 0. 
\]
By Theorem \ref{base}, a basis in $G_i:=S\otimes H_i(\underline{\mathbf{x}})$, $i>0$ is
	\[\{f(\sigma;u):=1\otimes(-1)^{\frac{(i-1)(i-2)}{2}}[\varepsilon(u')e_{\sigma}\wedge e_{\max(u)}]:\ \sigma\subset\{1,\ldots,n\},\ |\sigma|=i-1,\]\[ u\in G(I),\ \max(\sigma)<\max(u)\},
\]
where $u'=u/x_{\max(u)}$.

Denote by $M(I)$ the set of all the monomials from the ideal $I$. Let $g:M(I)\rightarrow G(I)$ be the decomposition function. For all $j$, $1\leq j\leq n$, and for all $u\in G(I)$, we denote
	\[u_j=g(x_ju)\ \mbox{and}\ y_j=\frac{x_ju}{u_j}.
\]

\begin{Theorem1}[\cite{EK}, \cite{AH}]\label{map} The chain maps of the resolution $\mathbb{F}_{\bullet}$ are
	\[\nu_i(f(\sigma;u))=\sum_{t\in\sigma}(-1)^{\alpha(\sigma,t)}(-x_tf(\sigma\setminus t;u)+y_tf(\sigma\setminus t;u_t)),
\]
for all $i\geq2$ and $\nu_1(f(\emptyset;u))=u$. Set $f(\sigma;u)=0$ if $\max(\sigma)\geq\max(u)$.
\end{Theorem1}

\begin{Example1}\rm Let $I=(x_1^2,\ x_1x_2^2,\ x_1x_2x_3,\ x_2^3)$ be a monomial ideal in the polynomial ring $k[x_1,x_2,x_3]$. One may easily check that $I$ is stable. The minimal monomial system of generators is $G(I)=\{x_1^2,\ x_1x_2^2,\ x_1x_2x_3,\ x_2^3\}$. I has the minimal free resolution over $S$:
	\[0\longrightarrow G_3\stackrel{\nu_3}{\longrightarrow}G_2\stackrel{\nu_2}{\longrightarrow}G_1\stackrel{\nu_1}{\longrightarrow}I\longrightarrow0.
\]
By easy computations, using Theorem \ref{base}, one obtains that:

$G_1$ has the basis $\{f(\emptyset;u)\ :\ u\in G(I)\}$,

$G_2$ has the basis $\{f(1;x_1x_2^2),\ f(1;x_1x_2x_3),\ f(2;x_1x_2x_3),\ f(1;x_2^3)\}$, and

$G_3$ has the basis $\{f((1,2);x_1x_2x_3)\}$.

Using Theorem \ref{map}, one may obtain the maps of the resolution:

$\nu_1(f(\emptyset;u))=u$, for all $u\in G(I)$, that is 

	\[\nu_1=
\left(\begin{array}{cccc}
	x_1^2&x_1x_2^2&x_1x_2x_3&x_2^3
\end{array}\right).
\]

$\nu_2(f(1;x_1x_2^2))=-x_1f(\emptyset;x_1x_2^2)+x_2^2f(\emptyset;x_1^2)$,

$\nu_2(f(1;x_1x_2x_3))=-x_1f(\emptyset;x_1x_2x_3)+x_2x_3f(\emptyset;x_1^2)$,

$\nu_2(f(2;x_1x_2x_3))=-x_2f(\emptyset;x_1x_2x_3)+x_3f(\emptyset;x_1x_2^2)$,

$\nu_2(f(1;x_2^3))=-x_1f(\emptyset;x_2^3)+x_2f(\emptyset;x_1x_2^2)$,

that is
 
 \[\nu_2=
\left(\begin{array}{cccc}
	x_2^2&x_2x_3&0&0\\
	-x_1&0&x_3&x_2\\
	0&-x_1&-x_2&0\\
	0&0&0&-x_1
\end{array}\right).
\]

$\nu_3(f((1,2);x_1x_2x_3))=-x_1f(2;x_1x_2x_3)+x_2f(1;x_1x_2x_3)-x_3f(1;x_1x_2^2)$, that is

\[\nu_3=
\left(\begin{array}{c}
	-x_3\\
	x_2\\
	-x_1\\
	0
\end{array}\right).
\]
Thus, the minimal free resolution of $I$ is
\[0\longrightarrow S(-5)\stackrel{\nu_3}{\longrightarrow}S(-4)^4\stackrel{\nu_2}{\longrightarrow}S(-2)\oplus S(-3)^3\stackrel{\nu_1}{\longrightarrow}I\longrightarrow0.
\]
\end{Example1}

\section{Mapping cones}\label{mapping}

Many free resolutions arise as iterated mapping cones. G. Evans and H. Charalambous \cite{EC} proved that the Eliahou--Kervaire resolution of stable monomial ideals is one of them. For a brief description and some properties of the mapping cones, see D. Eisenbud \cite[pp. 650--655]{Ei}.

Let $\mathbb{F}_{\bullet}$ and $\mathbb{G}_{\bullet}$ be two complexes and $\alpha:\mathbb{F}_{\bullet}\rightarrow\mathbb{G}_{\bullet}$ be a complex homomorphism. We write $\varphi$ and $\psi$, respectively, for the chain maps of $\mathbb{F}_{\bullet}$ and $\mathbb{G}_{\bullet}$.

\begin{Definition1}\rm The \textit{mapping cone} $M(\alpha)$ of $\alpha$ is the complex such that $$M(\alpha)_i=G_i\oplus F_{i-1}$$ with the chain map $d$ given by $d_i:M(\alpha)_i\rightarrow M(\alpha)_{i-1}$
	\[d_i(g,f)=(\alpha_{i-1}(f)+\psi_i(g),-\varphi_{i-1}(f)).
\]
\end{Definition1}

We are going to use this construction in the particular case of monomial ideals with linear quotients. We follow J. Herzog and Y. Takayama \cite{HeTa}.

Let $I$ be a monomial ideal with $G(I)=\{u_1,\ldots,u_m\}$ and assume that $I$ has linear quotients with respect to the sequence $u_1,\ldots,u_m$. Set $I_j=(u_1,\ldots, u_j)$ and $L_j=(u_1,\ldots, u_{j-1}):u_j$. One has $I_{j+1}/I_j\cong S/L_{j+1}$. Therefore, we get the exact sequences
	\[0\longrightarrow S/L_{j+1}\stackrel{u_{j+1}}{\longrightarrow} S/I_j\longrightarrow S/I_{j+1}\longrightarrow 0.
\]
Let $F^{(j)}$ be a graded free resolution of $S/I_{j}$. Recall that the ideals $L_j$ are generated by a sequence of variables, for $j=2,\ldots,m$. Let $K^{(j)}$ be the Koszul complex associated to the regular sequence $x_{k_1},\ldots,x_{k_l}$, with $k_i\in\set(u_{j+1})$, $i=1,\ldots,l$ and $\psi^{(j)}:K^{(j)}\longrightarrow F^{(j)}$ be a graded complex homomorphism lifting $S/L_{j+1}\stackrel{u_{j+1}}{\longrightarrow} S/I_j$. Then the mapping cone $M(\psi^{(j)})$ of $\psi^{(j)}$ yields a free resolution of $S/I_{j+1}$. By iterated mapping cones, one may obtain step by step a graded free resolution of $S/I$.

\begin{Lemma1}[\cite{HeTa}] Suppose $\deg (u_1) \leq \deg( u_2) \leq \ldots\leq \deg (u_m).$ Then the iterated mapping
cone $F$, derived from the sequence $u_1,\ldots,u_m$ is a minimal graded free resolution
of $S/I$, and for all $i > 0$ the symbols
\[f(\sigma; u)\ \mbox{with}\ u\in G(I),\ \sigma \subset \set(u),\  |\sigma| = i - 1
\]
form a homogeneous basis of the $S-$module $F_i$. Moreover $\deg(f(\sigma; u)) = |\sigma| +\deg(u)$.
\end{Lemma1} 

The following result generalizes the theorem of Eliahou and Kervaire for stable ideals.

\begin{Theorem1}[\cite{HeTa}]\label{mapm} Let $I$ be a monomial ideal of $S$ with linear quotients and $\mathbb{F}_{\bullet}$ the graded
minimal free resolution of $S/I$. Suppose that the decomposition function $g : M(I) \rightarrow G(I)$ is regular. Then the chain map $\partial$ of $\mathbb{F}_{\bullet}$ is given by
\[\partial(f(\sigma; u)) = -\sum_{t\in\sigma}(-1)^{\alpha(\sigma;t)}x_tf(\sigma\setminus t;u)+\sum_{t\in\sigma}(-1)^{\alpha(\sigma;t)}\frac{x_tu}{g(x_tu)}f(\sigma\setminus t;g(x_tu)),\]
if $\sigma\neq\emptyset$, and
\[\partial(f(\emptyset; u)) = u\] otherwise.
Here $\alpha(\sigma;t)=|\{s\in\sigma\ |\ s<t\}|$ and set $f(\sigma;u)=0$ if $\sigma\nsubseteq\set(u)$.
\end{Theorem1} 

\begin{Example1}\rm Let $I=(x_1x_2,\ x_2x_3x_4,\ x_2x_3^2)$ be a monomial ideal in the polynomial ring $k[x_1,\ldots,x_4]$. Denote by $u_1=x_1x_2,\ u_2=x_2x_3x_4,\ u_3=x_2x_3^2$. One may easily check that $I$ has linear quotients with respect to this order of the minimal monomial generators. Moreover, the decomposition function for this order is regular. Since $\max\{|\set(u)|\ :\ u\in G(I)\}=2$ we have that $F_i=0$, for all $i\geq4$. Hence 
	\[0\longrightarrow F_3\stackrel{\partial_3}{\longrightarrow}F_2\stackrel{\partial_2}{\longrightarrow}F_1\stackrel{\partial_1}{\longrightarrow}S\longrightarrow S/I\longrightarrow0.
\]
It is easy to see that:

$F_1$ has the basis $\{f(\emptyset;u)\ :\ u\in G(I)\}$,

$F_2$ has the basis $\{f(1;x_2x_3x_4),\ f(1;x_2x_3^2),\ f(4;x_2x_3^2)\}$, and

$F_3$ has the basis $\{f((1,4);x_2x_3^2)\}$.

Using Theorem \ref{mapm}, one may obtain the maps of the resolution:

$\partial_1(f(\emptyset(u)))=u$, for all $u\in G(I)$, that is 

	\[\partial_1=
\left(\begin{array}{ccc}
	x_1x_2&x_2x_3x_4&x_2x_3^2
\end{array}\right).
\]

$\partial_2(f(1;x_2x_3x_4))=-x_1f(\emptyset;x_2x_3x_4)+x_3x_4f(\emptyset;x_1x_2)$,

$\partial_2(f(1;x_2x_3^2))=-x_1f(\emptyset;x_2x_3^2)+x_3^2f(\emptyset;x_1x_2)$,

$\partial_2(f(4;x_2x_3^2))=-x_4f(\emptyset;x_2x_3^2)+x_3f(\emptyset;x_2x_3x_4)$

that is
 
 \[\partial_2=
\left(\begin{array}{ccc}
	x_3x_4&x_3^2&0\\
	-x_1&0&x_3\\
	0&-x_1&-x_4
	\end{array}\right).
\]

$\partial_3(f((1,4);x_2x_3^2))=-x_1f(4;x_2x_3^2)+x_4f(1;x_2x_3^2)-x_3f(1;x_2x_3x_4)$.

Therefore

\[\partial_3=
\left(\begin{array}{c}
	-x_3\\
	x_4\\
	-x_1
\end{array}\right).
\]
Thus, the minimal free resolution of $S/I$ is
\[0\longrightarrow S(-5)\stackrel{\partial_3}{\longrightarrow}S(-4)^3\stackrel{\partial_2}{\longrightarrow}S(-2)\oplus S(-3)^2\stackrel{\partial_1}{\longrightarrow}S\longrightarrow S/I\longrightarrow0.
\]
\end{Example1}

\section{Simplicial complexes}

\subsection{Basic notions}

Let $n>0$ be an integer. We denote by $[n]$ the set $\{1,\ldots,n\}$.

\begin{Definition1}\rm A \textit{simplicial complex} $\Delta$ on $[n]$ is a collection of subsets of $[n]$ such that the following conditions hold:
\begin{enumerate}
	\item $\{i\}\in\Delta$ for all $i\in[n]$; 
	\item if $F\in\Delta$ and $G\subseteq F$, then $G\in\Delta$. 
\end{enumerate}
\end{Definition1}

The set $[n]$ is the \textit{vertex set} of $\Delta$. The elements of $\Delta$ are called \textit{faces} and the \textit{dimension} of a face $F\in\Delta$ is denoted by $\dim(F)$ and $\dim(F)=|F|-1$. The faces of dimension $0$ are called \textit{vertices} and the faces of dimension $1$ are called \textit{edges}.  A \textit{facet} is a maximal face (with respect to the inclusion). We denote by $\mathcal{F}(\Delta)$ the set of all the facets of $\Delta$.
 
Let $d=\max\{|F|\ : \ F\in\Delta\}$. Then, the \textit{dimension of $\Delta$} is $\dim(\Delta)=d-1$. 

\begin{Example1}\label{sc}\rm Let $\Delta$ be the following simplicial complex on the vertex set $\{1,2,3,4\}$: 
\begin{figure}[h]
\begin{center}
\unitlength 1mm 
\linethickness{0.4pt}
\ifx\plotpoint\undefined\newsavebox{\plotpoint}\fi 
\begin{picture}(75,40.75)(0,0)
\multiput(7.5,8.5)(.033737024,.053633218){578}{\line(0,1){.053633218}}
\put(7.5,8.25){\line(1,0){37}}
\multiput(44.5,8.25)(-.033681214,.059297913){527}{\line(0,1){.059297913}}
\put(44.5,8.5){\line(1,0){28.5}}
\put(26.5,40.75){$1$}
\put(4.25,6.75){$2$}
\put(46,11){$3$}
\put(75,6.75){$4$}
\multiput(17.25,12.25)(.03361345,.05777311){238}{\line(0,1){.05777311}}
\multiput(21.5,12.5)(.03365385,.06089744){156}{\line(0,1){.06089744}}
\multiput(26.25,12)(.03350515,.06443299){97}{\line(0,1){.06443299}}
\end{picture}
\caption{}\label{f1}
\end{center}
\end{figure}

The facets of $\Delta$ are $\mathcal{F}(\Delta)=\{\{1,2,3\}, \ \{3,4\}\}$.
\end{Example1}

A simplicial complex $\Delta$ with the facets $F_1,\ldots,F_r$, is often denoted by $\Delta=\langle F_1,\ldots,F_r\rangle$. A \textit{simplex} is a simplicial complex with only one facet.

\begin{Example1}\rm Let $\Delta$ be the simplicial complex with the vertex set $\{1,2,3\}$:
	\[\Delta=\{\emptyset,\ \{1\},\ \{2\},\ \{3\},\ \{1,2\},\ \{1,3\},\ \{2,3\},\ \{1,2,3\}\}.
\]
Therefore, $\Delta$ is the simplex $\Delta=\langle\{1,2,3\}\rangle$, that is

\begin{figure}[h]
\begin{center}
\unitlength 1mm 
\linethickness{0.4pt}
\ifx\plotpoint\undefined\newsavebox{\plotpoint}\fi 
\begin{picture}(54,40.5)(0,0)
\multiput(15.25,8.75)(.033702213,.056841046){497}{\line(0,1){.056841046}}
\put(15.25,8.5){\line(1,0){33.25}}
\multiput(48.5,8.5)(-.033673469,.057142857){490}{\line(0,1){.057142857}}
\multiput(32,25.5)(-.03372093,-.05465116){215}{\line(0,-1){.05465116}}
\multiput(35.5,20.5)(-.03365385,-.04967949){156}{\line(0,-1){.04967949}}
\multiput(38.75,14.75)(-.0333333,-.0433333){75}{\line(0,-1){.0433333}}
\put(30.75,39.5){$1$}
\put(11,8.5){$2$}
\put(52,8.5){$3$}
\end{picture}
\caption{}
\label{f2}
\end{center}
\end{figure}
\end{Example1}

For a simplicial complex $\Delta$ on the vertex set $[n]$ and a field $k$, one may consider the polynomial ring $k[x_1,\ldots,x_n]$ and the square-free monomial ideal  
	\[I_{\Delta}=(x_F\ :\ F\notin\Delta)
\]
called the \textit{Stanley--Reisner ideal} of $\Delta$. Here $x_F$ is the square-free monomial $$x_F=x_{i_1}\cdots x_{i_s},$$ where $F=\{i_1,\ldots,i_s\}$.

For the simplicial complex considered in Figure \ref{f1}, the Stanley--Reisner ideal is
	\[I_{\Delta}=(x_1x_4,\ x_2x_4).
\]
The factor ring
	\[k[\Delta]=\frac{k[x_1,\ldots,x_n]}{I_{\Delta}}
\]
is called the \textit{Stanley--Reisner ring} of $\Delta$.

For each subset $F\subset[n]$, we set $p_F$ to be the prime monomial ideal generated by the variables $x_i$ such that $i\in F$, that is
	\[p_F=(x_i\ :\ i\in F).
\]

The standard primary decomposition of the Stanley--Reisner ideal of a simplicial complex can be written just looking to its facets.

\begin{Proposition1}[\cite{BH}]\label{primdec} Let $\Delta$ be a simplicial complex. The standard primary decomposition of $I_{\Delta}$ is
	\[I_{\Delta}=\bigcap_{F\in\mathcal{F}(\Delta)}p_{F^c}.
\]
\end{Proposition1}

\begin{Example1}\rm For the simplicial complex from Figure \ref{f1}, the standard primary decomposition of $I_{\Delta}$ is
	\[I_{\Delta}=(x_4)\cap(x_1,x_2).
\]
\end{Example1}
For a simplicial complex $\Delta$, the dimension of its Stanley--Reisner ring can be easily determined.

\begin{Proposition1}[\cite{BH}] Let $\Delta$ be a simplicial complex. Then $$\dim( k[\Delta])=\dim(\Delta)+1.$$
\end{Proposition1}

In \cite{Fa}, S. Faridi introduced the notion of facet ideal of a simplicial complex. The \textit{facet ideal}, denoted by $I(\Delta)$, is the square-free monomial ideal generated by the monomials corresponding to the facets of $\Delta$, that is
	\[I(\Delta)=(x_F\ :\ F\in\mathcal{F}(\Delta)).
\]
If $\Delta$ is the simplicial complex from Figure \ref{f1}, then the facet ideal is
	\[I(\Delta)=(x_1x_2x_3,\ x_3x_4).
\]

\begin{Definition1}\rm Let $\Delta$ be a simplicial complex and $F$ a face of $\Delta$. The simplicial complex
	\[\del(F,\Delta)=\{G\in\Delta\ \mid\ G\cap F=\emptyset\}
\]
is called the \textit{deletion of $F$ from $\Delta$}.
\end{Definition1}

\begin{Example1}\rm Let $\Delta$ be the following simplicial complex: 
\newpage
\begin{figure}[ht]
\begin{center}
\unitlength 1mm 
\linethickness{0.4pt}
\ifx\plotpoint\undefined\newsavebox{\plotpoint}\fi 
\begin{picture}(87.5,40.25)(0,0)
\multiput(15.25,8.75)(.033702213,.056841046){497}{\line(0,1){.056841046}}
\put(15.25,8.5){\line(1,0){33.25}}
\multiput(48.5,8.5)(-.033673469,.057142857){490}{\line(0,1){.057142857}}
\multiput(32,25.5)(-.03372093,-.05465116){215}{\line(0,-1){.05465116}}
\multiput(35.5,20.5)(-.03365385,-.04967949){156}{\line(0,-1){.04967949}}
\multiput(38.75,14.75)(-.0333333,-.0433333){75}{\line(0,-1){.0433333}}
\put(32.25,36.5){\line(1,0){34}}
\multiput(66.25,36.5)(-.033681214,-.053605313){527}{\line(0,-1){.053605313}}
\put(48.75,8.75){\line(1,0){35}}
\multiput(83.75,8.75)(-.03371869,.053949904){519}{\line(0,1){.053949904}}
\multiput(59.25,14.75)(.03360215,.05645161){186}{\line(0,1){.05645161}}
\multiput(65,14.25)(.03348214,.0625){112}{\line(0,1){.0625}}
\multiput(70.25,14)(.0333333,.0722222){45}{\line(0,1){.0722222}}
\put(12,6.5){$1$}
\put(30,39){$2$}
\put(66.5,39){$4$}
\put(47.5,5){$3$}
\put(87,6.5){$5$}
\end{picture}
\end{center}
\caption{}
\label{f5}
\end{figure}
Then $\del(\{1,3\},\Delta)=\langle\{2,4\},\ \{4,5\}\rangle$.
\end{Example1}

\begin{Definition1}\rm Let $\Delta$ be a simplicial complex and $F$ a face of $\Delta$. The simplicial complex
	\[\lk(F,\Delta)=\{G\in\Delta\ \mid\ G\cap F=\emptyset\ \mbox{and}\ G\cup F\in\Delta\}
\]
is called the \textit{link of $F$ in $\Delta$}.
\end{Definition1}

\begin{Example1}\rm Let $\Delta$ be the simplicial complex from Figure \ref{f5}. Then $$\lk(\{1,3\},\Delta)=\langle\{2\}\rangle.$$ 
\end{Example1}

\begin{Remark1}\rm Let $\Delta$ be a simplicial complex. Then $\lk(\emptyset,\Delta)=\Delta$. 
\end{Remark1}
A particular class of simplicial complexes is the class of pure simplicial complexes. Recall that a simplicial complex is \textit{pure} if all its facets have the same dimension.

The simplicial complex from Figure \ref{f1} is not pure since $\dim(\{1,2,3\})=2$ and $\dim(\{3,4\})=1$. Next, we consider an example of a pure simplicial complex.
\begin{Example1}\label{sc1}\rm Let $\Delta$ be the simplicial complex

\begin{figure}[h]
\begin{center}
\unitlength 1mm 
\linethickness{0.4pt}
\ifx\plotpoint\undefined\newsavebox{\plotpoint}\fi 
\begin{picture}(78.5,34.75)(0,0)
\multiput(25,32)(-.0337136929,-.0529045643){482}{\line(0,-1){.0529045643}}
\multiput(57.5,31.75)(-.0337136929,-.0529045643){482}{\line(0,-1){.0529045643}}
\put(8.75,6.5){\line(1,0){32.5}}
\put(41.25,6.25){\line(1,0){32.5}}
\multiput(41.25,6.5)(-.0336842105,.0542105263){475}{\line(0,1){.0542105263}}
\multiput(73.75,6.25)(-.0336842105,.0542105263){475}{\line(0,1){.0542105263}}
\put(25,23.25){\line(-2,-3){8}}
\put(57.5,23){\line(-2,-3){8}}
\multiput(27.75,19.25)(-.033707865,-.046348315){178}{\line(0,-1){.046348315}}
\multiput(60.25,19)(-.033707865,-.046348315){178}{\line(0,-1){.046348315}}
\multiput(30.25,15)(-.033653846,-.043269231){104}{\line(0,-1){.043269231}}
\multiput(62.75,14.75)(-.033653846,-.043269231){104}{\line(0,-1){.043269231}}
\put(32.5,9.5){$F_1$}
\put(65.5,9.5){$F_2$}
\put(5.25,6){$1$}
\put(24.75,34.5){$2$}
\put(40.5,10){$3$}
\put(57.75,34.5){$4$}
\put(77,6){$5$}
\end{picture}
\end{center}
\caption{}
\label{f3}
\end{figure}

One may note that $\Delta$ is pure since both its facets are of dimension $2$.
\end{Example1}
	\[
\]
\subsection{Classes of pure simplicial complexes}
We recall some known classes of pure simplicial complexes. We give examples and counter-examples of such simplicial complexes and we establish a hierachy on them.

The most larger and important class of pure simplicial complexes is the class of Cohen--Macaulay simplicial complexes. Recall that a simplicial complex $\Delta$ is Cohen--Macaulay if the Stanley--Reisner ring $k[\Delta]$ is Cohen--Macaulay. 

For a proof of the following result, see for instance W. Bruns and J. Herzog \cite[p. 202]{BH} or R.H. Villareal \cite[p. 145]{Vi}:

\begin{Proposition1}[\cite{BH}] A Cohen--Macaulay simplicial complex is pure.
\end{Proposition1}

We consider the simplicial complex $\Delta=\langle\{1,2\},\ \{3,4\}\rangle$ with the vertex set $\{1,2,3,4\}$, that is
\begin{figure}[h]
\begin{center}
\unitlength 1mm 
\linethickness{0.4pt}
\ifx\plotpoint\undefined\newsavebox{\plotpoint}\fi 
\begin{picture}(28.5,33.75)(0,0)
\put(5.5,30.25){\line(0,-1){19.75}}
\put(21.5,30.25){\line(0,-1){19.75}}
\put(6,6.5){$1$}
\put(6,32.5){$2$}
\put(21.5,6.5){$3$}
\put(21.5,32.5){$4$}
\end{picture}
\end{center}
\caption{}
\label{f6}
\end{figure}

\noindent One may note that $\Delta$ is a pure simplicial complex, but the Stanley--Reisner ring
	\[k[\Delta]=\frac{k[x_1,\ldots,x_4]}{(x_1x_3,\ x_1x_4,\ x_2x_3,\ x_2x_4)}
\]
is not Cohen--Macaulay since $\dim(k[\Delta])=2$ and $\depth(k[\Delta])=1$.

The following result is a characterization of the Cohen--Macaulay simplicial complexes in terms of simplicial homology.
\begin{Theorem1}[Reisner \cite{Rei}] Let $\Delta$ be a simplicial complex and $k$ a field. The following conditions are equivalent:
\begin{itemize}
	\item[(i)] $\Delta$ is Cohen--Macaulay over $k$.
	\item[(ii)] $\widetilde{H}_i(\lk(F,\Delta);k)=0$ for all $F\in\Delta$ and all $i\leq\dim(\lk(F,\Delta))$.
\end{itemize}
\end{Theorem1}

A simplicial complex $\Delta$ is called \it disconnected\rm$\ $ if the vertex set $V$ of $\Delta$ is a disjoint union $V=V_1\cup V_2$ such that no face of $\Delta$ has vertices in both $V_1$ and $V_2$. Otherwise, $\Delta$ is \it connected\rm.

A useful corollary of the Reisner's Theorem is the following.

\begin{Corollary1}[\cite{BH}] Let $\Delta$ be a simplicial complex and $k$ a field. If $\Delta$ is Cohen--Macaulay over $k$, then $\Delta$ is connected.
\end{Corollary1}

The simplicial complex from Figure \ref{f6} is not connected, hence, it is not Cohen--Macaulay.

In general, the property of a simplicial complex to be Cohen--Macaulay depends on the characteristic of the base field. In  the following, we consider a simplicial complex for which the property of being Cohen--Macaulay depends on the characteristic.

\begin{Example1}\rm Let $\Delta$ be the triangulation of the projective plane from Figure \ref{f14}. Then $k[\Delta]$ is Cohen--Macaulay if $\chr(k)\neq2$ and it is not Cohen--Macaulay if $k$ is of characteristic $2$ (see, for instance, W. Bruns and J. Herzog \cite[p. 228]{BH}).
\begin{figure}[h]
\begin{center}
\unitlength 1mm 
\linethickness{0.4pt}
\ifx\plotpoint\undefined\newsavebox{\plotpoint}\fi 
\begin{picture}(75.75,69.5)(0,0)
\multiput(23.75,65.25)(-.03373494,-.061445783){415}{\line(0,-1){.061445783}}
\multiput(72.75,40)(-.03373494,.061445783){415}{\line(0,1){.061445783}}
\multiput(72.75,40.5)(-.03373494,-.061445783){415}{\line(0,-1){.061445783}}
\put(23.75,65.25){\line(1,0){35}}
\put(23.75,15.25){\line(1,0){35}}
\multiput(40.25,48)(-.033636364,-.055454545){275}{\line(0,-1){.055454545}}
\put(31,32.75){\line(1,0){17.5}}
\multiput(48.5,32.75)(-.033730159,.0625){252}{\line(0,1){.0625}}
\multiput(40,48.5)(-.033713693,.034751037){482}{\line(0,1){.034751037}}
\multiput(72.75,40.25)(-.11395349,-.03372093){215}{\line(-1,0){.11395349}}
\multiput(48.25,33)(.033653846,-.055288462){312}{\line(0,-1){.055288462}}
\multiput(23.75,15)(.03372093,.08372093){215}{\line(0,1){.08372093}}
\put(60.75,67.25){$1$}
\put(20.75,67.25){$2$}
\put(5.75,40.5){$3$}
\put(20.75,11.5){$1$}
\put(60.75,11.5){$2$}
\put(74.5,40.5){$3$}
\put(51.75,30.5){$5$}
\put(26.25,30.5){$4$}
\put(39.75,51){$6$}
\multiput(9.75,39.75)(.033711217,-.059069212){419}{\line(0,-1){.059069212}}
\multiput(23.63,65.13)(.0336758,-.14840183){219}{\line(0,-1){.14840183}}
\multiput(10,39.5)(.1050995,-.03358209){201}{\line(1,0){.1050995}}
\multiput(58.75,65.13)(-.036458333,-.033730159){504}{\line(-1,0){.036458333}}
\multiput(40.25,48.13)(.13782051,-.03365385){234}{\line(1,0){.13782051}}
\multiput(48.38,32.75)(-.046726755,-.033681214){527}{\line(-1,0){.046726755}}
\end{picture}
\end{center}
\caption{}
\label{f14}
\end{figure}
\end{Example1}
An example of a simplicial complex which is Cohen--Macaulay over any field is the Dunce Hat, Figure \ref{f7} (see E.C. Zeeman \cite{Ze}):
\begin{figure}[h]
\begin{center}
\unitlength 1mm 
\linethickness{0.4pt}
\ifx\plotpoint\undefined\newsavebox{\plotpoint}\fi 
\begin{picture}(68.75,63.5)(0,0)
\multiput(37.25,58.5)(-.0337045721,-.0565650645){853}{\line(0,-1){.0565650645}}
\put(8.5,10.25){\line(1,0){57.5}}
\multiput(66,10.25)(-.0337278107,.0573964497){845}{\line(0,1){.0573964497}}
\put(37.5,58.75){\line(0,-1){23}}
\multiput(37.5,35.75)(-.0472561,-.03353659){164}{\line(-1,0){.0472561}}
\multiput(29.75,30.25)(.03350515,-.06443299){97}{\line(0,-1){.06443299}}
\put(33,24){\line(1,0){9.5}}
\put(42.5,24){\line(2,3){4}}
\multiput(46.5,30)(-.05487805,.03353659){164}{\line(-1,0){.05487805}}
\multiput(37.5,35.5)(-.03373016,-.08928571){126}{\line(0,-1){.08928571}}
\multiput(33.25,24.25)(-.059036145,-.03373494){415}{\line(-1,0){.059036145}}
\multiput(37.75,35)(.03358209,-.08022388){134}{\line(0,-1){.08022388}}
\multiput(42.25,24.25)(-.039156627,-.03373494){415}{\line(-1,0){.039156627}}
\multiput(26,10.25)(.03365385,.06610577){208}{\line(0,1){.06610577}}
\multiput(33,24)(-.12946429,.03348214){112}{\line(-1,0){.12946429}}
\multiput(18.5,27.75)(.15,.0333333){75}{\line(1,0){.15}}
\multiput(29.75,30.25)(-.0335821,.1865672){67}{\line(0,1){.1865672}}
\multiput(37.75,58.25)(-.03369565,-.12065217){230}{\line(0,-1){.12065217}}
\multiput(42.25,23.75)(.03348214,-.11830357){112}{\line(0,-1){.11830357}}
\multiput(42.5,24)(.057598039,-.03370098){408}{\line(1,0){.057598039}}
\multiput(66,10.25)(-.033703072,.034129693){586}{\line(0,1){.034129693}}
\multiput(46.25,30.25)(.1730769,-.0336538){52}{\line(1,0){.1730769}}
\multiput(46.5,29.75)(-.03125,1.84375){8}{\line(0,1){1.84375}}
\multiput(46.25,44.5)(-.033730159,-.034722222){252}{\line(0,-1){.034722222}}
\put(37,61.5){$1$}
\put(24.5,43.5){$3$}
\put(15.5,28.75){$2$}
\put(3.75,8){$1$}
\put(25.75,6.5){$3$}
\put(46.5,6.5){$2$}
\put(68.75,8){$1$}
\put(57.5,29.25){$2$}
\put(48.25,45){$3$}
\put(34.5,36.5){$6$}
\put(26.75,31.75){$5$}
\put(33,21){$4$}
\put(45,23){$8$}
\put(44,32){$7$}
\end{picture}
\end{center}
\caption{}
\label{f7}
\end{figure}

A smaller class of pure simplicial complexes is the class of constructible simplicial complexes.

\begin{Definition1}\rm A pure simplicial complex $\Delta$ is \textit{constructible} if it can be obtained by the following recursive procedure:
\begin{itemize}
	\item[(i)] Any simplex is constructible.
	\item[(ii)] If $\Delta_1,\ \Delta_2$ are $d-$dimensional constructible simplicial complexes and $\Delta_1\cap\Delta_2$ is a constructible simplicial complex of dimension $d-1$, then $\Delta_1\cup\Delta_2$ is a $d-$dimensional constructible simplicial complex. 
\end{itemize}
\end{Definition1}

The following simplicial complex is given by M. Hachimori \cite{Ha} and it is a constructible simplicial complex:
	
\begin{figure}[h]
\begin{center}
\unitlength 1mm 
\linethickness{0.4pt}
\ifx\plotpoint\undefined\newsavebox{\plotpoint}\fi 
\begin{picture}(52.25,50.25)(0,0)
\put(3.5,19.75){\line(1,0){46.25}}
\put(3.5,28.5){\line(1,0){46.25}}
\multiput(27,48.25)(-.074519231,-.033653846){312}{\line(-1,0){.074519231}}
\put(3.75,37.75){\line(0,-1){17.5}}
\multiput(3.75,20.25)(.073636364,-.033636364){275}{\line(1,0){.073636364}}
\multiput(24,11)(.098076923,.033653846){260}{\line(1,0){.098076923}}
\put(49.5,19.75){\line(0,1){18}}
\multiput(49.5,37.75)(-.069749216,.03369906){319}{\line(-1,0){.069749216}}
\multiput(27.25,48)(-.033687943,-.045508274){423}{\line(0,-1){.045508274}}
\multiput(13,28.75)(.033639144,-.054281346){327}{\line(0,-1){.054281346}}
\thicklines
\put(4,37.5){\line(1,0){45.5}}
\multiput(13,29)(.08508403,.03361345){238}{\line(1,0){.08508403}}
\multiput(4,37.25)(.035714286,-.033673469){245}{\line(1,0){.035714286}}
\multiput(12.75,28.25)(-.034693878,-.033673469){245}{\line(-1,0){.034693878}}
\multiput(18.75,19.75)(-.03991597,-.03361345){119}{\line(-1,0){.03991597}}
\multiput(18.75,20)(.075961538,.033653846){260}{\line(1,0){.075961538}}
\multiput(38.25,28.5)(.043071161,.033707865){267}{\line(1,0){.043071161}}
\multiput(38.5,28.25)(.04456522,-.03369565){230}{\line(1,0){.04456522}}
\multiput(48.75,20.5)(.032609,-.032609){23}{\line(0,-1){.032609}}
\multiput(31.75,19.5)(.04910714,-.03348214){112}{\line(1,0){.04910714}}
\put(26.75,50.25){$0$}
\put(15.25,45){$3$}
\put(.5,38.25){$2$}
\put(.5,28.75){$3$}
\put(.5,19){$0$}
\put(23.5,7.5){$2$}
\put(12.75,12.25){$3$}
\put(38.75,13.5){$1$}
\put(51.5,19){$0$}
\put(52,28.5){$1$}
\put(52.25,38){$2$}
\put(38.75,44.75){$1$}
\put(20,34.5){$9$}
\put(11.5,24.5){$8$}
\put(32.25,33){$4$}
\put(39.25,31){$5$}
\put(19.5,21.5){$7$}
\put(30,21.5){$6$}
\thinlines
\multiput(33.5,37.75)(.03355705,.03691275){149}{\line(0,1){.03691275}}
\multiput(19.25,37.75)(-.03348214,.046875){112}{\line(0,1){.046875}}
\thicklines
\multiput(27.75,48)(.033653846,-.0625){312}{\line(0,-1){.0625}}
\multiput(28,48.25)(.033653846,-.0625){312}{\line(0,-1){.0625}}
\multiput(38.25,28.5)(-.033687943,-.042553191){423}{\line(0,-1){.042553191}}
\multiput(38.5,28.75)(-.033687943,-.042553191){423}{\line(0,-1){.042553191}}
\multiput(38.75,28.5)(-.033687943,-.042553191){423}{\line(0,-1){.042553191}}
\end{picture}
\end{center}
\caption{}
\label{f8}
\end{figure}

In between Cohen--Macaulay simplicial complexes and constructible simplicial complexes there exists the following relation:

\begin{Theorem1}$($\cite[p. 228]{BH}$)$\label{constCM} A constructible simplicial complex is Cohen--Macaulay over any field.
\end{Theorem1}

Not all Cohen--Macaulay simplicial complexes are constructible. The Dunce Hat (Figure \ref{f7}) is such a simplicial complex (see M. Hachimori \cite{Ha}). 

A class of pure simplicial complexes with many important combinatorial properties and that it is often used to prove the Cohen-Macaulayness is the class of pure shellable simplicial complexes. There are several equivalent definitions for a pure shellable simplicial complex. For the equivalence of these definitions, see for instance W. Bruns and J. Herzog \cite[p. 207]{BH}.

\begin{Definition1}\rm A pure simplicial complex $\Delta$ is called \textit{shellable} if on the facets of $\Delta$ can be given a linear order $F_1,\ldots, F_m$ such that one of the following equivalent conditions is satisfied: 
\begin{itemize}
	\item[(a)] The simplicial complex $\langle F_1,\ldots,F_{i-1}\rangle\cap\langle F_i\rangle$ is generated by a nonempty set of maximal proper faces of $\langle F_i\rangle$ for all $i$, $2\leq i\leq m$.
	\item[(b)] The set $\{F\ :\ F\in\langle F_1,\ldots,F_{i}\rangle,\ F\notin\langle F_1,\ldots,F_{i-1}\rangle\}$ has a unique minimal element for all $i$, $2\leq i\leq m$.
	\item[(c)] For all $i,j$, $1\leq j<i\leq m$, there exist some $v\in F_i\setminus F_j$ and some $k\in\{1,2,\ldots, i-1\}$ with $F_i\setminus F_k=\{v\}$.
\end{itemize}
\end{Definition1}

A linear order of the facets satisfying one of the equivalent conditions (a), (b) or (c) is called \textit{a shelling} of $\Delta$.  

\begin{Example1}\rm Let $\Delta$ be the simplicial complex given in Figure \ref{f4}
\begin{figure}[h]
\begin{center}
\unitlength 1mm 
\linethickness{0.4pt}
\ifx\plotpoint\undefined\newsavebox{\plotpoint}\fi 
\begin{picture}(69.75,31.25)(0,0)
\multiput(24.25,28.25)(-.0337259101,-.0524625268){467}{\line(0,-1){.0524625268}}
\multiput(53.5,28.25)(-.0337259101,-.0524625268){467}{\line(0,-1){.0524625268}}
\put(8.5,3.75){\line(1,0){29}}
\put(37.75,3.75){\line(1,0){29}}
\multiput(37.5,3.75)(-.0336787565,.064119171){386}{\line(0,1){.064119171}}
\multiput(66.75,3.75)(-.0336787565,.064119171){386}{\line(0,1){.064119171}}
\put(25.25,28.25){\line(1,0){28.5}}
\put(5,3.5){$1$}
\put(23.5,29.75){$2$}
\put(36.5,6.75){$3$}
\put(55.25,29.75){$4$}
\put(69.75,3.5){$5$}
\multiput(24,21)(-.0337301587,-.0496031746){252}{\line(0,-1){.0496031746}}
\multiput(52.5,21)(-.0337301587,-.0496031746){252}{\line(0,-1){.0496031746}}
\multiput(25.75,17.25)(-.033707865,-.049157303){178}{\line(0,-1){.049157303}}
\multiput(54.25,17.25)(-.033707865,-.049157303){178}{\line(0,-1){.049157303}}
\multiput(27.75,13.25)(-.033482143,-.046875){112}{\line(0,-1){.046875}}
\multiput(56.25,13.25)(-.033482143,-.046875){112}{\line(0,-1){.046875}}
\multiput(30,9.5)(-.03358209,-.03731343){67}{\line(0,-1){.03731343}}
\multiput(58.5,9.5)(-.03358209,-.03731343){67}{\line(0,-1){.03731343}}
\put(31.75,22){\line(1,0){13.5}}
\put(33.75,18.5){\line(1,0){9.5}}
\put(36.5,15.5){\line(1,0){4.5}}
\put(29,25.25){\line(1,0){18.25}}
\put(15.25,6){$F_1$}
\put(46.5,24){$F_2$}
\put(60.5,6){$F_3$}
\end{picture}
\end{center}
\caption{}
\label{f4}
\end{figure}

One may easily check that $F_1,\ F_2,\ F_3$ (as in the figure) is a shelling of $\Delta$.
\end{Example1}
\begin{Theorem1}[\cite{Hi}] A shellable simplicial complex is Cohen--Macaulay over any field.
\end{Theorem1}

Moreover, we have the following stronger result:

\begin{Theorem1}[\cite{St}] Any shellable simplicial complex is constructible.
\end{Theorem1}

Not every constructible simplicial complex is shellable. However, it seems to be difficult to find examples of such simplicial complexes. M. Hachimori \cite{Ha} proved that the simplicial complex from Figure \ref{f8} is a constructible simplicial complex which is not shellable. 

The following class of pure simplicial complexes was introduced by L.J. Billera and J.S. Prova \cite{BiPr}.
\begin{Definition1}\rm Let $\Delta$ be a pure simplicial complex. $\Delta$ is called \it{vertex-decomposable }\rm if $\Delta=\{\emptyset\}$ or there exists a vertex $v\in\Delta$ such that both $\del(v,\Delta)$ and $\lk(v,\Delta)$ are vertex-decomposable.
\end{Definition1}

\begin{Example1}\rm Let $\Delta$ be the simplicial complex from Figure \ref{f9} with the facets  $$\mathcal{F}(\Delta)=\{\{1,2,3\},\ \{1,3,4\}\}.$$

\begin{figure}[ht]
\begin{center}
\unitlength 1mm 
\linethickness{0.4pt}
\ifx\plotpoint\undefined\newsavebox{\plotpoint}\fi 
\begin{picture}(64,43)(0,0)
\put(35.25,37.75){\line(0,-1){26.25}}
\multiput(35.25,11.5)(-.058210784,.03370098){408}{\line(-1,0){.058210784}}
\multiput(59.5,25.25)(-.063471503,.033678756){386}{\line(-1,0){.063471503}}
\multiput(59.25,25.5)(-.057328605,-.033687943){423}{\line(-1,0){.057328605}}
\multiput(12,25.25)(.060233161,.033678756){386}{\line(1,0){.060233161}}
\put(30.25,29){\line(0,-1){9}}
\put(40.75,29){\line(0,-1){9}}
\put(25,26.5){\line(0,-1){4.75}}
\put(46.75,26.5){\line(0,-1){4.75}}
\put(20.5,25.5){\line(0,-1){2.25}}
\put(51.75,25.5){\line(0,-1){2.25}}
\put(35.25,41){$1$}
\put(8.75,25){$2$}
\put(35.25,6.75){$3$}
\put(62,25){$4$}
\end{picture}
\end{center}
\caption{}
\label{f9}
\end{figure}
One may note that, for the vertex $\{1\}$, we have 
	\[\del(\{1\},\Delta)=\langle\{2,3\},\ \{3,4\}\rangle=\lk(\{1\},\Delta).
\]
Let us denote $\Delta'=\langle\{2,3\},\ \{3,4\}\rangle$. For the vertex $\{3\}$ in $\Delta'$, we obtain that 
\[\del(\{3\},\Delta)=\lk(\{3\},\Delta)=\langle\{2\},\ \{4\}\rangle.
\]
One may easily check that the last simplicial complex is vertex-decomposable.
\end{Example1}
\newpage
In between shellable simplicial complexes and vertex-decomposable simplicial complexes there exists the following relation.

\begin{Theorem1}[\cite{BiPr}] Any vertex-decomposable simplicial complex is shellable.
\end{Theorem1}

We consider a $2-$dimensional simplicial complex which is shellable and it is not vertex-decomposable. This example is also due to M. Hachimori \cite{Ha1}.

\begin{figure}[h]
\begin{center}
\unitlength 1mm 
\linethickness{0.4pt}
\ifx\plotpoint\undefined\newsavebox{\plotpoint}\fi 
\begin{picture}(54,38.5)(0,0)
\multiput(25.25,36.5)(-.033739837,-.050813008){615}{\line(0,-1){.050813008}}
\put(4.5,5.25){\line(1,0){47}}
\multiput(51.5,5.25)(-.0336970475,.0401155327){779}{\line(0,1){.0401155327}}
\multiput(25.25,36.5)(-.03370787,-.09129213){178}{\line(0,-1){.09129213}}
\put(19.5,20){\line(0,-1){7.5}}
\put(19.5,12.5){\line(1,0){11.5}}
\put(31,12.5){\line(0,1){8.25}}
\put(19.25,20){\line(1,0){12}}
\multiput(31.25,20)(-.03360215,.09005376){186}{\line(0,1){.09005376}}
\multiput(31.5,19.75)(.03571429,.03361345){119}{\line(1,0){.03571429}}
\multiput(31.5,19.75)(.045930233,-.03372093){430}{\line(1,0){.045930233}}
\multiput(51.25,5.25)(-.09534884,.03372093){215}{\line(-1,0){.09534884}}
\multiput(30.75,12.5)(-.05269058,.03363229){223}{\line(-1,0){.05269058}}
\multiput(19,20)(-.133333,.033333){30}{\line(-1,0){.133333}}
\multiput(19.25,20.25)(-.03372093,-.034883721){430}{\line(0,-1){.034883721}}
\multiput(4.75,5.25)(.07091346,.03365385){208}{\line(1,0){.07091346}}
\multiput(19.75,12.25)(-.03353659,-.04115854){164}{\line(0,-1){.04115854}}
\multiput(30.5,12.75)(.03361345,-.06302521){119}{\line(0,-1){.06302521}}
\multiput(34.5,5.25)(-.06860465,.03372093){215}{\line(-1,0){.06860465}}
\put(24.5,37){$3$}
\put(12.75,22){$2$}
\put(1.75,4){$1$}
\put(13,2){$3$}
\put(34.75,2){$2$}
\put(53,4){$1$}
\put(38,23){$2$}
\put(21,21){$5$}
\put(17,13){$4$}
\put(32,13){$7$}
\put(28,21){$6$}
\end{picture}
\end{center}
\caption{}
\label{f10}
\end{figure}

The smallest example which is shellable but not vertex-decomposable was found by S. Moriyama and F. Takeuchi \cite{MT} and it has $6$ vertices and $9$ facets.

The smallest class of pure simplicial complexes that we recall here is that of shifted simplicial complexes.
\begin{Definition1}\rm A pure simplicial complex with the vertex set $\{v_1,\ldots, v_n\}$ is \textit{shifted} if there exists a labelling of the vertices by $1$ to $n$ such that, for any face $F$, replacing any vertex $v_i\in F$
by a vertex with a smaller label that does not belong to $F$, it results a set which is also a face.
\end{Definition1}

\begin{Example1}\rm Let $\Delta$ be the following simplicial complex of dimension $1$ from Figure \ref{f11}.
One may easily check that $\Delta$ is a shifted simplicial complex.
\newpage
\begin{figure}[h]
\begin{center}
\unitlength 1mm 
\linethickness{0.4pt}
\ifx\plotpoint\undefined\newsavebox{\plotpoint}\fi 
\begin{picture}(44.75,38.25)(0,0)
\multiput(7,9.25)(.033695652,.05923913){460}{\line(0,1){.05923913}}
\put(7,9){\line(1,0){35.75}}
\multiput(42.5,8.75)(-.033717105,.045230263){608}{\line(0,1){.045230263}}
\put(3.75,7.25){$1$}
\put(22,38.25){$2$}
\put(44.75,7.25){$3$}
\end{picture}
\end{center}
\caption{}
\label{f11}
\end{figure}
\end{Example1}

\begin{Theorem1}[\cite{BK}] Any shifted simplicial complex is vertex-decomposable.
\end{Theorem1}

Not any vertex-decomposable simplicial complex is shifted. Indeed, let us consider the simplicial complex $\Delta=\langle\{v_1,v_2\},\ \{v_1,v_3\},\ \{v_2,v_4\},\ \{v_3,v_4\}\rangle$, that is
	\[
\]
\begin{figure}[h]
\begin{center}
\unitlength 1mm 
\linethickness{0.4pt}
\ifx\plotpoint\undefined\newsavebox{\plotpoint}\fi 
\begin{picture}(40.5,38.25)(0,0)
\put(7.75,13.5){\framebox(28.5,24.25)[cc]{}}
\put(3.5,38.25){$v_1$}
\put(38,38.25){$v_2$}
\put(38,13.5){$v_4$}
\put(3.5,13.5){$v_3$}
\end{picture}
\end{center}
\caption{}
\label{f12}
\end{figure}

It is not hard to note that this is a vertex-decomposable simplicial complex. But this simplicial complex is not shifted. Indeed, if we consider a total order on the vertex set such that $v_1<v_3$, then looking to the face $\{v_3,v_4\}$ and replacing $v_3$ by $v_1$ we obtain $\{v_1,v_4\}$ which is not a face in $\Delta$. If we consider a total order on the vertices such that $v_3<v_1$, then looking to the face $\{v_1,v_2\}$ and replacing $v_1$ by $v_3$ we obtain the set $\{v_3,v_2\}$ which is not a face in $\Delta$.

	\[
\]
We conclude that, for pure simplicial complexes, we have the following strict implications:
	\[
\]
	\[\mbox{shifted}\Rightarrow\mbox{vertex-decomposable}\Rightarrow\mbox{shellable}\Rightarrow\mbox{constructible}\Rightarrow\mbox{Cohen--Macaulay}.
\]
	\[
\]
\subsection{Simplicial complexes and Alexander duality}
Many properties of a simplicial complex can be obtained by studying the Alexander dual and its Stanley--Reisner ideal.

\begin{Definition1}\rm Let $\Delta$ be a simplicial complex on the vertex set $[n]$. \textit{The Alexander dual} is the simplicial complex
	\[\Delta^{\vee}=\{F^c\ \mid\ F\notin\Delta\},
\]
where we denote $F^c=[n]\setminus F$. 
\end{Definition1}

\begin{Example1}\rm Let $\Delta$ be the following simplicial complex on the vertex set $\{1,2,3,4\}$: 

\begin{figure}[h]
\begin{center}
\unitlength 1mm 
\linethickness{0.4pt}
\ifx\plotpoint\undefined\newsavebox{\plotpoint}\fi 
\begin{picture}(75,40.75)(0,0)
\multiput(7.5,8.5)(.033737024,.053633218){578}{\line(0,1){.053633218}}
\put(7.5,8.25){\line(1,0){37}}
\multiput(44.5,8.25)(-.033681214,.059297913){527}{\line(0,1){.059297913}}
\put(43.5,8.5){\line(1,0){28.5}}
\put(27.5,40.75){$1$}
\put(4.25,6.75){$2$}
\put(44.75,10){$3$}
\put(75,6.75){$4$}
\multiput(17.25,12.25)(.03361345,.05777311){238}{\line(0,1){.05777311}}
\multiput(21.5,12.5)(.03365385,.06089744){156}{\line(0,1){.06089744}}
\multiput(26.25,12)(.03350515,.06443299){97}{\line(0,1){.06443299}}
\multiput(27,39.25)(.0498903509,-.0337171053){912}{\line(1,0){.0498903509}}
\end{picture}

\end{center}
\caption{}
\label{f13}
\end{figure} 
Then, the Alexander dual of $\Delta$ is $\Delta^{\vee}=\langle\{2\},\ \{1,3\}\rangle$, that is
	\[
\]
\begin{figure}[h]
\begin{center}
\unitlength 1mm 
\linethickness{0.4pt}
\ifx\plotpoint\undefined\newsavebox{\plotpoint}\fi 
\begin{picture}(42.25,25)(0,0)
\put(9.75,8.5){\line(1,0){29.75}}
\put(20.25,21.5){\circle*{.71}}
\put(20.5,25){$2$}
\put(7,9.25){$1$}
\put(41.25,9.25){$3$}
\end{picture}

\end{center}
\caption{}
\label{fdual}
\end{figure}
\end{Example1}

\begin{Remark1}\rm Let $\Delta$ be a simplicial complex. Then $(\Delta^{\vee})^{\vee}=\Delta$.
\end{Remark1}

For the simplicial complex $\Delta$ with the vertex set $[n]$, one may also consider the simplicial complex denoted by $\Delta^c$ whose facets are $\mathcal{F}(\Delta^c)=\{F^c\ \mid\ F\in\mathcal{F}(\Delta)\}$.

The Alexander dual of a simplicial complex $\Delta$ and the simplicial complex $\Delta^c$ are closely related, as can be seen in the following proposition.

\begin{Proposition1}[\cite{HeHiZh}] Let $\Delta$ be a simplicial complex. Then
	\[I_{\Delta^{\vee}}=I(\Delta^c).
\]
\end{Proposition1}
The above proposition allows us to write the Stanley--Reisner ideal of the Alexander dual of a simplicial complex $\Delta$ just looking to the facets of $\Delta$. For the simplicial complex in Figure \ref{f13}, the Stanley--Reisner ideal of the Alexander dual associated to $\Delta$ is 
	\[I_{\Delta^{\vee}}=(x_4, x_2x_3, x_1x_2).
\]

Between some invariants of a simplicial complex and its Alexander dual there exists some strong connections.

\begin{Theorem1}[Terai \cite{T}]\label{Terai} Let $\Delta$ be a simplicial complex. Then
	\[\projdim(I_{\Delta})=\reg(k[\Delta^{\vee}]).
\]
\end{Theorem1}
The connection between Cohen--Macaulay simplicial complex and the Alexander dual is given by the Eagon--Reiner theorem:

\begin{Theorem1}[Eagon-Reiner \cite{EaRe}] Let $k$ be a field and $\Delta$ be a pure simplicial complex on the vertex set $[n]$. $\Delta$ is Cohen--Macaulay if and only if $I_{\Delta^{\vee}}$ has a linear resolution.
\end{Theorem1}

Square-free monomial ideals with linear quotients generated in one degree are closely related to shellable simplicial complexes.

\begin{Theorem1}[\cite{HeHiZh}]\label{linquot} Let $k$ be a field and $\Delta$ be a pure simplicial complex on the vertex set $[n]$. $\Delta$ is shellable if and only if $I_{\Delta^{\vee}}$ has linear quotients.
\end{Theorem1}

	\[
\]
Therefore, if we consider the connection with the Stanley--Reisner ideal of the Alexander dual, for a pure simplicial complex $\Delta$ we have the following diagram:
	\[\begin{array}{ccccc}
	\Delta\ \mbox{is shellable}&\Longrightarrow&\Delta\ \mbox{is constructible}&\Longrightarrow&\Delta\ \mbox{is Cohen--Macaulay}\\
	\Updownarrow&&&&\Updownarrow\\
	I_{\Delta^{\vee}}\ \mbox{is an ideal with}&&&&I_{\Delta^{\vee}}\ \mbox{is an ideal with a}\\
	\mbox{linear quotients}&&&&\mbox{linear resolution} 
\end{array}.
\]
	\[
\]
%
%
%
%
%
%
%
\section{Coxeter groups}

In this section we recall some basic notions related to Coxeter groups, following J.E. Humphreys \cite{Hu}. They will be needed in the last chapter of our thesis. See A. Bj\"orner and F. Brenti \cite{BjBr} or J.E. Humphreys \cite{Hu} for more details. 
\begin{Definition1}\rm A \textit{Coxeter system} is a pair $(W,S)$ consisting of a group $W$ and a set of generators $S\subset W$ subject only to relations of the form
	\[(\sigma\sigma\,')^{m(\sigma,\sigma\,')}=1,
\]
where $m(\sigma,\sigma)=1$ and $m(\sigma,\sigma\,')=m(\sigma\,',\sigma)\geq 2$ for all $\sigma\neq \sigma\,'$ in $S$. By convention, $m(\sigma,\sigma\,')=\infty$, if no relation occurs for a pair $(\sigma,\sigma\,')$.
\end{Definition1}

The following example will be often used in Chapter \ref{sub}. For a proof, see for instance A. Bj\"orner and F. Brenti \cite{BjBr} or G. Lusztig \cite{L}.

\begin{Example1}\label{sym}\rm Let $S_n$ be the symmetric group and denote by $S$ the set of all the adjacent transpositions $s_i:=(i,i+1)$ for $1\leq i\leq n-1$. Then $(S_n,S)$ is a Coxeter system.
\end{Example1}

The elements of $S$ are called \textit{simple reflections}. Note that any simple reflection $\sigma\in S$ has order $2$ in $W$. Hence, each $w\neq 1$ in $W$ can be written in the form $w=\sigma_1 \sigma_2\cdots\sigma_r$, for some $\sigma_i\in S$ not necessarily distinct.

If $r$ is as small as possible, $r$ will be called the \textit{length} of $w$ and it will be denoted by $\ell(w)$. Any expression of $w$ as a product of $\ell(w)$ elements of $S$ is called a \textit{reduced expression}. By convention, $\ell(1)=0$. 

More formally, a reduced expression of $w$, $w=\sigma_1 \sigma_2\cdots \sigma_{\ell(w)}$, can be viewed as an ordered $\ell(w)-$tuple, $(\sigma_1, \sigma_2,\ldots ,\sigma_{\ell(w)})$.

\begin{Example1}\rm Let $(S_4,S)$ be the Coxeter system and let $w=(1,2,4)$ be a permutation in $S_4$. A reduced expression of $w$ is $w=(1,2)(2,3)(3,4)(2,3)$. In the notations from Example \ref{sym}, we have $w=s_1s_2s_3s_2$. Note that $s_1s_3s_2s_3$ is also a reduced expression for $w$. We have that $\ell(w)=4$.
\end{Example1}

Next, we recall some properties of the length function.

\begin{Proposition1}[\cite{Bo}] Let $(W,S)$ be a Coxeter system. For all $u,w\in W$ the following hold:
\begin{enumerate}
	\item $\ell(w)=1$ if and only if $w\in S$.
	\item $\ell(w\sigma)=\ell(w)\pm1$ for all $\sigma\in S$ and $w\in W$.
	\item $\ell(\sigma w)=\ell(w)\pm1$ for all $\sigma\in S$ and $w\in W$,
	\item $|\ell(u)-\ell(w)|\leq \ell(uw)\leq\ell(u)+\ell(w)$.
	\item$\ell(w^{-1})=\ell(w)$.
\end{enumerate}
\end{Proposition1}

The \textit{"Exchange Property"} is a fundamental combinatorial property of a Coxeter system.

\begin{Theorem1}[Exchange Property \cite{Bo}] Suppose that $w=\sigma_1\sigma_2\cdots \sigma_k$ is a reduced expression, $\sigma_i\in S$ for $i=1,\ldots,k$ and $\sigma\in S$ a simple reflection. If $\ell(\sigma w)<\ell(w)$, then $$\sigma w=\sigma_1\cdots \widehat{\sigma_i}\cdots \sigma_k$$ for some $i\in [k]$, where we denote by $\widehat{\sigma_i}$ the deletion of the simple reflection $\sigma_i$ from the reduced expression.
\end{Theorem1}

An important consequence of the Exchange Property is the \textit{"Deletion Property"}:

\begin{Proposition1}[Deletion Property \cite{Bo}] If $w=\sigma_1\sigma_2\cdots \sigma_k$ is an element in $W$ and $\ell(w)<k$ then $w=\sigma_1\cdots \widehat{\sigma_i}\cdots\widehat{\sigma_j}\cdots \sigma_k$ for some $1\leq i<j\leq k$.
\end{Proposition1}

\begin{Corollary1} [\cite{Bo}]$\ $

\begin{itemize}
	\item[(i)] Any expression $w=\sigma_1\sigma_2\cdots \sigma_k$, $\ell(w)<k$, contains a reduced expression of $w$ as a subword, obtainable by deleting an even number of letters.
	\item[(ii)] Suppose $w=\sigma_1\sigma_2\cdots \sigma_k=\sigma'_1\sigma'_2\cdots \sigma'_k$ are two reduced expressions. Then, the set of letters appearing in the word $\sigma_1\sigma_2\cdots \sigma_k$ equals the set of letters appearing in $\sigma'_1\sigma'_2\cdots \sigma'_k$.
	\item[(iii)] $S$ is a minimal generating set for $W$; that is, no Coxeter generator can be expressed in terms of the others.
\end{itemize}
\end{Corollary1}

We denote by "$\prec$" the Bruhat order on the Coxeter system $(W,S)$. Recall that, if $w=\sigma_1\cdots \sigma_r\in W$ is a fixed but arbitrary reduced expression for $w$ and $w\,'$ is an element in $W$, then $w'\preceq w$ if and only if $w'$ can be obtained as a subexpression of this reduced expression (see, for instance, J.E. Humphreys \cite{Hu}).

\section{Subword complexes in Coxeter groups}
Subword complexes were introduced by A. Knutson and E. Miller \cite{KM} for the study of Schubert polynomials and combinatorics of determinantal ideals. Many properties of the subword complexes in Coxeter groups are obtained by using the Demazure algebra and the Demazure product.

Let $(W,S)$ be an arbitrary Coxeter system consisting of a Coxeter group $W$ and a set of simple reflections $S$ that generates $W$.

\begin{Definition1}\rm $\ $A \textit{word} $Q$ of size $n$ is an ordered sequence $Q=(\sigma_1,\ldots,\sigma_n)$ of elements of $S$. An ordered subsequence $P$ of $Q$ is called a \textit{subword} of $Q$. Let $\pi\in W$ be an element.
\begin{enumerate}
	\item $P$ \textit{represents} $\pi$ if the ordered product of the simple reflections in $P$ is a reduced expression for $\pi$.
	\item $P$ \textit{contains} $\pi$ if some subsequence of $P$ represents $\pi$.
\end{enumerate}
\end{Definition1}

\begin{Definition1}\rm The \it{subword complex}\rm $\ \Delta(Q,\pi)$ is the set of subwords $Q\setminus P$ such that $P$ contains $\pi$.
\end{Definition1}

\begin{Lemma1}[\cite{KnMi}] The subword complex $\Delta(Q,\pi)$ is a pure simplicial complex whose facets are the subwords $Q\setminus P$ such that $P\subseteq Q$ represents $\pi$.
\end{Lemma1}

\begin{Example1}\rm Let $(S_4,S)$ be the Coxeter system, and $Q$ the word of size $8$, $$Q=(s_1,\ s_2,\ s_1,\ s_3,\ s_1,\ s_2,\ s_3,\ s_1)$$
Let $\pi=(1,2,4)\in S_4$ with $\ell(\pi)=4$. The set of all the reduced expressions of $\pi$ is $$\{s_1s_2s_3s_2,\ s_1s_3s_2s_3,\ s_3s_1s_2s_3\}.$$ Let us denote $Q=(\sigma_1,\ \sigma_2,\ \sigma_3,\ \sigma_4,\ \sigma_5,\ \sigma_6,\ \sigma_7,\ \sigma_8).$ The set of all the subwords of $Q$ that represent $\pi$ is $$\{(\sigma_1,\sigma_2,\sigma_4,\sigma_6),\ (\sigma_1,\sigma_4,\sigma_6,\sigma_7),\  (\sigma_3,\sigma_4,\sigma_6,\sigma_7),\ (\sigma_4,\sigma_5,\sigma_6,\sigma_7)\}.$$ Therefore, the subword complex $\Delta=\Delta(Q,\pi)$ is the simplicial complex with the facets $$\mathcal{F}(\Delta)=\{\{\sigma_3,\sigma_5,\sigma_7,\sigma_8\},\ \{\sigma_2,\sigma_3,\sigma_5,\sigma_8\},\ \{\sigma_1,\sigma_2,\sigma_5,\sigma_8\},\ \{\sigma_1,\sigma_2,\sigma_3,\sigma_8\}\}.$$
\end{Example1}

\begin{Theorem1}[\cite{KM}]\label{vdshell} Subword complexes $\Delta(Q,\pi)$ are vertex-decomposable,\newline hence they are shellable.
\end{Theorem1}

A useful tool in the study of the subword complexes in Coxeter groups is the Demazure algebra. We recall, following E. Miller and A. Knutson \cite{KnMi}, the notion of Demazure algebra of an arbitrary Coxeter system $(W,S)$ over a commutative ring $R$.

\begin{Definition1}\rm$\ $ Let $R$ be a commutative ring and $\mathcal{D}$ be a free $R$--module with basis $\{e_{\pi}\mid\pi\in W\}$. Defining a multiplication on $\mathcal{D}$ by
\[e_{\pi}e_{\sigma}=
\left\{\begin {array}{cc}
			e_{\pi \sigma}, & \mbox{if}\ \ell(\pi \sigma)>\ell(\pi)\\
			e_{\pi}, & \mbox{if}\ \ell(\pi \sigma)<\ell(\pi),
	\end{array}\right.
\]  
for $\sigma\in S$ yields the \textit{Demazure algebra} of $(W,S)$ over $R$. It can be defined the \textit{Demazure product}, $\delta(Q)$, of the word $Q=(\sigma_1,\ldots,\sigma_n)$ by $e_{\sigma_1}\cdots e_{\sigma_n}=e_{\delta(Q)}$.
\end{Definition1}

Henceforth we write $Q\setminus \sigma_i$ for the word of size $n-1$ obtained from $Q$ by omitting $\sigma_i$, that is $Q\setminus \sigma_i=(\sigma_1,\ldots,\sigma_{i-1},\sigma_{i+1},\ldots,\sigma_n).$ We denote also by $"\succ"$ the Bruhat order on $W$.

\begin{Lemma1}[\cite{KnMi}]\label{pi} Let $P$ be a word in $W$ and let $\pi$ be an element in $W$.
\begin{itemize}
	\item[(a)] $\delta(P)\succeq\pi$ if and only if $P$ contains $\pi$.
	\item[(b)] If $\delta(P)=\pi$, then every subword of $P$ containing $\pi$ has the Demazure product $\pi$.
	\item[(c)] If $\delta(P)\succ\pi$, then $P$ contains a word $T$ representing an element $\tau\succ\pi$ satisfying $|T|=\ell(\tau)=\ell(\pi)+1$. 
\end{itemize}
\end{Lemma1}

\begin{Lemma1}[\cite{KnMi}]\label{m-1} Let $T$ be a word in $W$ and let $\pi$ be an element in $W$ such that $|T|=\ell(\pi)+1$.
\begin{itemize}
	\item[(a)] There are at most two elements $\sigma\in T$ such that $T\setminus \sigma$ represents $\pi$.
	\item[(b)] If $\delta(P)=\pi$, then there are two distinct $\sigma\in T$ such that $T\setminus \sigma$ represents $\pi$.
	\item[(c)] If $T$ represents $\tau\succ\pi$, then $T\setminus \sigma$ represents $\pi$ for exactly one $\sigma\in T$. 
\end{itemize}
\end{Lemma1}

The following theorem gives a complete description of the structure of the subword complexes in Coxeter groups.

\begin{Theorem1}[\cite{KnMi}]\label{sphere} The subword complex $\Delta(Q,\pi)$ is a sphere if $\delta(Q)=\pi$ and a ball otherwise.
\end{Theorem1}

Recall that, if we consider $k[x_1,\ldots,x_n]$ with the fine grading and if $I$ is a monomial ideal of $S$, then the Hilbert series of $I$ has the form
	\[H_{I}(t_1,\ldots,t_n)=\frac{\mathcal{K}_{I}(t_1,\ldots,t_n)}{\prod\limits_{i=1}^n (1-t_i)}
\]
and $\mathcal{K}_I(t_1,\ldots,t_n)$ is called the \textit{K--polynomial} or the \textit{Hilbert numerator}, that is $\mathcal{K}_I(t_1,\ldots, t_n) = \sum\limits_\twoline{P\subset\mathbb{Z}^n}{j\geq0}(-1)^j \beta_{j,P}(I) \cdot \mathbf{t}^P$ where $\mathbf{t}^P=\prod\limits_{p_i\in P}t_i^{p_i}$ and $\beta_{j,P}(I)$ is the $j$th Betti number of $I$ in degree $P$.

\begin{Lemma1}[\cite{KnMi}]\label{hilb} If $\Delta$ is the subword complex $\Delta(Q,\pi)$, then the Hilbert numerator of the Alexander dual is
\begin{eqnarray*}
  \mathcal{K}_{I_{\Delta^\vee}}( t_1,\ldots,t_n) &=& \sum_\twoline{P \subseteq Q}
  {\delta(P)=\pi} (-1)^{|P|-\ell(\pi)} \mathbf{t}^P.
\end{eqnarray*}
\end{Lemma1} 

\begin{Theorem1} [\cite{KnMi}]
If $\Delta$ is the subword complex $\Delta(Q,\pi)$, then the Hilbert series of the Stanley--Reisner ring
$k[\Delta]$ is
\begin{eqnarray*}
  H(k[\Delta]; t_1,\ldots,t_n) &=& \frac{\sum\limits_{\delta(P)=\pi}
  (-1)^{|P| - \ell(\pi)}(\mathbf{1}-\mathbf{t})^P}{\prod\limits_{i=1}^n(1-t_i)},
\end{eqnarray*}
where $(\mathbf{1}-\mathbf{t})^P = \prod\limits_{\sigma_i \in P}(1-t_i)$ and the sum is over the subwords $P \subseteq Q$.
\end{Theorem1}

%
%
%
%
%
%
%
%
\chapter{Lexsegment ideals}

In this chapter we show that any lexsegment ideal with a linear resolution has linear quotients with respect to a suitable order of the generators. For the completely lexsegment ideals with a linear resolution it will turn out that their decomposition function with respect to the suitable defined ordering is regular. Therefore, one may apply the procedure developed in J. Herzog and Y. Takayama \cite{HeTa} to get the explicit resolutions for this class of ideals.

In the last section of this chapter we study the depth and the dimension of lexsegment ideals. Our results show that one may compute these invariants just looking at the ends of the lexsegment. As an application, we characterize the Cohen--Macaulay lexsegment ideals.

\section{Completely lexsegment ideals with a linear resolution}\label{Section1}
Let $S=k[x_1,\ldots, x_n]$ be the polynomial ring in $n$ variables over a field $k$. We order lexicographically the monomials of $S$ such that $x_1>x_2>\ldots> x_n$. Let $d\geq 2$ be an integer and $\mathcal{M}_d$ be the set of all the monomials of $S$ of degree $d$.

\begin{Theorem1}\label{colex}
Let $u=x_1^{a_1}\cdots x_n^{a_n},$ with  $a_1>0,$ and $ v=x_1^{b_1}\cdots x_n^{b_n}$ be monomials of degree $d$ with $u\geq_{lex} v,$ such that $I=({\mathcal L}(u,v))$ is a completely lexsegment ideal. Then $I$ has a linear resolution if and only if $I$ has linear quotients.
\end{Theorem1}

\begin{proof}
We have to prove that if $I$ has a linear resolution then $I$ has linear quotients, since the other implication is known (Proposition \ref{linquotlinres}). By Theorem \ref{completelylex}, since $I$ has a linear resolution, one of the conditions (a), (b), (c) holds.

We define on the set of the monomials of  degree $d$ from $S$ the following total order: For $$w=x_1^{\alpha_1}\cdots x_n^{\alpha_n},\ w^{\prime}=x_1^{\beta_1}\cdots x_n^{\beta_n},$$ we set $$w\prec w^{\prime}\text{ if }\alpha_1 < \beta_1 \text{ or }\alpha_1=\beta_1 \text{ and } w>_{lex} w^{\prime}.$$ Let $$\mathcal{L}(u,v)=\{w_1,\ldots,w_r\}, \text{ where }w_1\prec w_2\prec\ldots\prec w_r.$$ We shall prove that $I=(\mathcal{L}(u,v))$ has linear quotients with respect to this ordering of the generators. 

Assume that $u,v$ satisfy the condition (a) and $a<d$ (the case $a=d$ is trivial). Then $I$ is isomorphic as $S$--module with the ideal generated by the final lexsegment ${\mathcal L}^f(x_2^{d-a}) \subset S$ and the ordering $\prec$ of its minimal generators coincides with the lexicographical ordering $>_{lex}.$ The ideal $(\mathcal{L}^f(x_2^{d-a}))\cap k[x_2,\ldots,x_n]$ is the initial ideal in $k[x_2,\ldots,x_n]$ defined by $x_n^{d-a},$ which has linear quotients with respect to $>_{lex}.$ Hence $I$ has linear quotients with respect to $\prec$ since it is obvious that the extension in the ring $k[x_1,\ldots, x_n]$ of a monomial  ideal with linear quotients in $k[x_2,\ldots,x_n]$ has linear quotients, too.

Next we assume that $u,v$ satisfy either the condition (b) or the condition (c).

By definition, $I$ has linear quotients with respect to the monomial generators $w_1,\ldots,w_r$ if the colon ideals $(w_1,\ldots,w_{i-1}):w_i$ are generated by variables for all $i\geq 2,$ that is for all $j < i$ there exist an integer $1\leq k <i$ and an integer $l\in [n]$ such that $w_k/\gcd(w_k,w_i)=x_l \text{ and }x_l\text{ divides } w_j/\gcd(w_j,w_i).$ 

In other words, for any $w_j\prec w_i, w_j,w_i\in \mathcal{L}(u,v),$ we have to find a monomial $w^{\prime}\in \mathcal{L}(u,v)$ such that 
\begin{eqnarray}
w^{\prime}\prec w_i,\ \frac{w^{\prime}}{\gcd(w^{\prime},w_i)}=x_l, \text{ for some }l\in [n]  \text{, and }x_l\text{ divides } \frac{w_j}{\gcd(w_j,w_i)}. \eqname{$*$}\label{*}
\end{eqnarray}

Let us fix  $w_i=x_1^{\alpha_{1}}\cdots x_n^{\alpha_{n}}$ and $w_j=x_1^{\beta_{1}}\cdots x_n^{\beta_{n}}, \ w_i,w_j\in \mathcal{L}(u,v),$ such that $w_j\prec w_i.$ By the definition of the  ordering $\prec,$ we must have one of the relations $$\beta_{1} < \alpha_{1} \text { or }\beta_{1} = \alpha_{1} \text{ and } w_j>_{lex} w_i.$$ We consider them in turn.

{\textit{Case $1$}}. Let $\beta_1<\alpha_1$. One may find an integer $l$, $2\leq l\leq n$, such that $\alpha_s\geq \beta_s$ for all $s<l$ and $\alpha_l<\beta_l$ since, otherwise, $\deg(w_i)>\deg(w_j)=d$ which is impossible. We obviously have $\max(w_j)\geq l$. If $l\geq \max(w_i),$ one may take $\bar{w}=x_lw_i/x_1$ which satisfies the condition (\ref{*}) since the inequalities $\bar{w}\prec w_i,$ $\bar{w}\leq_{lex}w_i\leq_{lex} u$ hold, and we will show that $\bar{w}\geq_{lex}w_j$. This will imply that $\bar{w}\geq_{lex}v$, hence $\bar{w}\in\mathcal{L}(u,v)$.

The inequality $\bar{w}\geq_{lex}w_j$ is obviously fulfilled if $\alpha_1-1>\beta_1$ or if $\alpha_1-1=\beta_1$ and at least one of the inequalities $\alpha_s\geq\beta_s$ for $2\leq s<l$ is strict. If $\alpha_1-1=\beta_1$ and $\alpha_s=\beta_s$ for all $s<l$, comparing the degrees of $w_i$ and $w_j$, it results $$d=\alpha_1+\ldots+\alpha_l=\beta_1+1+\beta_2+\ldots+\beta_{l-1}+\alpha_l<(\beta_1+1)+\beta_2+\ldots+\beta_l.$$ It follows that $d\geq\beta_1+\beta_2+\ldots+\beta_l>d-1,$ that is $\beta_1+\beta_2+\ldots+\beta_l=d.$ This implies that $l=\max(w_j)$ and $\beta_l=\alpha_l+1$, that is $$\bar{w}=\frac{x_lw_i}{x_1}=x_1^{\alpha_1-1}x_2^{\alpha_2}\cdots x_l^{\alpha_l+1}=x_1^{\beta_1}\cdots x_l^{\beta_l}=w_j.$$ 

From now on, in the Case $1$, we may assume that $l<\max(w_i)$. We will show that at least one of the following monomials $$w'=\frac{x_lw_i}{x_{\max(w_i)}},\ w''=\frac{x_lw_i}{x_1}$$ belongs to $\mathcal{L}(u,v)$. It is clear that both monomials are strictly less than $w_i$ with respect to the ordering $\prec.$ Therefore one of the monomials $w',\ w''$ will satisfy the condition (\ref{*}).

The following inequalities are fulfilled: $$w'>_{lex}w_i\geq_{lex}v,\ \text{and}$$ $$w''<_{lex}w_i\leq_{lex}u.$$

Let us assume, by contradiction, that $w'>_{lex}u$ and $w''<_{lex}v$. Comparing the exponents of the variable $x_1$, we obtain $a_1-1\leq\alpha_1-1\leq b_1$. Since the ideal generated by $\mathcal{L}(u,v)$ has a linear resolution, we must have $b_1=a_1-1$. Let $z$ be the greatest monomial of degree $d$ such that $z<_{lex}v$. Then, by our assumption on $w''$, we also have the inequality $w''\leq_{lex}z$.

Now we need the following

\begin{Lemma1}\label{lema1} Let $m=x_1^{a_1}\cdots x_n^{a_n},\ m'=x_1^{b_1}\cdots x_n^{b_n}$ be two monomials of degree $d$. If $m\leq_{lex}m'$ then $$\frac{m}{x_{\max(m)}}\leq_{lex}\frac{m'}{x_{\max(m')}}.$$
\end{Lemma1}
\begin{proof} Let $m <_{lex}m'$. Then there exists $s\geq1$ such that $a_1=b_1,\ldots, a_{s-1}=b_{s-1}$ and $a_s<b_s$. It is clear that $\max(m')\geq s$. Comparing the degrees of $m$ and $m'$, we get  $\max(m)>s.$ 

If $\max(m)>s$ and $\max(m')>s$, the required inequality is obvious.

Let $\max(m)>s$ and $\max(m')=s$. Let us suppose, by contradiction, that $$\frac{m}{x_{\max(m)}}>_{lex}\frac{m'}{x_{\max(m')}}=\frac{m'}{x_s}.$$ This implies that $a_s\geq b_s-1$, and, since $a_s<b_s$, we get $a_s=b_s-1$.
Looking at the degree of $m'$ we obtain $d=b_1+b_2+\ldots+b_s=a_1+a_2+\ldots+ a_{s-1}+a_s+1,$ that is $a_1+\ldots+a_s=d-1$. It follows that $a_{\max(m)}=1$ and $$\frac{m}{x_{\max(m)}}=x_1^{a_1}\cdots x_s^{a_s}=x_1^{b_1}\cdots x_{s-1}^{b_{s-1}}x_s^{b_s-1}=\frac{m'}{x_{\max(m')}},$$ contradiction.
\end{proof}
Going back to the proof of our theorem, 
we apply the above lemma for the monomials $w''$ and $z$ and we obtain  $$\frac{w''}{x_{\max(w'')}}\leq_{lex}\frac{z}{x_{\max(z)}},$$ which implies that $$\frac{x_1w''}{x_{\max(w'')}}\leq_{lex}\frac{x_1z}{x_{\max(z)}}.$$
By using  condition (c) in the Theorem \ref{completelylex} it follows that $x_1w''/x_{\max(w'')}\leq_{lex}u.$ On the other hand, $$\frac{x_1w''}{x_{\max(w'')}}=\frac{x_1x_lw_i}{x_1x_{\max(w_i)}}=\frac{x_lw_i}{x_{\max(w_i)}}=w'.$$ Therefore, it results $w'\leq_{lex}u$, which contradicts our assumption on $w'$.

Consequently, we have $w'\leq_{lex}u$ or $w''\geq_{lex}v$, which proves that at least one of the monomials $w',\ w''$ belongs to $\mathcal{L}(u,v)$.

{\textit{Case $2$.}} Let $\beta_1=\alpha_1$ and $w_j>_{lex}w_i$. Then there exists $l$, $2\leq l\leq n$, such that $\alpha_s=\beta_s$, for all $s<l$ and $\alpha_l<\beta_l$. If $\max(w_i)\leq l$, then, looking at the degrees of $w_i$ and $w_j$, we get $d=\alpha_1+\alpha_2+\ldots+\alpha_l<\beta_1+\beta_2+\ldots+\beta_l,$ contradiction. Therefore, $l<\max(w_i)$. We proceed in a similar way as in the previous case. Namely, exactly as in the Case $1,$ it results that at least one of the following two monomials $w'=x_lw_i/x_{\max(w_i)},\ w''=x_lw_i/x_1$ belongs to $\mathcal{L}(u,v)$. It is clear that both monomials are strictly less than $w_i$ with respect to the order $\prec$.  
\end{proof}

\begin{Example1}\rm 
Let $S=k[x_1,x_2,x_3]$. We consider the monomials: $u=x_1x_2x_3$ and $v=x_2x_3^2$, $u>_{lex} v,$ and let $I$ be the monomial ideal generated by $\mathcal{L}(u,v)$. The minimal system of generators of the ideal $I$ is
	$$G(I)=\mathcal{L}(u,v)=\{x_1x_2x_3,\ x_1x_3^2,\ x_2^3,\ x_2^2x_3,\ x_2x_3^2\}.$$

Since $I$ verifies the condition (c) in Theorem \ref{completelylex}, it follows that $I$ is a completely lexsegment ideal with a linear resolution. We denote the monomials from $G(I)$ as follows: $u_1=x_1x_2x_3,\ u_2=x_1x_3^2,\ u_3=x_2^3,\ u_4=x_2^2x_3,\ u_5=x_2x_3^2$, so $u_1>_{lex}u_2>_{lex}\ldots>_{lex}u_5.$ The colon ideal $(u_1,u_2):u_3=(x_1x_3)$ is not generated by a subset of $\{x_1,x_2,x_3\}$. This shows that $I$ is not with linear quotients with respect to lexicographical order.

We consider now the  order $\prec$ and check by direct computation that $I$ has a linear quotients.  We label the monomials from $G(I)$ as follows: $u_1=x_2^3,\ u_2=x_2^2x_3,\ u_3=x_2x_3^2,\ u_4=x_1x_2x_3,\ u_5=x_1x_3^2$, so $u_1\prec u_2\prec\ldots\prec u_5$. Then $(u_1):u_2=(x_2),\ (u_1,u_2):u_3=(x_2),\ (u_1,u_2,u_3):u_4=(x_2,x_3),$ $(u_1,u_2,u_3,u_4):u_5=(x_2)$. 
\end{Example1}

We further study the decomposition function of a completely lexsegment ideal with a linear resolution. The decomposition function of a monomial ideal was introduced by J. Herzog and Y. Takayama in \cite{HeTa}. We show in the sequel that completely lexsegment ideals which have linear quotients with respect to $\prec$ have also regular decomposition functions.

In order to do this, we need some preparatory notations and results.

For an arbitrary lexsegment $\mathcal{L}(u,v)$ with the elements ordered by $\prec$, we denote by $I_{\prec w}$, the ideal generated by all the monomials $z\in\mathcal{L}(u,v)$ with $z\prec w$. $I_{\preceq}w$ will be the ideal generated by all the monomials $z\in\mathcal{L}(u,v)$ with $z\preceq w$.

\begin{Lemma1}\label{1nset} Let $I=(\mathcal{L}(u,v))$ be a  lexsegment ideal which has linear quotients with respect to the order $\prec$ of the generators. Then, for any $w\in\mathcal{L}(u,v)$, $1\notin\set(w)$.
\end{Lemma1}

\begin{proof} Let us assume that $1\in\set(w)$, that is $x_1w\in I_{\prec w}$. It follows that there exist $w'\in\mathcal{L}(u,v)$, $w'\prec w$, and a variable $x_j$ such that $x_1w=x_jw'$. Obviously, we have $j\geq2$. But this equality shows that $\nu_1(w')>\nu_1(w)$, which is impossible since $w'\prec w$.
\end{proof}

\begin{Lemma1}\label{desc} Let $I=(\mathcal{L}(u,v))$ be a completely lexsegment ideal which has linear quotients with respect to the ordering $\prec$ of the generators and let $g:M(I)\rightarrow G(I)$ be the decomposition function of $I$ with respect to the ordering $\prec$. If $w\in\mathcal{L}(u,v)$ and $s\in\set(w)$, then
	\[g(x_sw)=
\left\{\begin {array}{cc}
			\frac{x_s w}{x_1}, & \mbox{if}\ x_sw\geq_{lex} x_1v,\\
			&\\
			\frac{x_s w}{x_{\max(w)}}, & \mbox{if}\ x_sw<_{lex} x_1v.
	\end{array}\right. \]
\end{Lemma1}
\begin{proof} Let $u=x_1^{a_1}\cdots x_n^{a_n},\ v=x_1^{b_1}\cdots x_n^{b_n},\ a_1>0,$ and $w=x_1^{\alpha_1}\cdots x_n^{\alpha_n}$. 

In the first place we consider $$x_sw\geq_{lex}x_1v.$$ Since, by Lemma \ref{1nset}, we have $s\geq 2$, the above inequality shows that $\nu_1(w)\geq1$. We have to show that $g(x_sw)=x_s w/x_1$, that is $$\frac{x_s w}{x_1}=\mbox{min} _{\prec}\{w'\in\mathcal{L}(u,v)\ |\ x_sw\in I_{\preceq w'}\}.$$ It is clear that $v\leq_{lex} x_sw/x_1<_{lex}w\leq_{lex}u$, hence $x_sw/x_1\in\mathcal{L}(u,v)$. Let $w'\in\mathcal{L}(u,v)$ such that $x_sw\in I_{\preceq w'}$. We have to show that $x_sw/x_1\preceq w'$. Let $w''\in\mathcal{L}(u,v),\ w''\preceq w'$ such that $x_sw=w''x_j$, for some variable $x_j$. Then
$w''=x_sw/x_j\succeq x_sw/x_1
$
by the definition of our ordering $\prec$. This implies that $w'\succeq x_sw/x_1$.

Now we have to consider the second inequality 
\begin{eqnarray}x_sw<_{lex}x_1v.\label{star}
\end{eqnarray}

Since $s\in\set(w)$, we have $x_sw\in I_{\prec w}$, that is there exists $w'\in\mathcal{L}(u,v)$, $w'\prec w$, and a variable $x_j,\ j\neq s$, such that 
\begin{eqnarray}x_sw=x_jw'.\label{1}
\end{eqnarray}

If $j=1$, then $x_sw=x_1w'\geq_{lex}x_1v$, contradiction. Hence $j\geq2$. We also note that $x_j|w$ since $j\neq s$, thus $j\leq\max(w)$. The following inequalities hold: 
\begin{eqnarray}
\frac{x_s w}{x_{\max(w)}}\geq_{lex} \frac{x_s w}{x_j}= w'\geq_{lex} v. \label{4}
\end{eqnarray}

If $\nu_1(w)<a_1$, we obviously get $x_sw/x_{\max(w)}\leq_{lex}u$. Let $\nu_1(w)=a_1$. From the inequality (\ref{star}) we obtain  $a_1\leq b_1+1$.  

If $a_1=b_1$, then $u=x_1^{a_1}x_2^{d-a_1}$ and $v=x_1^{a_1}x_n^{d-a_1}$ by Theorem \ref{completelylex}. Since $w\leq_{lex} u$, by using Lemma \ref{lema1}, we have
$$\frac{x_sw}{x_{\max(w)}}\leq_{lex}\frac{x_su}{x_{\max(u)}}=\frac{x_su}{x_2}\leq_{lex} u,
$$ the last inequality being true by Lemma \ref{1nset}. Therefore, $x_s w/x_{\max(w)}\in \mathcal{L}(u,v).$

If $a_1=b_1+1$ then the condition (c) in Theorem \ref{completelylex} holds. Let $z$ be the greatest monomial with respect to the lexicographical order such that $z<_{lex}v$. Since $x_sw/x_1<_{lex}v$ by hypothesis, we also have $x_sw/x_1\leq_{lex}z$. By Lemma \ref{lema1} we obtain 
$$\frac{x_sw}{x_1x_{\max\left(\frac{x_sw}{x_1}\right)}}\leq_{lex}\frac{z}{x_{\max(z)}}.
$$
Next we apply the condition (c) from Theorem \ref{completelylex} and get the following inequa- lities:

\begin{eqnarray}
x_1\frac{x_sw}{x_1x_{\max\left(\frac{x_sw}{x_1}\right)}}\leq_{lex}x_1\frac{z}{x_{\max(z)}}\leq_{lex}u.\label{2}
\end{eqnarray}
From the equality (\ref{1}) we have $w'=x_sw/x_j$.  As $j\neq1,\ \nu_1(w')=\nu_1(w),$ and the  inequality $w'\prec w$   gives $w'>_{lex} w,$ that is $x_s w/x_j>_{lex}w,$ which implies that $x_s>_{lex} x_j$. This shows that  $s<j\leq\max(w)$. Now looking at the inequalities (\ref{2}), we have
\begin{eqnarray}\frac{x_sw}{x_{\max(w)}}\leq_{lex}u.\label{3}
\end{eqnarray}
From (\ref{3}) and (\ref{4}) we obtain  $x_sw/x_{\max(w)}\in\mathcal{L}(u,v)$. 

It remains to show that $x_s w/x_{\max(w)}=\min_{\prec}\{w'\in\mathcal{L}(u,v)\mid x_sw\in I_{\preceq w'}\}.$ Let $\tilde{w}=\min_{\prec}\{w'\in\mathcal{L}(u,v)\mid x_sw\in I_{\preceq w'}\}.$ We obviously have $\tilde{w}\preceq x_s w/x_{\max(w)}\prec w.$ By the choice of $\tilde{w}$ we have $$x_sw=x_t\tilde{w},$$ for some variable $x_t.$ 

If $t=s$ we get $w=\tilde{w}$ which is impossible since $\tilde{w}\prec w.$ Therefore $t\neq s.$
Then $x_t| w,$  so $t\leq\max(w)$. It follows that 
$$\tilde{w}=\frac{x_sw}{x_t}\leq_{lex}\frac{x_sw}{x_{\max(w)}}.$$
 If $t=1$ we have $x_1\tilde{w}=x_sw <_{lex}x_1 v,$ which implies that $\tilde{w}<_{lex}v,$ contradiction. Therefore $t\neq 1$ and, moreover, $\tilde{w}\succeq x_sw/x_{\max(w)}$, the  inequality being true by the definition of the ordering $\prec$. 
 This yields $\tilde{w}=x_sw/x_{\max(w)}.$ Therefore we have proved that $$\frac{x_sw}{x_{\max(w)}}=g(x_sw).$$
\end{proof}

After this preparation, we  prove the following 

\begin{Theorem1} Let $u=x_1^{a_1}\cdots x_n^{a_n},\ v=x_1^{b_1}\cdots x_n^{b_n},\ u,v\in\mathcal{M}_d$, with $u\geq_{lex}v$, and $I=(\mathcal{L}(u,v))$ be a completely lexsegment ideal which has a linear resolution. Then the decomposition function $g:M(I)\rightarrow G(I)$ associated to the ordering  $\prec$ of the generators from $G(I)$ is regular.
\end{Theorem1}

\begin{proof} Let $w\in\mathcal{L}(u,v)$ and $s\in\set(w)$. We have to show that $\set(g(x_sw))\subset\set(w)$.

Let $t\in\set(g(x_sw))$. In order to prove that $t\in\set(w)$, that is $x_tw\in I_{\prec w}$, we will consider the following two cases:

{\textit{Case 1}.} Let $x_sw\geq_{lex}x_1v$. By Lemma \ref{desc}, $g(x_sw)=x_sw/x_1$. Since $t\in\set(g(x_sw))$, we have $$\frac{x_tx_sw}{x_1}\in I_{\prec\frac{x_sw}{x_1}},$$ so there exists $w'\prec x_sw/x_1$, $w'\in\mathcal{L}(u,v)$, and a variable $x_j$, such that
$$\frac{x_t x_sw}{x_1}=x_jw',
$$
that is
	\begin{eqnarray}x_tx_sw=x_1x_jw'.\label{6}
\end{eqnarray}
By Lemma \ref{1nset}, $s,t\neq1$ and, since $w'\prec\ x_sw/x_1$, we have $j\neq t$. Note also that $w'\prec w$ since $\nu_1(w')<\nu_1(w)$. If $j=s$ then $x_tw=x_1w'\in I_{\prec w}$ and $t\in\set(w)$.

Now let $j\neq s$. If $j=1$, we have 
$x_tx_sw=x_1^2w',
$
which implies that $\nu_1(w')=\nu_1(w)-2$. The following inequalities hold:
$$v<_{lex}\frac{x_1w'}{x_s}<_{lex}w\leq_{lex}u,
$$
the first one being true since $v\leq_{lex}w'$, so $\nu_1(v)\leq\nu_1(w')$. These inequalities show that $x_1w'/x_s\in\mathcal{L}(u,v)$. But we also have $x_1w'/x_s\prec w$, hence $x_1w'/x_s\in I_{\prec w}$.

To finish this case we only need  to treat the case $j\neq1,\ j\neq s$. We are going to show that at least one of the monomials $x_1w'/x_s$ or $x_jw'/x_s$ belongs to $I_{\prec w}.$ In any case this will lead to the conclusion that $x_t w\in I_{\prec w}$ by using (\ref{6}).

From the equality (\ref{6}), we have $x_j|w$, hence $j\leq\max(w)$, and $\nu_1(w')=\nu_1(w)-1$. Since $w'\prec\ x_sw/x_1$ and $\nu_1(w')=\nu_1(w)-1=\nu_1(x_sw/x_1)$, we get 
\begin{eqnarray}w'>_{lex}\frac{x_sw}{x_1},\label{7} 
\end{eqnarray} which gives
	\[\frac{x_1w'}{x_s}>_{lex}v.
\]
If the inequality 
\begin{eqnarray}\frac{x_1w'}{x_s}\leq_{lex}u\label{8}
\end{eqnarray}
holds, then we get $x_1w'/x_s\in\mathcal{L}(u,v)$. We also note that $\nu_1(x_1w'/x_s)=\nu_1(w)$ and $x_1w'/x_s>_{lex}w$ (by (\ref{7})). Therefore $x_1w'/x_s\prec w$ and we may write
$$x_tw=x_j\frac{x_1w'}{x_s}\in I_{\prec w}.
$$
This implies that $t\in\set(w)$.

Now we look at the monomial $x_jw'/x_s$ for which we have $\nu_1(x_jw'/x_s)=\nu_1(w')<\nu_1(w)$, so
$x_jw'/x_s<_{lex}w\leq_{lex}u.
$
If the inequality 
\begin{eqnarray}\frac{x_jw'}{x_s}\geq_{lex}v\label{9}
\end{eqnarray}
holds, we obtain  $x_jw'/x_s\in\mathcal{L}(u,v)$. Obviously we have $x_jw'/x_s\prec w$. By using (\ref{6}), we may write 
$$x_tw=x_1\frac{x_jw'}{x_s}\in I_{\prec w},
$$
which shows that $t\in\set(w)$.

To finish the proof in the Case $1$ we have to consider the situation when both inequalities (\ref{8}) and (\ref{9}) fail. Hence, let
	\[\frac{x_1w'}{x_s}>_{lex}u\ \mbox{and}\ \frac{x_jw'}{x_s}<_{lex}v.
\]
We will show that this inequalities cannot hold simultaneously. Comparing the exponents of $x_1$ in the monomials involved in the above inequalities, we obtain  $\nu_1(w')=b_1\geq a_1-1$. Since, by hypothesis, $x_sw>_{lex} x_1v,$ we have $\nu_1(w)>b_1.$ On the other hand, $w\leq_{lex} u$ implies that $\nu_1(w)\leq a_1.$ So $b_1=a_1-1$ and $\mathcal{L}(u,v)$ satisfies the condition (c) in Theorem \ref{completelylex}. Let, as usually, $z$ be the largest monomial with respect to the lexicographical order such that $z<_{lex}v$.

Since $x_jw'/x_s<_{lex}v$, we have $x_jw'/x_s\leq_{lex}z$. By Lemma \ref{lema1} and using the condition $x_1z/x_{\max(z)}\leq_{lex}u$, we obtain:
$$\frac{x_1x_jw'}{x_sx_{\max\left(\frac{x_jw'}{x_s}\right)}}\leq_{lex}u.
$$
But our  assumption was that
$$u<_{lex}\frac{x_1w'}{x_s}.
$$
Therefore, combining the last two inequalities, after cancellation, one obtains that
$$x_j<_{lex}x_{\max\left(\frac{x_jw'}{x_s}\right)}=x_{\max\left(\frac{x_tw}{x_1}\right)}=x_{\max(x_tw)}.
$$
This leads to the inequality $j>\max(x_tw)$ and, since $j\leq\max(w)$, we get $\max(w)>\max(x_tw)$, which is impossible.

{\textit{Case 2}.} Let $x_sw<_{lex}x_1v$. Then $g(x_sw)=x_sw/x_{\max(w)}$. In particular, we have $x_sw/x_{\max(w)}\prec w$. Indeed, since $s\in \set(w),$ we have  $x_sw\in I_{\prec w}$, that is there exists $w'\in \mathcal{L}(u,v), w'\prec w,$ such that $x_sw\in I_{\preceq w'}.$ By the definition of the decomposition function we have $g(x_sw)\preceq w'$ and next we get  $g(x_sw)\prec w.$ Since $\nu_1(x_sw/x_{\max(w)})=\nu_1(w)$, the above inequality implies that $$\frac{x_sw}{x_{\max(w)}}>_{lex}w,$$ that is $x_s>_{lex}x_{\max(w)}$ which means that $s<\max(w)$.

As $t\in\set(g(x_sw))$, there exist $w'\prec\ x_sw/x_{\max(w)}$, $w'\in\mathcal{L}(u,v)$, and a variable $x_j$, such that
	\[\frac{x_tx_sw}{x_{\max(w)}}=x_jw',
\]
that is
	\begin{eqnarray}x_tx_sw=x_jx_{\max(w)}w'.\label{10}
\end{eqnarray}
As in the previous case, we would like to show that one of the monomials $x_{\max(w)}w'/x_s$ or $x_jw'/x_s$ belongs to $\mathcal{L}(u,v)$ and it is strictly less than $w$ with respect to $\prec$. In this way we  obtain  $x_tw\in I_{\prec w}$ and $t\in\set(w)$.

We begin our proof noticing that $s,t\neq1$, by Lemma \ref{1nset}. The equality $j=t$ is impossible since $w'\neq\ x_sw/x_{\max(w)}$. If $j=s$, then $x_tw=w'x_{\max(w)}\in I_{\preceq w'}$. But $w'\prec\ x_sw/x_{\max(w)}\prec w$, hence $x_tw\in I_{\prec w}$. 

Let $j\neq s,t$. From the equality (\ref{10}) we have $x_j|w$, so $j\leq\max(w)$. We firstly consider $j=1$. Then the equality (\ref{10}) becomes
	\begin{eqnarray}x_tx_sw=x_1x_{\max(w)}w'.\label{11}
\end{eqnarray}
Since $s< \max(w)$, we have
$$\frac{x_{\max(w)}w'}{x_s}<_{lex}w'\leq_{lex}u.
$$
If the inequality $x_{\max(w)}w'/x_s$ $\geq_{lex}v$ holds too, then $x_{\max(w)}w'/x_s\in\mathcal{L}(u,v)$ and, as $\nu_1(w')<\nu_1(w)$, it follows that $x_{\max(w)}w'/x_s\prec w$. From (\ref{11}), we have $$x_tw=x_1\frac{x_{\max(w)}w'}{x_s}\in I_{\prec w},$$ hence $t\in\set(w)$. 

From the inequality $x_sw<_{lex}x_1v$, we get
	\[x_sw<_{lex}x_1w',
\]
so
	\[\frac{x_1w'}{x_s}>_{lex}w.
\]
Let us assume that $x_1w'/x_s\leq_{lex}u$. Since $\nu_1(x_1w'/x_s)=\nu_1(w)$, by using the definition of the ordering $\prec$ we get  $x_1w'/x_s\in I_{\prec w}$. Then we may write $$x_tw=x_{\max(w)}\frac{x_1w'}{x_s}\in I_{\prec w}.$$

It remains to consider that 
$x_{\max(w)}w'/x_s<_{lex}v\ \mbox{and}\ x_1w'/x_s>_{lex}u.
$
Proceeding as in the case 1 we show that we reach a contradiction and this ends the proof for $j=1$.
We only need to notice that we have to consider $b_1\leq a_1-1.$ Indeed, we can not have $b_1=a_1$ since one may find in $\mathcal{L}(u,v)$ at least two monomials, namely $w$ and $w',$ with $\nu_1(w')<\nu_1(w).$ 

Finally, let $j\neq1$. Recall that in the equality (\ref{10})  we have $j\neq 1,t,s $ and $s<\max(w)$. From (\ref{10}) we obtain  $\nu_1(w)=\nu_1(w')$.  Since $w'\prec \ x_s w/x_{\max(w)}$, we have $w'>_{lex} x_sw/x_{\max(w)}$, that is 
\begin{eqnarray}w'x_{\max(w)}>_{lex}x_sw.\label{12}
\end{eqnarray}

Replacing $w'x_{\max(w)}$ by $x_tx_sw/x_j$ in (\ref{12}), we get $x_t>_{lex}x_j$, which means $t<j$. It follows that:
$$\frac{x_{\max(w)}w'}{x_s}=\frac{x_tw}{x_j}>_{lex} w\geq_{lex}v.
$$
Since $s<\max(w)$, as in the proof for $j=1$, we have $x_{\max(w)}w'/x_s\leq_{lex}u$. Therefore $x_{\max(w)}w'/x_s\in\mathcal{L}(u,v)$. In addition, from (\ref{12}), $x_{\max(w)}w'/x_s>_{lex}w$ and $\nu_1(x_{\max(w)}w'/x_s)=\nu_1(w)$, so $x_{\max(w)}w'/x_s\prec w$. In other words, we have got that
$$x_tw=x_j\frac{x_{\max(w)}w'}{x_s}\in I_{\prec w}
$$
and $t\in\set(w)$.
\end{proof}

The general problem of determining the resolution of arbitrary lexsegment ideals is not completely solved. 
In our specific context, one may apply the results from the first chapter, Section \ref{mapping}. We get the following result:

\begin{Theorem1} Let $I=(\mathcal{L}(u,v))\subset S$ be a completely lexsegment ideal with linear quotients with respect to $\prec$ and $\mathbb{F}_{\bullet}$ be the graded minimal free resolution of $S/I$. Then the chain map of $\mathbb{F}_{\bullet}$ is given by
\[\partial(f(\sigma; w)) = -\sum_{s\in\sigma}(-1)^{\alpha(\sigma;s)}x_sf(\sigma\setminus s;w)+\sum_{\stackrel{s\in\sigma:}{x_sw\geq_{lex}}x_1v}(-1)^{\alpha(\sigma;s)}x_1f\left(\sigma\setminus s;\frac{x_sw}{x_1}\right)+\]\[+\sum_{\stackrel{s\in\sigma:}{x_sw<_{lex}}x_1v}(-1)^{\alpha(\sigma;s)}x_{\max(w)}f\left(\sigma\setminus s;\frac{x_sw}{x_{\max(w)}}\right),\ 
\mbox{if}\ \sigma\neq\emptyset,\] and
\[\partial(f(\emptyset; w)) = w\ \mbox{otherwise}.\] For convenience we set $f(\sigma;w)=0$ if $\sigma\nsubseteq\set(w)$.
\end{Theorem1}

\begin{Example1}\rm
Let $u=x_1^2x_2$ and $v=x_2^3$  be monomials in the polynomial ring $S=k[x_1,x_2,x_3]$. Then $$\mathcal{L}(u,v)=\{x_2^3,\ x_1x_2^2,\ x_1x_2x_3,\ x_1x_3^2,\ x_1^2x_2\}.$$ The ideal $I=(\mathcal{L}(u,v))$ is a completely lexsegment ideal with linear quotients with respect to this ordering of the generators. We denote $u_1=x_2^3,\ u_2=x_1x_2^2,\ u_3=x_1x_2x_3,\ u_4=x_1x_3^2,\ u_5=x_1^2x_2$. We have  $\set(u_1)=\emptyset,\ \set(u_2)=\{2\},\ \set(u_3)=\{2\},\ \set(u_4)=\{2\},\ \set(u_5)=\{2,3\}$. Let $\mathbb{F}_{\bullet}$ be the minimal graded free resolution of $S/I$.

Since $\max\{|\set(w)|\mid w\in\mathcal{L}(u,v)\}=2$, we have $F_i=0$, for all $i\geq4$.

A basis for the $S-$module $F_1$ is $\{f(\emptyset;u_1),\ f(\emptyset;u_2),\ f(\emptyset;u_3),\ f(\emptyset;u_4),\ f(\emptyset;u_5)\}$.
	
A basis for the $S-$module $F_2$ is $$\{f(\{2\};u_2),\ f(\{2\};u_3),\ f(\{2\};u_4),\ f(\{2\};u_5),\ f(\{3\};u_5)\}.$$

A basis for the $S-$module $F_3$ is $\{f(\{2,3\};u_5)\}$.

We have the minimal graded free resolution $\mathbb{F}_{\bullet}$:
	\[0\longrightarrow S(-5)\stackrel{\partial_2}{\longrightarrow} S(-4)^5\stackrel{\partial_1}{\longrightarrow} S(-3)^5\stackrel{\partial_0}{\longrightarrow} S\longrightarrow S/I\longrightarrow 0
\]
where the maps  are
	\[\partial_0(f(\emptyset;u_i))=u_i,\ \mbox{for}\ 1\leq i\leq 5, 
\]
so \[\partial_0=
\left(\begin {array}{ccccc}
			x_2^3& x_1x_2^2& x_1x_2x_3& x_1x_3^2& x_1^2x_2
	\end{array}\right) .\]
\[\begin{array}{lll}
\partial_1(f(\{2\};u_2))&= &-x_2f(\emptyset;u_2)+x_1f(\emptyset;u_1),\\
\partial_1(f(\{2\};u_3))&= &-x_2f(\emptyset;u_3)+x_3f(\emptyset;u_2),\\
\partial_1(f(\{2\};u_4))&= &-x_2f(\emptyset;u_4)+x_3f(\emptyset;u_3),\\
\partial_1(f(\{2\};u_5))&= &-x_2f(\emptyset;u_5)+x_1f(\emptyset;u_2),\\
\partial_1(f(\{3\};u_5))&= &x_3f(\emptyset;u_5)-x_1f(\emptyset;u_3),\\
\end{array}\]
so
\[\partial_1=
\left(\begin {array}{ccccc}
			x_1& 0& 0& 0& 0\\
			-x_2& x_3& 0& x_1& 0\\
			0& -x_2& x_3& 0& -x_1\\
			0& 0& -x_2& 0& 0\\
			0& 0& 0& -x_2& x_3
	\end{array}\right) .\]

	\[\partial_2(f(\{2,3\};u_5))=-x_2f(\{3\};u_5)+x_3f(\{2\};u_5)+x_1f(\{3\};u_2)-x_1f(\{2\};u_3)=\]\[=-x_2f(\{3\};u_5)+x_3f(\{2\};u_5)-x_1f(\{2\};u_3),
\]
since $\{3\}\nsubseteq\set(u_2)$, so\[\partial_2=
\left(\begin {array}{c}
			0\\
			-x_1\\
			0\\
			x_3\\
			-x_2
	\end{array}\right) .\]
\end{Example1}
\section{Non-completely lexsegment ideals with a linear resolution}

\begin{Theorem1}\label{noncompletely}
Let $u=x_1^{a_1}\cdots x_n^{a_n},\ v=x_2^{b_2}\cdots x_n^{b_n}$ be monomials of degree $d$ in $k[x_1,\ldots,x_n],$with $a_1\neq 0.$ Suppose that the ideal $I=({\mathcal L}(u,v))$ is not a completely lexsegment ideal. Then $I$ has a linear resolution if and only if $I$ has linear quotients.
\end{Theorem1}

\begin{proof}
We only have to proof that if $I$ has a linear resolution then $I$ has linear quotients. By Theorem \ref{noncompletelylex}, since $I$ has a linear resolution, $u$ and $v$ have the form: $$u=x_1x_{l+1}^{a_{l+1}}\cdots x_n^{a_n},\ v=x_lx_n^{d-1}, \text {for some } l\geq 2.$$ Then the ideal $I=({\mathcal L}(u,v))$ can be written as a sum of ideals $I=J+K,$ where $J$ is the ideal generated by all the monomials of ${\mathcal L}(u,v)$ which are not divisible by $x_1$ and $K$ is generated by all the monomials of ${\mathcal L}(u,v)$ which are  divisible by $x_1.$  More precise, we have $$J=(\{w\ :\ x_2^d\geq_{lex}w\geq_{lex}v\})$$ and $$K=(\{w\ :\ u\geq_{lex}w\geq_{lex}x_1x_n^{d-1}).$$ One may see that $J$ is generated by
 the initial lexsegment $\mathcal{L}^i(v)\subset k[x_2,\ldots,x_n],$ and hence it has linear quotients with respect to
 lexicographical order $>_{lex}.$ Let $G(J)=\{g_1\prec \ldots\prec g_m\},$ where $g_i\prec g_j\ \text{if and only if}\
 g_i>_{lex} g_j.$ The ideal $K$ is isomorphic with the ideal generated by the final lexsegment of degree $d-1$
 $$\mathcal{L}^f(u/x_1)=\{w\mid u/x_1\geq_{lex} w\geq_{lex} x_n^{d-1},\ \deg(w)=d-1\}.$$ Since final lexsegments
  are stable with respect to the order  $x_n>\ldots >x_1$ of the variables, it follows that the ideal $K$ has linear quotients
 with respect to $>_{\overline{lex}},$ where by $\overline{lex}$ we mean the lexicographical order corresponding to 
 $x_n>\ldots >x_1.$ Let $G(K)=\{h_1\prec \ldots\prec h_p\},$ where $h_i\prec h_j\ \text{if and only if}\
 h_i>_{\overline{lex}} h_j.$ We consider the following ordering of the monomials of $G(I)$ $$G(I)=\{g_1\prec \ldots\prec
 g_m\prec h_1\prec \ldots\prec h_p\}.$$ We claim that, for this ordering of its minimal monomial generators, $I$ has linear quotients. In order to check this, we firstly notice that $I_{\prec g}:g=J_{\prec g}:g$ for every  $g\in G(J).$ Since $J$
 has linear quotients with respect to $\prec$, it follows that $J_{\prec g}:g$ is generated by variables. Now it is enough
to show that, for any generator $h$ of $K,$ the colon ideal $I_{\prec h}:h$ is generated by variables. 
We note that $$I_{\prec h}:h=J:h + K_{\prec h}:h.$$ Since $K$ is with linear quotients, we already know that
$K_{\prec h}:h$ is generated by variables. Therefore we only need to prove that $J:h$ is generated by variables. We
will show that $J:h=(x_2,\ldots,x_l)$ and this will end our proof. Let $m\in J:h$ be a monomial. It follows that $mh\in J.$ Since $h$ is a generator of $ K,$ $h$ is of the form $h=x_1x_{l+1}^{\alpha_{l+1}}\cdots x_n^{\alpha_n},$ that is $h\not\in (x_2,\ldots,x_l).$ But this implies that $m$ must be in the ideal $(x_2,\ldots,x_l).$ For the reverse inclusion, let $2\leq t\leq  l.$ Then $x_t h=x_1 \gamma$ for some monomial $\gamma,$ of degree $d.$ Replacing $h$ in the equality we get  $\gamma=x_tx_{l+1}^{\alpha_{l+1}}\cdots x_n^{\alpha_n}$ which shows that $\gamma$ is a generator of $J.$ Hence $x_th\in J.$
\end{proof}

\begin{Example1}\rm 
Let $I=(\mathcal{L}(u,v))\subset k[x_1,\ldots,x_6]$ be the lexsegment ideal of degree $4$ determined by the monomials $u=x_1x_3^2x_5$ and $v=x_2x_6^3.$  $I$ is not a completely lexsegment ideal as it follows applying Theorem \ref{condcompletely}, but $I$ has a linear resolution by Theorem \ref{noncompletelylex}. $I$ has linear quotients if we order its minimal monomial generators as indicated in the proof of the above theorem. On the other hand, if we order the generators of $I$ using the order relation defined in the proof of Theorem \ref{colex}, then we can easy see that $I$ does not have linear quotients. Indeed, following the definition of the order relation from Theorem \ref{colex}, we should take $$G(I)=\{x_2^4\prec x_2^3x_3\prec\ldots\prec x_2x_6^3\prec x_1x_3^2x_5\prec x_1x_3^2x_6\prec x_1x_3x_4^2\prec\ldots\prec x_1x_6^3\}.$$ For $h=x_1x_3x_4^2$ one may easy check that  $I_{\prec h}:h$ is not generated by variables.
\end{Example1}

\begin{Example1}\rm Let $u=x_1x_3x_4$, $v=x_2x_4^2$ be monomials in $k[x_1,\ldots, x_4]$. The ideal $I=(\mathcal{L}(u,v))\subset k[x_1,\ldots,x_4]$ is a non-completely lexsegment, since it does not verify the Theorem \ref{condcompletely}(b). By Theorem \ref{noncompletelylex}, $I$ has a linear resolution and by the proof of Theorem \ref{noncompletely}, $I$ has linear quotients with respect to the following ordering  of its minimal  monomial generators:
$$x_2^3,\ x_2^2x_3,\ x_2^2x_4,\ x_2x_3^2,\ x_2x_3x_4,\ x_2x_4^2,\ x_1x_4^2,\ x_1x_3x_4.$$
We note that $\set(x_1x_4^2)=\{2\}$ and $\set(g(x_1x_2x_4^2))=\set(x_2x_4^2)=\{2,3\}\nsubseteq\set(x_1x_4^2)$, so the decomposition function is not regular for this ordering of the generators.
\end{Example1}
\section{Cohen-Macaulay lexsegment ideals} \label{Section3}

In this section we study the dimension and the depth of arbitrary lexsegment ideals. These results are applied to describe the lexsegments ideals which are Cohen-Macaulay. We begin with the study of the dimension. As in the previous sections, let $d\geq 2$ be an integer. We denote $\frak{m}=(x_1,\ldots,x_n).$ It is clear that if $I=(\mathcal{L}(u,v))\subset S$ is a lexsegment ideal of degree $d$, then $\dim(S/I)=0$ if and only if $I=\frak{m}^d.$

\begin{Proposition1}\label{dim} Let $u=x_1^{a_1}\cdots x_n^{a_n},\ v=x_q^{b_q}\cdots x_n^{b_n}$, $1\leq q\leq n$, $a_1,b_q>0,$ be two monomials of degree $d$ such that $u\geq_{lex}v$ and let $I$ be the lexsegment ideal generated by $\mathcal{L}(u,v)$. We assume that $I\neq\frak{m}^d$. Then
\[\dim(S/I)=
\left\{\begin {array}{ll}
			n-q, & \mbox{if}\ 1\leq q<n,\\
			1, & \mbox{if}\ q=n.
	\end{array}\right. \]
\end{Proposition1}

\begin{proof} For $q=1$, we have $I\subset(x_1)$. Obviously $(x_1)$ is a minimal prime of $I$ and $\dim(S/I)=n-1$.

Let $q=n$, that is $v=x_n^d$ and $\mathcal{L}(u,v)=\mathcal{L}^f(u)$. We may write the ideal $I$ as a sum of two ideals, 
$I=J+K,$ where
$J=(x_1\mathcal{L} (u/x_1,x_n^{d-1}))
$
and
$K=(\mathcal{L}(x_2^d,x_n^d)).
$
Let $p\supset I$ be a monomial prime ideal. If $x_1\in p$, then $J\subseteq p$. Since $p$ also contains $K$, we have  $p\supset (x_2,\ldots,x_n)$. Hence $p=(x_1,x_2,\ldots,x_n)$. If $x_1\notin p$, we obtain  $(x_2,\ldots,x_n)\subset p$. Hence, the only minimal prime ideal of $I$ is $(x_2,\ldots,x_n)$. Therefore, $\dim(S/I)=1$.

Now we consider $1<q<n$ and write $I$ as before,
$I=J+K,
$
where $J=(x_1\mathcal{L}(u/x_1,x_n^{d-1}))$ and $K=(\mathcal{L}(x_2^d,v))$. 

Firstly we consider $u=x_1^d.$ Let $p\supset I$ be a monomial prime ideal. Then $p\ni x_1$ and, since $p\supset K,$ we also have $p\supset (x_2,\ldots,x_q).$ Hence $(x_1,\ldots,x_q)\subset p.$ Since $I\subset (x_1,\ldots,x_q),$ it follows that $(x_1,\ldots,x_q)$ is the only minimal prime ideal of $I.$ Therefore $\dim(S/I)=n-q.$

Secondly, let $a_1>1$ and $u\neq x_1^d.$ The lexsegment $\mathcal{L}(u/x_1,x_n^{d-1})$ contains the lexsegment $\mathcal{L}(x_2^{d-1},x_n^{d-1})$. Let $p$ be a monomial prime ideal which contains $I$ and such that $x_1\not\in p$. Then $p\supset\mathcal{L}(x_2^{d-1},x_n^{d-1})$ which implies that $(x_2,\ldots,x_n)\subset p$. Obviously we also have $I\subset(x_2,\ldots,x_n)$, hence $(x_2,\ldots,x_n)$ is a minimal prime ideal of $I$.

Let $p\supset I$ be a monomial prime ideal which contains $x_1$. Since $p\supset K$, we also have $(x_2,\ldots,x_q)\subset p$. This shows that $(x_1,\ldots, x_q)$ is a minimal prime ideal of $I$. In conclusion, for $a_1>1$, the minimal prime ideals of $I$ are $(x_1,\ldots,x_q)$ and $(x_2,\ldots, x_n)$. Since $q\leq n-1$, we get  $\height(I)=q$ and $\dim(S/I)=n-q$.

Finally, let $a_1=1$, that is $u=x_1 x_l^{a_l}\cdots x_n^{a_n}$, for some $a_l>0$, $l\geq2$. As in the previous case, we obtain that $(x_1,\ldots,x_q)$ is a minimal prime ideal of $I$. Now we look for those minimal prime ideals of $I$ which do not contain $x_1$.

If $a_l=d-1$, the ideal $J=(x_1\mathcal{L}(u/x_1,x_n^{d-1}))$ becomes $J=(x_1\mathcal{L}(x_l^{d-1},x_n^{d-1}))$. If $p\supset I$ is a monomial prime ideal such that  $x_1\notin p$, we get  $(x_l,\ldots,x_n)\subset p$, and, since $p$ contains $K$, we obtain  $(x_2,\ldots,x_q)\subset p$. This shows that, if $q<l$, then $(x_2,\ldots,x_q,x_l,\ldots,x_n)$ is a minimal prime ideal of $I$ of height $q+n-l\geq q$, and if $q\geq l$, then $(x_2,\ldots,x_n)$ is a minimal prime ideal of height $n-1\geq q$. In both cases we may draw the conclusion that $\height(I)=q$ and, consequently, $\dim(S/I)=n-q$.

The last case we have to consider is $a_l<d-1$. Then $l<n$ and, with similar arguments as above, we obtain  $\dim(S/I)=n-q$.
\end{proof}

In order to study the depth of arbitrary lexsegment ideals, we  note that one can restrict to those lexsegments defined by monomials of the form $u=x_1^{a_1}\cdots x_n^{a_n},\ v=x_1^{b_1}\cdots x_n^{b_n}$ of degree $d$ with $a_1>0$ and $b_1=0$.

Indeed, if $a_1=b_1$, then $I=(\mathcal{L}(u,v))$ is isomorphic, as an $S-$module, with the ideal generated by the lexsegment $\mathcal{L}(u/x_1^{a_1},v/x_1^{b_1})$ of degree $d-a_1$. This lexsegment may be studied in the polynomial ring in a smaller number of variables.

If $a_1>b_1$, then $I=(\mathcal{L}(u,v))$ is isomorphic, as an $S-$module, with the ideal generated by the lexsegment $\mathcal{L}(u',v')$, where $u'=u/x_1^{b_1}$ has $\nu_1(u')=a_1-b_1>0$ and $v'=v/x_1^{b_1}$ has $\nu_1(v')=0$.

Taking into account these remarks, from now on, we  consider lexsegment ideals of ends $u=x_1^{a_1}\cdots x_n^{a_n}$, $v=x_q^{b_q}\cdots x_n^{b_n}$, for some $q\geq2$, $a_1,b_q>0$.

The first step in the depth's study  is the next

\begin{Proposition1}\label{depthzero} Let $I=(\mathcal{L}(u,v))$, where $u=x_1^{a_1}\cdots x_n^{a_n}$, $v=x_q^{b_q}\cdots x_n^{b_n}$, $q\geq2$, $a_1,b_q>0$. Then $\depth(S/I)=0$ if and only if $x_nu/x_1\geq_{lex} v$.
\end{Proposition1}
\begin{proof} Let $x_nu/x_1\geq_{lex} v$. We claim that
$(I\colon (u/x_1))=(x_1,\ldots,x_n).
$
Indeed, for $1\leq j\leq n$, the  inequalities 
$$u\geq_{lex}\frac{x_ju}{x_1}\geq_{lex}\frac{x_nu}{x_1}\geq_{lex}v
$$ hold.
They show that $x_ju/x_1\in I$ for $1\leq j\leq n$. Therefore $(x_1,\ldots,x_n)\subseteq(I\colon (u/x_1))$. The other inclusion is obvious. We conclude that $(x_1,\ldots,x_n)\in \Ass(S/I)$, hence $\depth(S/I)=0$.

For the converse, let us assume, by contradiction, that $x_nu/x_1<_{lex}v$. We will show that $x_1-x_n$ is regular on $S/I$. This will imply that $\depth(S/I)>0$, which contradicts our hypothesis. We firstly notice that, from the above inequality, we have $a_1-1=0$, that is $a_1=1$. Therefore, $u$ is of the form $u=x_1x_l^{a_l}\cdots x_n^{a_n}$, $l\geq2,\ a_l>0$. Moreover, we have $l\geq q.$

Let us suppose that $x_1-x_n$ is not regular on $S/I$, that is there exists at least a polynomial $f\notin I$ such that $f(x_1-x_n)\in I$. One may assume that all monomials of $\supp(f)$ do not belong to $I$. Let us choose such a polynomial $f=c_1w_1+\ldots+c_tw_t,\ c_i\in k,\ 1\leq i\leq t$, with $w_1>_{lex}w_2>_{lex}\ldots>_{lex}w_t,\ w_i\notin I,\ 1\leq i\leq t$.

Then $\ini_{lex}((x_1-x_n)f)=x_1w_1\in I$. It follows that there exists $\alpha\in G(I)$ such that
\begin{eqnarray}x_1w_1=\alpha\cdot\alpha\,'.\label{3.1}
\end{eqnarray}
for some monomial $\alpha\,'$. We have $x_1\nmid\alpha\,'$ since, otherwise, $w_1\in I$, which is false. Hence $\alpha$ is a minimal generator of $I$ which is divisible by $x_1$, that is $\alpha$ is of the form $\alpha=x_1\gamma$, for some monomial $\gamma$ such that $x_n^{d-1}\leq_{lex}\gamma\leq_{lex}u/x_1$. Looking at $(\ref{3.1})$, we get  $w_1=\gamma \alpha\,'$. This equality shows that $x_1\nmid w_1$. We claim that the monomial $x_nw_1$ does not cancel in the expansion of $f(x_1-x_n)$. Indeed, it is clear that $x_nw_1$ cannot cancel by some monomial $x_nw_i$,  $i\geq2$. But it also cannot cancel by some monomial of the form $x_1w_i$ since $x_nw_1$ is not divisible by $x_1$. Now we may draw the conclusion that there exists a monomial $w\notin I$ such that $w(x_1-x_n)\in I$, that is $wx_1,\ wx_n\in I$.

Let $w\notin I$ be a monomial  such that $wx_1,\ wx_n\in I$, let $\alpha,\ \beta\in\mathcal{L}(u,v)$ and $\alpha\,',\ \beta\,'$ monomials such that
\begin{eqnarray}x_1w=\alpha\cdot\alpha\,'\label{3.2}
\end{eqnarray}
and
\begin{eqnarray}x_nw=\beta\cdot\beta\,'.\label{3.3}
\end{eqnarray}
As before, we get  $x_1\nmid w$, hence $\beta$ must be a minimal generator of $I$ such that
$x_2^d\geq_{lex}\beta\geq_{lex}v.
$ By using (\ref{3.3}), we can see that $x_n$ does not divide $\beta\,',$ hence $x_n|\beta.$  It follows that $w$ is divisible by $\beta/ x_n$. The monomial $w$ is also divisible by $\alpha/x_1$. Therefore,
 $\delta=\lcm(\alpha/x_1,\beta/x_n) | w.$ If $\deg \delta\geq d$ there exists a variable $x_j,$ with $j\geq 2,$ such that $(x_j\beta/x_n )| \delta,$ thus $(x_j\beta/x_n) | w.$ It is obvious that $x_2^d\geq_{lex} x_j\beta/x_n \geq_{lex} \beta\geq_{lex} v,$ hence $x_j\beta/x_n $ is a minimal generator of $I$ which divides $w,$ contradiction. This implies that $\delta$ has the degree $d-1.$ This yields $\alpha/x_1=\beta/x_n.$ Then $\beta=x_n\alpha /x_1\leq_{lex} x_n u/x_1<_{lex} v,$ contradiction. 
\end{proof}
 Using the Auslander--Buchsbaum Theorem, we obtain the following corollary.
 
\begin{Corollary1}
Let $I=(\mathcal{L}(u,v))$, where $u=x_1^{a_1}\cdots x_n^{a_n}$, $v=x_q^{b_q}\cdots x_n^{b_n}$, $q\geq2$, $a_1,b_q>0$. Then $\projdim(S/I)=n$ if and only if $x_nu/x_1\geq_{lex} v$.
\end{Corollary1}

We can compute the depth in the case of a final lexsegment ideal.
\begin{Corollary1}Let $I=\left(\mathcal{L}^f(u)\right)$ be the ideal generated by the final lexsegment defined by $u=x_1^{a_1}\cdots x_n^{a_n},\ a_1>0$. Then $\depth(S/I)=0$.
\end{Corollary1}

The following corollary is obvious.
\begin{Corollary1}
Let $I=\left(\mathcal{L}^i(v)\right)$ be the ideal generated by the initial lexsegment defined by the monomial $v$. Then $\depth(S/I)=0$ if and only if $v\leq_{lex}x_1^{d-1}x_n$.
\end{Corollary1}

Next we are going to characterize the lexsegment ideals $I$ such that $\depth\ S/I>0$, that is $x_nu/x_1<_{lex}v$, which implies that $u$ has the form $u=x_1x_l^{a_l}\cdots x_n^{a_n}$, for some $l\geq2,\ a_l>0$ and $l>q,$ or $l=q$ and $a_q\leq b_q$. We  denote $u\,'=u/x_1=x_l^{a_l}\cdots x_n^{a_n}$. Then we have $x_nu\,'<_{lex}v$. From the proof of Proposition \ref{depthzero} we know that $x_1-x_n$ is regular on $S/I$. Therefore 
	\[\depth(S/I)=\depth(S'/I')+1,
\]
where $S'=k[x_2,\ldots,x_n]$ and $I'$ is the ideal of $S'$ whose minimal monomial generating set is $G(I')=x_n\mathcal{L}(u\,',x_n^{d-1})\cup\mathcal{L}^i(v)$.

\begin{Lemma1}\label{depthneq0} In the above notations and  hypotheses on the lexsegment ideal $I$, the following statements hold:
\begin{itemize}
	\item[(a)] If $v=x_2^d$ and $l\geq4$, then $\depth(S'/I')=l-3$.
	\item[(b)] If $v=x_2^{d-1}x_j$ for some $3\leq j\leq n-2$ and $l\geq j+2$ then $\depth(S'/I')=l-j-1$.
	\item[(c)] $\depth(S'/I')=0$ in all the other cases. 
\end{itemize}
\end{Lemma1}

\begin{proof} (a) Let $v=x_2^d$ and $l\geq4$. The ideal $I'\subset S '$ is minimally generated by all the monomials $x_n\gamma$, where $x_n^{d-1}\leq_{lex}\gamma\leq_{lex}u\,'$, $\deg(\gamma)=d-1$, and by the monomial $x_2^d$. Then it is clear that $\{x_3,\ldots,x_{l-1}\}$ is a regular sequence on $ S'/I'$, hence
$$\depth S'/I'=\depth\frac{S'/I'}{(x_3,\ldots,x_{l-1})S'/I'}+l-3.
$$We have
$$\frac{S'/I'}{(x_3,\ldots,x_{l-1})S'/I'}\cong\frac{k[x_2,x_l,\ldots,x_n]}{I'\cap k[x_2,x_l,\ldots,x_n]}.
$$In this way we may reduce the computation of $\depth(S'/I')$ to the case (c).

(b) Let $v=x_2^{d-1}x_j$, for some $3\leq j\leq n-2$ and $l\geq j+2$. Hence $I'$ is minimally generated by the following set of monomials
	\[\{x_n\gamma\ |\ \gamma\ \mbox{is a monomial of degree}\ d-1\ \mbox{such that}\ x_n^{d-1}\leq_{lex}\gamma\leq_{lex}u\,' \}\cup\]\[\cup\{x_2^d,\ x_2^{d-1}x_3,\ldots,x_2^{d-1}x_j\}.
\]
Then $\{x_{j+1},\ldots,x_{l-1}\}$ is a regular sequence on $S'/I'$ and 
$$\depth S'/I'=\depth\frac{S'/I'}{(x_{j+1},\ldots,x_{l-1})S'/I'}+(l-j-1).
$$
Since
$$\frac{S'/I'}{(x_{j+1},\ldots,x_{l-1})S'/I'}\cong\frac{k[x_2,\ldots,x_j,x_l,\ldots,x_n]}{I'\cap k[x_2,\ldots,x_j,x_l,\ldots,x_n]},
$$
we may reduce the computation of $\depth(S'/I')$ to the case (c).

(c) In each of the cases that it remains to be treated, we will show that $$(x_2,\ldots,x_n)\in \Ass(S '/I '),$$ that is there exists a monomial $w\notin I'$ such that $I':w=(x_2,\ldots,x_n)$. This implies that $\depth(S'/I')=0$.

{\textit{Subcase $C_1$}.} $v=x_2^d,\ l=2$. Then $w=x_n^{d-1}\notin I'$ and 
$x_n^{d-1}\leq_{lex}x_jw/x_n=x_jx_n^{d-2}\leq_{lex}x_2x_n^{d-2}\leq_{lex}x_l^{a_l}\cdots x_n^{a_n}=u',
$
for all $2\leq j\leq n$. Hence $\gamma=x_jw/x_n$ has the property that $x_n\gamma\in G(I')$. Therefore, $x_j\in I':w$ for all $2\leq j\leq n$. It follows that $I':w=(x_2,\ldots,x_n)$.

{\textit{Subcase $C_2$}.} $v=x_2^d,\ l=3$. Then $w=x_2^{d-1}x_n^{d-1}\notin I'$. Indeed, $x_2^d\nmid w$ and if we assume that there exists $x_n^{d-1}\leq_{lex}\gamma\leq_{lex}u\,'$, $\deg\gamma=d-1,$ such that $x_n\gamma|w$,  we obtain  $x_n\gamma|x_n^{d-1}$ which is impossible.

We show that $x_jw\in I'$ for all $2\leq j\leq n$. Indeed, $x_2w=x_2^dx_n^{d-1}\in I'$. Let $3\leq j\leq n$. Then $x_n^{d-1}\leq_{lex}x_jx_n^{d-2}\leq_{lex}x_3x_n^{d-2}\leq_{lex} u\,'$. It follows that $\gamma=x_jx_n^{d-2}$ has the property that $x_n\gamma=x_jx_n^{d-1}\in G(I')$. Since $x_n\gamma|x_jw$, we have  $x_jw\in I'$. This arguments shows that $I':w=(x_2,\ldots,x_n)$.

{\textit{Subcase $C_3$}.} $v=x_2^{d-1}x_j$ for some $3\leq j\leq n-1$ and $2\leq l\leq j+1$. Let us consider again the monomial $w=x_2^{d-1}x_n^{d-1}$. It is clear that $x_tw\in I$ for all $2\leq t\leq j$. Let $t\geq j+1$. Then $x_tw$ is divisible by $x_tx_n^{d-1}$. Since $x_tx_n^{d-2}$ satisfies the  inequalities
$x_n^{d-1}\leq_{lex}x_tx_n^{d-2}\leq_{lex}u\,',
$
we have  $x_tx_n^{d-1}\in G(I')$. It follows that $x_t w\in I'$ for $t\geq j+1$. Assume that $w\in I'$. Since $x_2^{d-1}x_t\nmid w$ for $2\leq t\leq j$, we should have $x_n\gamma|w$ for some $\gamma$ of degree $d-1$ such that $x_n^{d-1}\leq_{lex}\gamma\leq_{lex}u\,'$. Since $\gamma|x_2^{d-1}x_n^{d-2}$ and $\gamma\leq_{lex}u\,'$, we get  $l=2$ and $a_2=\nu_2(u\,')\geq\nu_2(\gamma)$. Let $\gamma=x_2^ax_n^{d-1-a}$, for some $a\geq1$. In this case we change the monomial $w.$ Namely, we consider the monomial $w'=x_2x_n^{d-2}$ which does not belong to $G(I')$ since it has degree $d-1$.

 If $a_2\geq2$, for any $j$ such that $2\leq j\leq n$, we have 
$x_n^{d-1}<_{lex}x_jw'/x_n=x_2x_jx_n^{d-3}<_{lex}x_l^{a_l}\cdots x_n^{a_n}=u\,'.
$
This shows that $x_jw'\in I'$ for $2\leq j\leq n$ and hence, $I':w=(x_2,\ldots,x_n)$.

If $a_2=1$, we take $w\,''=x_n^{d-1}\notin I'$. For all $j$ such that $2\leq j\leq n$, we have 
$x_n^{d-1}\leq_{lex}x_j w\,''/x_n=x_jx_n^{d-2}\leq_{lex}x_2x_n^{d-2}\leq_{lex} u\,'.
$
Therefore $x_jw\,''\in I'$ for $2\leq j\leq n$, hence 
$I':w\,''=(x_2,\ldots,x_n).
$
In conclusion we have proved that in every case one may find a monomial $w\notin I'$ such that $I':w=(x_2,\ldots, x_n)$.

{\textit{Subcase $C_4$}.} Finally, let $v\leq_{lex} x_2^{d-1}x_n$. In this case, the ideal $I':x_2^{d-1}$ obviously contains $(x_2,\ldots,x_n)$. Since the other inclusion is trivial, we get $I':x_2^{d-1}=(x_2,\ldots,x_n)$. It is clear that $x_2^{d-1}\notin I'$.
\end{proof}

By using Lemma \ref{depthneq0} we get:

\begin{Proposition1}\label{depth} Let $I=(\mathcal{L}(u,v))$ be a lexsegment ideal defined by the monomials $u=x_1x_l^{a_l}\cdots x_n^{a_n},\ v=x_q^{b_q}\cdots x_n^{b_n}$ where $a_l,\ b_q>0$, $l,q\geq2$ and $x_n u/x_1<_{lex} v$. Then the following statements hold:
\begin{itemize}
	\item[(a)] If $v=x_2^d$ and $l\geq4$ then $\depth(S/I)=l-2$; 
	\item[(b)] If $v=x_2^{d-1}x_j$ for some $3\leq j\leq n-2$ and $l\geq j+2$ then $\depth(S/I)=l-j$;
	\item[(c)] $\depth(S/I)=1$ in all the other cases.
\end{itemize}
\end{Proposition1}

\begin{proof} Since $x_1-x_n$ is regular on $S/I$ if $x_nu/x_1<_{lex}v$, we have $\depth(S/I)=\depth(S'/I')+1$. The conclusion follows  applying Lemma \ref{depthneq0}.
\end{proof}

By using the Auslander--Buchsbaum Theorem, we get the following corollary:
\begin{Corollary1}
Let $I=(\mathcal{L}(u,v))$ be a lexsegment ideal defined by the monomials $u=x_1x_l^{a_l}\cdots x_n^{a_n},\ v=x_q^{b_q}\cdots x_n^{b_n}$ where $a_l,\ b_q>0$, $l,q\geq2$ and $x_nu/x_1<_{lex} v$. Then the following statements hold:
\begin{itemize}
	\item[(a)] If $v=x_2^d$ and $l\geq4$ then $\projdim(S/I)=n-l+2$; 
	\item[(b)] If $v=x_2^{d-1}x_j$ for some $3\leq j\leq n-2$ and $l\geq j+2$ then $\projdim(S/I)=n-l+j$;
	\item[(c)] $\projdim(S/I)=n-1$ in all the other cases.
\end{itemize}
\end{Corollary1}

As a consequence of the results of this section we may characterize the Cohen--Macaulay lexsegment ideals.

In the first place, we note that the only Cohen--Macaulay lexsegment ideal such that $\dim(S/I)=0$ is $ I=\frak m^d$. Therefore it remains to consider Cohen--Macaulay ideals $I$ with $\dim(S/I)\geq1$.

\begin{Theorem1} Let $n\geq3$ be an integer, let $u=x_1^{a_1}\cdots x_n^{a_n}$, $v=x_1^{b_1}\cdots x_n^{b_n}$, with $a_1>b_1\geq0,$   monomials of degree $d,$  and $I=(\mathcal{L}(u,v))\subset S$  the lexsegment ideal defined by $u$ and $v.$ We assume  that $\dim(S/I)\geq 1$. Then $I$ is Cohen--Macaulay if and only if one of the following conditions is fulfilled:
\begin{itemize}
	\item[(a)] $u=x_1x_n^{d-1}$ and $v=x_2^d$;
	\item[(b)] $v=x_{n-1}^ax_n^{d-a}$ for some $a>0$ and $x_n\ u/x_1<_{lex}v$. 
\end{itemize}
\end{Theorem1}
\begin{proof} Let $u,v$ be as in (a). Then $\dim(S/I)=n-2$, by Proposition \ref{dim} and $\depth(S/I)=n-2$ by using (a) in Proposition \ref{depth} for $n\geq4$ and (c) for $n=3$.

Let $u,\ v$ as in (b). Then $\dim(S/I)=1$ by Proposition \ref{dim}. By using Proposition \ref{depth}(c), we obtain  $\depth(S/I)=1$, hence $S/I$ is Cohen--Macaulay.

For the converse, in the first place, let us take $I$ to be Cohen--Macaulay of $\dim(S/I)=1$. By Proposition \ref{dim} we have $q=n$ or $q=n-1$. If $q=n$, then $v=x_n^d$ and $x_n u/x_1\geq_{lex}v$. By Proposition \ref{depthzero}, $\depth(S/I)=0$, so $I$ is not Cohen--Macaulay.

Let $q=n-1$, that is $v=x_{n-1}^ax_n^{d-a}$ for some $a>0$. By Proposition \ref{depthzero}, since $\depth(S/I)>0$, we must have $x_nu/x_1<_{lex}v$, thus we get (b).

Finally, let $\dim(S/I)\geq2$, that is $q\leq n-2$. By using Proposition \ref{depth}, we obtain  $q=2$. Therefore $\dim(S/I)=\depth(S/I)=n-2$. Using again Proposition \ref{depth} (a),(b), it follows that $u=x_1x_n^{d-1}$ and $v=x_2^d$. 
\end{proof}
%
%
%
%
%
%
%
\chapter{Constructible ideals}

We define a new class of monomial ideals, namely constructible ideals. We prove that a square-free constructible ideal is closely related to the notion of constructible simplicial complexes. Next we will determine some properties for constructible ideals and we give a formula for computing the Betti numbers.

\section{Constructible ideals and constructible simplicial complexes}
Let $S=k[x_1,\ldots,x_n]$ be a polynomial ring in $n$ variables over a field, $k$. For a monomial ideal $I$ of $S$, we will denote by $G(I)$ the minimal monomial set of generators of $I$.
\begin{Definition1}\label{cons1}\rm
A monomial ideal $I$ of $S$ generated in degree $q$ is a \it constructible ideal \rm if it can be obtained by the following recursive procedure:
\begin{itemize}
	\item [(i)] If $u$ is a monomial in $S$ and $I=(u)$, then $I$ is a constructible ideal;
	\item[(ii)] If $I_1$, $I_2$ are constructible ideals generated in degree $q$ and $I_1\cap I_2$ is a constructible ideal  generated in degree $q+1$, then $I_1+I_2$ is a  constructible ideal.
\end{itemize}
\end{Definition1}

We note that the recursion procedure will stop. Indeed, let $G(I)=\{u_1,\ldots,u_r\}$ be the minimal monomial system of generators of the ideal $I$ and we consider $I=I_1+I_2$. Hence, in $I_1\cap I_2$, the generators can contain each variable to a power which is less than or equal to the maximal power to which that variable appears in all the generators of $I$. 

Let $a_i$ be the maximum of the exponents of the variable $x_i$ in the minimal monomial generators of $I$ and let $\mathbf{\underline{a}}=(a_1,\ldots, a_n)$. The recursion procedure will stop after at most $|\mathbf{\underline{a}}|:=a_1+\ldots +a_n$ steps.

The above remarks show that we could consider also the following definition of constructible ideals.

Let $\mathbf{\underline{a}}=(a_1,\ldots, a_n)\in\mathbb{Z}_{>0}^n$. We denote \[\mathcal{M}_{\mathbf{\underline{a}}}=\{x_1^{b_1}\cdots x_n^{b_n}\ :\ 0\leq b_i\leq a_i\ \mbox{for all} \ 1\leq i\leq n\}\]
and
\[\mathcal{J}_{\mathbf{\underline{a}}}=\{I\ :\ I\ \mbox{is a monomial ideal of}\ S\ \mbox{with}\ G(I)\subseteq\mathcal{M}_{\mathbf{\underline{a}}}\}.\] 
We set $|\mathbf{\underline{a}}|=a_1+\ldots+a_n$.

We note that, if $I,\ J\in\mathcal{J}_{\mathbf{\underline{a}}}$, then $I\cap J\in\mathcal{J}_{\mathbf{\underline{a}}}$.

\begin{Definition1}\label{cons2}\rm Let $I\in\mathcal{J}_{\mathbf{\underline{a}}}$ be a monomial ideal generated in degree $q$. $I$ is an \it{$\mathbf{\underline{a}}$--constructible ideal}\rm$\ $ if it can be obtained by the following recursive procedure:
\begin{itemize}
	\item [(i)] If $u\in\mathcal{M}_{\mathbf{\underline{a}}}$ and $I=(u)$, then $I$ is an $\mathbf{\underline{a}}$--constructible ideal;
	\item[(ii)] If $I_1,\ I_2\in\mathcal{J}_{\mathbf{\underline{a}}}$ are $\mathbf{\underline{a}}$--constructible ideals generated in degree $q<|\mathbf{\underline{a}}|$ and $I_1\cap I_2\in\mathcal{J}_{\mathbf{\underline{a}}}$ is an $\mathbf{\underline{a}}$--constructible ideal  generated in degree $q+1$, then $I_1+I_2$ is an $\mathbf{\underline{a}}$--constructible ideal.
\end{itemize}
\end{Definition1}

Let $\mathbf{\underline{a}}=\mathbf{\underline{1}}$, where $\mathbf{\underline{1}}=(1,\ldots,1)\in\mathbb{Z}_{>0}^n$. One has
	\[\mathcal{M}_\mathbf{\underline{1}}=\{x_1^{b_1}\cdots x_n^{b_n}\ :\ b_i\in\{0,1\},\ 1\leq i\leq n\}
	\]
	and
	\[\mathcal{J}_\mathbf{\underline{1}}=\{I\ :\ I\ \mbox{is a monomial ideal of}\ S\ \mbox{with}\ G(I)\subseteq\mathcal{M}_\mathbf{\underline{1}}\}.
\]

Then, all the ideals from $\mathcal{J}_{\mathbf{\underline{1}}}$ are square-free monomial ideals. In particular, a $\mathbf{\underline{1}}-$constructible ideal is a square-free monomial ideal. We will say that a monomial ideal $I$ is a square-free constructible ideal if $I$ is a $\underline{\mathbf{1}}$--constructible ideal, that is $I$ is a square-free monomial ideal which is constructible.

It is also important to notice that the only ideal in $\mathcal{J}_{\mathbf{\underline{a}}}$ generated in
degree $|\mathbf{\underline{a}}|$ is the principal ideal $  I= (x_1^{a_1}\cdots x_n^{a_n}) $. This remark justifies that the recursion procedure of the above definition terminates.

It is obvious that a monomial ideal $I$ is a constructible ideal (in the sense of Definition \ref{cons1}) if and only if $I$ is an $\mathbf{\underline{a}}-$constructible ideal, for some $\mathbf{\underline{a}}\in\mathbb{Z}_{>0}^n$. Although Definition \ref{cons2} looks more technical, it will turn out that it is very useful in the proofs.

Let us consider the example of M. Hachimori \cite{Ha} of a constructible simplicial complex which is not shellable. The Stanley--Reisner ideal of its Alexander dual is a constructible ideal which does not have linear quotients (see Example \ref{constnonshell}). On the other hand, the Stanley--Reisner ideal of the Alexander dual of the Dunce Hat is an example of ideal with a linear resolution which is not constructible (see also Example \ref{CMnonconst}).

In the following lemma we will prove some properties of the simplicial complexes which are often used and which we will need in the proof of the main theorem of this section.

\begin{Lemma1}\label{lema} Let $\Delta$ be a simplicial complex on the vertex set $[n]$ and $\Delta_1,\Delta_2$ subcomplexes of $\Delta$. Then:
\begin{itemize}
	\item[(a)] $I_{\Delta_1^{\vee}}\cap I_{\Delta_2^{\vee}}=I_{(\Delta_1\cap\Delta_2)^{\vee}}$
	\item[(b)] $\Delta=\Delta_1\cup\Delta_2$ if and only if $I_{\Delta^{\vee}}=I_{\Delta_1^{\vee}}+I_{\Delta_2^{\vee}}$
\item[(c)] $\Delta=\Delta_1\cap\Delta_2$ if and only if $I_{\Delta}=I_{\Delta_1}+I_{\Delta_2}$
\item[(d)] $\Delta=\Delta_1\cap\Delta_2$ if and only if $\Delta^{\vee}=\Delta_1^{\vee}\cup\Delta_2^{\vee}$
\end{itemize}
\end{Lemma1}
\begin{proof}(a) The monomial $x_F$ is in $I_{\Delta_1^{\vee}}\cap I_{\Delta_2^{\vee}}$ if and only if $F$ is not a face in $\Delta_1^{\vee}$ and $F$ is not a face in $\Delta_2^{\vee}$, that is $F^c$ is a face in $\Delta_1\cap \Delta_2$. This is equivalent with the fact that $F$ is not a face in $(\Delta_1\cap \Delta_2)^{\vee}$. In other words, $x_F$ is in $I_{(\Delta_1\cap \Delta_2)^{\vee}}$. 

(b) "$\Rightarrow$" The monomial $x_F$ is in $I_{\Delta^{\vee}}$ if and only if $F$ is not a face in $\Delta^{\vee}$, that is  $F^c$ is a face in $\Delta=\Delta_1\cup \Delta_2$. This is equivalent with the fact that $F$ is a face in $\Delta_1$ or $F$ is a face in $\Delta_2$. This is possible if and only if $F$ is not in $\Delta_1^{\vee}$ or $F$ is not in $\Delta_2^{\vee}$. But this is equivalent with the fact that $x_F\in I_{\Delta_1^{\vee}}\cup I_{\Delta_2^{\vee}}$ so $x_F\in I_{\Delta_1^{\vee}}+I_{\Delta_2^{\vee}}$.

"$\Leftarrow$" Let $F$ be a face of $\Delta$, that is $F^c$ is not a face in $\Delta^{\vee}$ and $x_{F^{c}}$ is a monomial in $I_{\Delta^{\vee}}$. In other words, $x_{F^c}$ is in $I_{\Delta_1^{\vee}}$ or $x_{F^c}$ is in $I_{\Delta_2^{\vee}}$. But this is equivalent with the fact that $F$ is in $\Delta_1$ or $F$ is a face in $\Delta_{2}$, that is $F$ is a face of $\Delta_1\cup\Delta_2$.

(c) "$\Leftarrow$"  $F$ is a face in $\Delta$ if and only if $x_F\notin I_{\Delta}=I_{\Delta_1}+I_{\Delta_2}$ that is $x_F\notin I_{\Delta_1}$ and $x_F\notin I_{\Delta_2}$. But, this is equivalent with the fact that $F$ is a face in $\Delta_1\cap\Delta_2$. 

"$\Rightarrow$" Let $x_F$ be a monomial in $I_{\Delta}$, that is $F$ is not a face in $\Delta=\Delta_1\cap\Delta_2$. In other words, $F$ is not a face of $\Delta_1$ or $F$ is not a face of $\Delta_2$. But this is possible if and only if $x_F$ is a monomial in $I_{\Delta_1}$ or $x_F\in I_{\Delta_2}$. This is equivalent with the fact that $x_F$ is a monomial of $I_{\Delta_1}+I_{\Delta_2}$.

(d) By (c), we have that $\Delta=\Delta_1\cap\Delta_2$ if and only if $I_{\Delta}=I_{\Delta_1}+I_{\Delta_2}$, that is $I_{(\Delta^{\vee})^{\vee}}=I_{(\Delta_1^{\vee})^{\vee}}+I_{(\Delta_2^{\vee})^{\vee}}$. By (b), this is equivalent with the fact that $\Delta^{\vee}=\Delta_1^{\vee}\cup\Delta_2^{\vee}$.
\end{proof}
Now, we can relate the notion of constructible ideals to the concept of constructible simplicial complexes.
\begin{Theorem1}\label{constconst} Let $\Delta$ be a pure simplicial complex on the vertex set $[n]=\{1,\ldots,n\}$. The following statements are equivalent:
\begin{itemize}
	\item[(a)] $\Delta$ is constructible.
	\item[(b)]  $I_{\Delta^{\vee}}$ is a square-free constructible ideal.
\end{itemize}
\end{Theorem1}

\begin{proof} We show that $\Delta$ is a constructible simplicial complex if and only if $I_{\Delta^{\vee}}$ is a $\underline{\mathbf{1}}$--constructible ideal.

"(a) $\Rightarrow$ (b)" We use induction on the dimension of $\Delta$. 

If $\dim(\Delta)=0$, then the facets of $\Delta$ are the vertices and hence $I_{\Delta^{\vee}}$ is the monomial ideal generated by all square-free monomials of degree $n-1$ in the polynomial ring $k[x_1,\ldots,x_n]$, where $[n]$ is the vertex set of $\Delta$. One may easily check that this ideal has linear quotients with respect to the lexicographical order of the minimal monomial generators. We will prove in Proposition \ref{linquotconst} that  every monomial ideal with linear quotients generated in one degree is a constructible ideal. Hence $I_{\Delta^{\vee}}$ is a $\underline{\mathbf{1}}$--constructible ideal.

Assume that for all constructible simplicial complexes $\Delta$ of dimension at most $d-1$, $I_{\Delta^{\vee}}$ is a $\underline{\mathbf{1}}$--constructible ideal. Let $\Delta$ be a $d-$dimensional constructible simplicial complex. We prove by induction on the number of facets of $\Delta$ that $I_{\Delta^{\vee}}$ is a $\mathbf{\underline{1}}-$constructible ideal.

If $\Delta$ is a simplex, $\Delta=\langle F\rangle$, we have that $I_{\Delta^{\vee}}=(x_{F^c})$ and it is a $\mathbf{\underline{1}}-$constructible ideal, by definition.

Assume that for all constructible simplicial complexes of dimension at most $d-1$ and for all constructible simplicial complexes  of dimension at most $d$ with at most $r-1$ facets, $r\geq 2$, the Stanley--Reisner ideal of the Alexander dual is a $\underline{\mathbf{1}}$--constructible ideal. Let $\Delta$ be a $d$--dimensional constructible simplicial complex with the facets set $\mathcal{F}(\Delta)=\{F_1,\ldots, F_r\}$. Since $\Delta$ is constructible, there exist two $d$--dimensional constructible subcomplexes $\Delta_1$ and $\Delta_2$ such that $\Delta=\Delta_1\union\Delta_2$ and $\Delta_1\cap\Delta_2$ is a $(d-1)$--dimensional constructible simplicial complex. Since $\Delta_1$ and $\Delta_2$ are $d$--dimensional constructible simplicial complexes with at most $r-1$ facets, by the induction hypothesis, $I_{\Delta_1^{\vee}}$ and $ I_{\Delta_2^{\vee}}$ are $\mathbf{\underline{1}}$--constructible ideals generated in degree $n-d-1$. 

By Lemma \ref{lema}, we have that $I_{\Delta_1^{\vee}}\cap I_{\Delta_2^{\vee}}=I_{(\Delta_1\cap \Delta_2)^{\vee}}$. Since $\Delta_1\cap \Delta_2$ is a $(d-1)$--dimensional constructible simplicial complex, by the induction hypothesis we have that $I_{\Delta_1^{\vee}}\cap I_{\Delta_2^{\vee}}$ is a $\mathbf{\underline{1}}$--constructible ideal generated in degree $n-d$. Since $\Delta=\Delta_1\cup\Delta_2$, we have that $I_{\Delta^{\vee}}=I_{\Delta_1^{\vee}}+I_{\Delta_2^{\vee}}$.  Hence $I_{\Delta^{\vee}}=I_{\Delta_1^{\vee}}+I_{\Delta_2^{\vee}}$ is a $\mathbf{\underline{1}}$--constructible ideal.
	
"(b)$\Rightarrow $(a)" We use descending induction on the degree of the monomials from the minimal system of generators of the monomial ideal $I$.

If $I=(x_1\cdots x_n)$, then $I=I_{\Gamma}$ and $\Gamma^{\vee}=\langle\emptyset\rangle$, hence it is a constructible simplicial complex.

Assume that, for all $\underline{\mathbf{1}}$--constructible ideals $I$ generated in degree at least $q+1$, there exists a simplicial complex $\Gamma$ such that $I=I_{\Gamma}$ and $\Gamma^{\vee}$ is a constructible simplicial complex.

Let $I$ be a $\mathbf{\underline{1}}$--constructible ideal generated in degree $q<n$. We use induction on the number of monomials from the minimal system of generators of the ideal $I$.

If $u\in\mathcal{M}_{\mathbf{\underline{1}}}$ and $I=(u)$, then let $\Gamma$ be such that $I_{\Gamma}=(u)$ and let $F=\supp(u)$. Then $\Gamma^{\vee}$ is the simplex generated by $F^c$ and it is a constructible simplicial complex. 

Assume that, for all $\underline{\mathbf{1}}$--constructible ideals with at most $r-1$ monomials in the minimal system of generators and generated in degree $q$, there exists a simplicial complex $\Gamma$ such that $I_{\Gamma}=I$ and $\Gamma^{\vee}$ is a constructible simplicial complex.

Let $I$ be a $\mathbf{\underline{1}}$--constructible ideal with $|G(I)|=r$, $r\geq2$, generated in degree $q$ and let $\Gamma$ be the simplicial complex such that $I=I_{\Gamma}$. We have $I=I_1+I_2$, with $I_1,I_2$  $\mathbf{\underline{1}}$--constructible ideals generated in degree $q$ and $I_1\cap I_2$ is a $\mathbf{\underline{1}}$--constructible ideal generated in degree $q+1$. Let $\Gamma_1,\ \Gamma_2$ be simplicial complexes on the vertex set $[n]$ such that $I_{\Gamma_1}=I_1$ and $I_{\Gamma_2}=I_2$.
	
Since $I_{\Gamma}=I_{\Gamma_1}+I_{\Gamma_2}$ we have that $\Gamma=\Gamma_1\cap\Gamma_2$. $\Gamma=\Gamma_1\cap\Gamma_2$ implies that $\Gamma^{\vee}=\Gamma_1^{\vee}\union\Gamma_2 ^{\vee}$. We have to prove that $\Gamma^{\vee}$ is a constructible simplicial complex.

By induction hypothesis, $\Gamma_1^{\vee},\ \Gamma_2^{\vee}$ are constructible simplicial complex of dimension $n-q-1$.

We have that $I_1\cap I_2=I_{\Gamma_1}\cap I_{\Gamma_2}=I_{(\Gamma_1^{\vee}\cap\Gamma_2^{\vee})^{\vee}}$. So $\Gamma^{\vee}=\Gamma_1^{\vee}\union\Gamma_2 ^{\vee}$ is a $(n-q-1)-$dimensional constructible simplicial complex, which ends our proof.
\end{proof}

Now we can complete the diagram for pure simplicial complexes and the connections with the Stanley--Reisner ideal of the Alexander dual:
	\[
\]
\[\begin{array}{ccccc}
	\Delta\ \mbox{is shellable}&\Longrightarrow&\Delta\ \mbox{is constructible}&\Longrightarrow&\Delta\ \mbox{is Cohen--Macaulay}\\
	\Updownarrow&&\Updownarrow&&\Updownarrow\\
I_{\Delta^{\vee}}\	\mbox{is an ideal with}&\Longrightarrow&I_{\Delta^{\vee}}\ \mbox{is a constructible}&\Longrightarrow&I_{\Delta^{\vee}}\ \mbox{is an ideal with a}\\
	\mbox{linear quotients}&&\mbox{ideal}&&\mbox{linear resolution} 
\end{array}\cdot
\]

\section{Properties of constructible ideals}
From the above diagram, one may note that for square-free monomial ideals generated in one degree we have the following implications:
	\[
\]
\[\begin{array}{ccccc}
	
\mbox{ square-free monomial}&&\mbox{square-free }&&\mbox{ square-free monomial}\\
	\mbox{ideals with linear}&\Longrightarrow&\mbox{}&\Longrightarrow&\mbox{ideals with a linear }\\
	\mbox{quotients}&&\mbox{constructible ideals}&&\mbox{resolution} 
\end{array}.
\]
	\[
\]

Our aim is to prove that the same implications hold for monomial ideals generated in one degree, not necessarily square-free.
 
In this section we prove that every constructible ideal has a linear resolution and we compute the Betti numbers of a constructible ideal.

For this, we need the following two lemmas, the first one being given and proved in T. R\"omer \cite{Ro}. 

\begin{Lemma1}[\cite{Ro}]\label{roemer}
 Let $R$ be a standard graded $k$-algebra and
	\[0\rightarrow M'\rightarrow M\rightarrow M''\rightarrow0
\]
be an exact sequence of $\mathbb{Z}$-graded $R$-modules. If $M'$ and $M''$ have $q-$linear resolutions, then $M$ has a $q-$linear resolution.
\end{Lemma1}

\begin{Lemma1}\label{q,q+1}
 Let $R$ be a standard graded $k$-algebra and
	\[0\rightarrow M'\rightarrow M\rightarrow M''\rightarrow0
\]
be an exact sequence of $\mathbb{Z}$-graded $R$-modules. If $M'$ has a $(q+1)-$linear resolution and $M$ has a $q-$linear resolution, then $M''$ has a $q-$linear resolution.
\end{Lemma1}
\begin{proof} The exact sequence 
	\[0\rightarrow M'\rightarrow M\rightarrow M''\rightarrow0
\]
yields the exact sequence
	\[\ldots \rightarrow \Tor_{i}^{R}(M',k)_{i+j}\rightarrow  \Tor_{i}^{R}(M,k)_{i+j}\rightarrow \Tor_{i}^{R}(M'',k)_{i+j}\rightarrow \]
	\[\rightarrow \Tor_{i-1}^{R}(M',k)_{i+j}\rightarrow \Tor_{i-1}^{R}(M,k)_{i+j}\rightarrow\ldots
\]
Since $M$ has a $q-$linear resolution, $\Tor_{i}^{R}(M,k)_{i+j}=0$ for all $j\neq q$. In the same way, $\Tor_{i}^{R}(M',k)_{i+j}=0$ for all $j\neq q+1$.

For $j=q$ we obtain
\[0\rightarrow  \Tor_{i}^{R}(M,k)_{i+q}\rightarrow \Tor_{i}^{R}(M'',k)_{i+q}\rightarrow \Tor_{i-1}^{R}(M',k)_{i+q}\rightarrow 0.
\]
If $j=q+1$ we have $\Tor_{i}^{R}(M,k)_{i+q+1}=0$ and $$\Tor_{i-1}^{R}(M',k)_{i+q+1}=\Tor_{i-1}^{R}(M',k)_{(i-1)+(q+2)}=0,$$ so $\Tor_{i}^{R}(M'',k)_{i+q+1}=0$. Since $\Tor_{i}^{R}(M'',k)_{i+j}=0$ for all $j\neq q$, we have that $M''$ has a $q$--linear resolution.
\end{proof}

Now we may prove the main result of this section.

\begin{Theorem1}\label{constlinres}
\it Let $S=k[x_1,\ldots, x_n]$ be the polynomial ring over a field $k$ and $I$ be a constructible ideal of $S$ generated in degree $q$. Then $I$ has a $q$--linear resolution. 
\end{Theorem1}\begin{proof}Let $\mathbf{\underline{a}}\in\mathbb{Z}_{>0}^n$ such that $I$ is $\mathbf{\underline{a}}-$constructible. We use descending induction on the degree of monomials from the minimal system of generators of the monomial ideal $I$.

If $I=(x_1^{a_1}\cdots x_n^{a_n})$, $I$ has an $|\mathbf{\underline{a}}|$--linear resolution.

Assume that the statement holds for all $\underline{\mathbf{a}}$--constructible ideals generated in degree at least $q+1$, where $q<|\mathbf{\underline{a}}|$.

Let $I\in\mathcal{J}_{\mathbf{\underline{a}}}$ be an $\mathbf{\underline{a}}$--constructible ideal generated in degree $q$. Now we use induction on the number of monomials from the minimal system of generators.

If $u\in\mathcal{M}_{\mathbf{\underline{a}}}$, $\deg(u)=q$ and $I=(u)$, then $I$ has a $q-$linear resolution. 

Let $I$ be an $\mathbf{\underline{a}}$--constructible ideal generated in degree $q$ with $G(I)=\{u_1,\ldots,u_r\}$, $r\geq2$. There exist $\mathbf{\underline{a}}$--constructible ideals $I_1$ and $I_2$, generated in degree $q$, such that $I=I_1+I_2$ and $I_1\cap I_2$ is an $\mathbf{\underline{a}}$--constructible ideal generated in degree $q+1$. By the induction hypothesis, $I_1$ and $I_2$ have $q$--linear resolutions and $I_1\cap I_2$ has a $(q+1)$--linear resolution. From the exact sequence
 \[0\rightarrow I_1\rightarrow I_1\dirsum I_2\rightarrow I_2\rightarrow0
\]
we have, by Lemma \ref{roemer}, that $I_1\dirsum I_2$ has a $q$--linear resolution and from the exact sequence:
\[0\rightarrow I_1\cap I_2\rightarrow I_1\dirsum I_2\rightarrow I_1+I_2\rightarrow0,
\]
$I_1+I_2$ has a $q$--linear resolution, by Lemma \ref{q,q+1}. So $I$ has a $q$--linear resolution.
\end{proof}
One may note that Theorem \ref{constCM} is a special case of Theorem \ref{constlinres}, in view of Theorem \ref{constconst}.

 Analysing the proof of the above theorem, we get a method for computing the Betti numbers of a constructible ideal, as we describe in the following corollary.
 \begin{Corollary1} \label{betticonst}
\it Let $S=k[x_1,\ldots, x_n]$ be the polynomial ring over a field $k$ and $I$ a constructible ideal generated in degree $q$. Assume that $I_1$ and $I_2$ are constructible ideals generated in degree $q$ such that $I_1\cap I_2$ is a constructible ideal generated in degree $q+1$ and $I=I_1+I_2$. Then
	\[\beta_i(I)=\beta_i(I_1)+\beta_i(I_2)+\beta_{i-1}(I_1\cap I_2).
\]
\end{Corollary1}
\begin{proof}
The exact sequence
\[0\rightarrow I_1\cap I_2\rightarrow I_1\dirsum I_2\rightarrow I_1+I_2\rightarrow0
\]
yields the exact sequence
	\[0\rightarrow \Tor_{i}^S(I_1\dirsum I_2,k)_{i+q}\rightarrow \Tor_{i}^S(I_1+ I_2,k)_{i+q}\rightarrow  \Tor_{i-1}^S(I_1\cap I_2,k)_{i+q}\rightarrow  0.\]
From this sequence we get
	\[\beta_i(I)=\beta_i(I_1\dirsum I_2)+\beta_{i-1}(I_1\cap I_2).
\]
and, next,
\[\beta_i(I)=\beta_i(I_1)+\beta_i(I_2)+\beta_{i-1}(I_1\cap I_2).
\]
\end{proof}
\section{Polarization of constructible ideals}
We prove that the property of constructibility is preserved in the polarization process.

The polarization of a monomial ideal was introduced by R. Fr\"oberg in \cite{Fr}. In the polarization process, homological properties of a monomial ideal are preserved. Since the polarization of a monomial ideal is a square-free monomial ideal, we can apply specific techniques suited for these classes of ideals. 

First, we recall the notion of polarization of a monomial ideal and some concepts related to it, following A. Soleyman Jahan \cite{Ali}.

Let $S=k[x_1,\ldots,x_n]$ be the polynomial ring in $n$ variables over a field $k$ and $u=x_1^{\alpha_1}\cdots x_n^{\alpha_n}$ be a monomial of $S$. The \textit{polarization} of $u$ is the monomial
	\[u'=\prod_{i=1}^{n}\prod_{j=1}^{\alpha_i}x_{ij},
\]
where $u'\in k[x_{11},\ldots,x_{1\alpha_1},\ldots, x_{n1},\ldots,x_{n\alpha_n}]$. 

Let $I$ be a monomial ideal of $S$ and $u_1,\ldots, u_m$ be a system of monomial generators for $I$. Then, the ideal generated by monomials $u_1',\ldots, u_m'$ is called \textit{a polarization} of $I$. Since the polarization seems to depend on the system of generators, we consider another system of monomial generators, $v_1,\ldots, v_k$, for the monomial ideal $I$. Let $S'$ be a polynomial ring with sufficiently many variables such that for all $1\leq i\leq m$ and all $1\leq j\leq k$, $u_i'$ and $v_j'$ are monomials in $S'$. Then $(u_1',\ldots,u_m')=(v_1',\ldots,v_k')$ in $S'$. So we note that, in a common polynomial ring extension, all polarizations of a monomial ideal are the same. It follows that we can denote any polarization of a monomial ideal $I$ by $P(I)$. If $I$ and $J$ are two monomial ideals of $S$, we write $P(I)=P(J)$ if a polarization of $I$ and a polarization of $J$ coincide in a common polynomial ring extension.

\begin{Proposition1} Let $I$ be a constructible ideal of $S$. Then $P(I)$ is a square-free constructible ideal.
\end{Proposition1}

\begin{proof} Let $\mathbf{\underline{a}}\in\mathbb{Z}_{>0}^n$ such that $I$ is $\mathbf{\underline{a}}$--constructible.  We use descending induction on the degree of monomials from the minimal system of generators of the monomial ideal $I$.

If $I=(x_1^{a_1}\cdots x_n^{a_n})$, then $P(I)$ is a principal square-free monomial ideal, hence it is a $\mathbf{\underline{1}}$--constructible ideal.

Let $I\in\mathcal{J}_{\mathbf{\underline{a}}}$ be an $\mathbf{\underline{a}}$--constructible ideal generated in degree $q<|\mathbf{\underline{a}}|$. We use induction on the number of monomials from the minimal system of generators of the ideal $I$. 

If $u\in\mathcal{M}_{\mathbf{\underline{a}}}$, $\deg(u)=q$ and $I=(u)$, the statement is obvious. 

Let $I\in\mathcal{J}_{\mathbf{\underline{a}}}$ be an $\mathbf{\underline{a}}$--constructible ideal generated in degree $q$ with $|G(I)|=r$, $r\geq 2$. There exist $I_1,\ I_2\in\mathcal{J}_{\mathbf{\underline{a}}}$ $\ \mathbf{\underline{a}}$--constructible ideals generated in degree $q$ such that $I=I_1+I_2$ and $I_1\cap I_2$ is an $\mathbf{\underline{a}}$--constructible ideal generated in degree $q+1$. By induction hypothesis $P(I_1),\ P(I_2)$ and $P(I_1)\cap P(I_2)=P(I_1\cap I_2)$ are $\mathbf{\underline{1}}-$constructible ideals. Hence, $P(I)$ is a $\mathbf{\underline{1}}$--constructible ideal.
\end{proof}

\section{Ideals with linear quotients}
In this section we describe the relation between monomial ideals with linear quotients and constructible ideals.
\begin{Proposition1}\label{linquotconst}
\it Let $I$ be a monomial ideal of $S$ with linear quotients generated in degree $q$. Then $I$ is a constructible ideal. 
\end{Proposition1}

\begin{proof} We prove by induction on the number of monomials in the minimal system of generators.

If $u$ is a monomial in $S$ and $I=(u)$, then $I$ is a constructible ideal, by definition.

Let $I$ be a monomial ideal, $G(I)=\{u_1,\ldots,u_r\}$, $r\geq2$, be its minimal system of generators and assume that $I$ has linear quotients with respect to the sequence $u_1,\ldots, u_r$. Denote $I_1=(u_1,\ldots, u_{r-1})$ and $I_2=(u_r)$. $I_1,\ I_2$ are constructible ideals, by induction hypothesis.
	\[I_1\cap I_2=\left(u_i u_r/[u_i,u_r]:\ 1\leq i\leq r-1\right)=(x_{l_1}u_r,\ldots,x_{l_t}u_r),
\]
for some $l_i\in\{1,\ldots,n\},\ 1\leq i\leq t,\ t\leq r-1 $ where $[u_i,u_r]:=\gcd(u_i,u_r)$. The last equality holds since the ideal $I$ has linear quotients. So $I_1\cap I_2$ is a monomial ideal with linear quotients with at most $r-1$ monomials in the minimal system of generators. By the induction hypothesis, we have that $I_1\cap I_2$ is a constructible ideal and then $I$ is a constructible ideal.
\end{proof}

We may obtain the same formula as in Corollary \ref{Bettilinquot} for the Betti numbers of a monomial ideal with linear quotients generated in one degree  in the next proposition.

\begin{Proposition1}
\it
Let $I$ be a monomial ideal of $S$ generated in degree $q$, $G(I)=\{u_1,\ldots, u_m\}$ be its minimal system of generators and assume that $I$ has linear quotients with respect to the sequence $u_1,\ldots, u_m$. Denote $I_{k}=(u_1,\ldots, u_k)$ and let $d_k$ be the number of generators of the monomial ideal $I_{k-1}:(u_k)$. Then
	\[\beta_i(I)=\sum_{k=2}^m\left(
\begin{array}{c}
	d_k\\
	i
\end{array}
\right),
\]
for all $i\geq1$.
\end{Proposition1}

\begin{proof}Since $I_k$ is a monomial ideal with linear quotients with respect to the sequence $u_1,\ldots, u_k$, by Proposition \ref{linquotconst}, $I_k$ is a constructible ideal.

By Corollary \ref{betticonst},
	\[\beta_i(I_k)=\beta_{i}(I_{k-1})+\beta_{i}((u_k))+\beta_{i-1}(I_{k-1}\cap(u_k))\ \ \ \ \ (*)
\]
for all $k\geq2$, $i\geq1$. 

The multiplication by $u_k$ defines an isomorphism between $I_{k-1}:(u_k)$ and $I_{k-1}\cap(u_k)$. Therefore
	\[|G(I_{k-1}\cap(u_k))|=|G(I_{k-1}:(u_k))|=d_k.
\]
Since $I_k$ has linear quotients, the ideal $I_{k-1}:(u_k)$ is generated by a regular sequence of length $d_k$, and then 
	\[\beta_{i-1}(I_{k-1}:(u_k))=\left(\begin{array}{c}
	d_k\\
	i
\end{array}
\right).
\]
Summing in $(*)$ for $k=2,3,\ldots,m$, we get
	\[\beta_{i}(I)=\sum_{k=2}^{m}\beta_{i-1}(I_{k-1}\cap(u_k))=\sum_{k=2}^m\left(\begin{array}{c}
	d_k\\
	i
\end{array}\right).
\]
\end{proof}

We need a property of prime ideals given by A. Soleyman Jahan \cite{Ali}.
\begin{Lemma1}[\cite{Ali}]\label{ali}Let $I = (u_1, \ldots ,u_m)$ be a monomial ideal of $S$, and $u$ be a monomial in $S$.
Then $I : u$ is a prime ideal if and only if $P(I)  : u'$ is a prime ideal.
\end{Lemma1}

\begin{Proposition1}\label{pol} Let $I=(u_1,\ldots,u_r)$ be a monomial ideal of $S$. Then $I$ has linear quotients with respect to the sequence $u_1,\ldots, u_r$ if and only if $P(I)$ has linear quotients with respect to the sequence $u_1',\ldots, u_r'$.
\end{Proposition1}

\begin{proof} Let $2\leq k\leq r$. By Lemma \ref{ali}, $(u_1,\ldots,u_{k-1}):u_k$ is a prime ideal if and only if $(u_1',\ldots,u_{k-1}'):u_k'$ is a prime ideal. Since any monomial prime ideal is generated by a sequence of variables, the statement follows.
\end{proof}
\section{Examples}
In the sequel, we analise some examples. First two examples arise from the Stanley--Reisner ideal of suitable simplicial complexes.

\begin{Example1}\label{constnonshell}\rm  The following example of constructible and non-shellable simplicial complex is due to M. Hachimori \cite{Ha}.
\begin{figure}[h]
\begin{center}
\unitlength 1mm 
\linethickness{0.4pt}
\ifx\plotpoint\undefined\newsavebox{\plotpoint}\fi 
\begin{picture}(52.25,50.25)(0,0)
\put(3.5,19.75){\line(1,0){46.25}}
\put(3.5,28.5){\line(1,0){46.25}}
\multiput(27,48.25)(-.074519231,-.033653846){312}{\line(-1,0){.074519231}}
\put(3.75,37.75){\line(0,-1){17.5}}
\multiput(3.75,20.25)(.073636364,-.033636364){275}{\line(1,0){.073636364}}
\multiput(24,11)(.098076923,.033653846){260}{\line(1,0){.098076923}}
\put(49.5,19.75){\line(0,1){18}}
\multiput(49.5,37.75)(-.069749216,.03369906){319}{\line(-1,0){.069749216}}
\multiput(27.25,48)(-.033687943,-.045508274){423}{\line(0,-1){.045508274}}
\multiput(13,28.75)(.033639144,-.054281346){327}{\line(0,-1){.054281346}}
\thicklines
\put(4,37.5){\line(1,0){45.5}}
\multiput(13,29)(.08508403,.03361345){238}{\line(1,0){.08508403}}
\multiput(4,37.25)(.035714286,-.033673469){245}{\line(1,0){.035714286}}
\multiput(12.75,28.25)(-.034693878,-.033673469){245}{\line(-1,0){.034693878}}
\multiput(18.75,19.75)(-.03991597,-.03361345){119}{\line(-1,0){.03991597}}
\multiput(18.75,20)(.075961538,.033653846){260}{\line(1,0){.075961538}}
\multiput(38.25,28.5)(.043071161,.033707865){267}{\line(1,0){.043071161}}
\multiput(38.5,28.25)(.04456522,-.03369565){230}{\line(1,0){.04456522}}
\multiput(48.75,20.5)(.032609,-.032609){23}{\line(0,-1){.032609}}
\multiput(31.75,19.5)(.04910714,-.03348214){112}{\line(1,0){.04910714}}
\put(26.75,50.25){$0$}
\put(15.25,45){$3$}
\put(.5,38.25){$2$}
\put(.5,28.75){$3$}
\put(.5,19){$0$}
\put(23.5,7.5){$2$}
\put(12.75,12.25){$3$}
\put(38.75,13.5){$1$}
\put(51.5,19){$0$}
\put(52,28.5){$1$}
\put(52.25,38){$2$}
\put(38.75,44.75){$1$}
\put(20,34.5){$9$}
\put(11.5,24.5){$8$}
\put(32.25,33){$4$}
\put(39.25,31){$5$}
\put(19.5,21.5){$7$}
\put(30,21.5){$6$}
\thinlines
\multiput(33.5,37.75)(.03355705,.03691275){149}{\line(0,1){.03691275}}
\multiput(19.25,37.75)(-.03348214,.046875){112}{\line(0,1){.046875}}
\thicklines
\multiput(27.75,48)(.033653846,-.0625){312}{\line(0,-1){.0625}}
\multiput(28,48.25)(.033653846,-.0625){312}{\line(0,-1){.0625}}
\multiput(38.25,28.5)(-.033687943,-.042553191){423}{\line(0,-1){.042553191}}
\multiput(38.5,28.75)(-.033687943,-.042553191){423}{\line(0,-1){.042553191}}
\multiput(38.75,28.5)(-.033687943,-.042553191){423}{\line(0,-1){.042553191}}
\end{picture}
\end{center}
\caption{}
\end{figure}

The simplicial complex is constructible because we can split it by the bold line and we obtain two shellable simplicial complexes $\Delta_1,\ \Delta_2$ of dimension $2$ whose intersection is a shellable $1-$dimensional simplicial complex.
	
The shelling order of the facets for the simplicial complex $\Delta_1$ is:
	\[\{0,3,9\},\{2,3,9\},\{2,8,9\},\{2,3,8\},\{0,3,8\},\{0,7,8\},\{0,3,7\},\{2,3,7\},\{2,6,7\},\]
	\[\{5,6,7\},\{5,7,8\},\{4,5,8\},\{4,8,9\},\{0,4,9\}.
\]

For the simplicial complex $\Delta_2$, the shelling order of the facets is
	\[\{0,1,4\},\{1,2,4\},\{2,4,5\},\{1,2,5\},\{0,1,5\},\{0,5,6\},\{0,1,6\},\{1,2,6\}.
\]

The simplicial complex $\Delta_1\cap\Delta_2$ has the shelling order of the facets $\{0,4\},\{4,5\}$ $\{5,6\},\{2,6\}$.

For the Alexander dual of $\Delta_1$, the Stanley--Reisner ideal is
	\[I_{\Delta_1^{\vee}}=(x_1 x_2 x_4 x_5x_6 x_7 x_8,\ x_0x_1x_4x_5x_6x_7x_8,\ x_0x_1x_3x_4x_5x_6x_7,\]\[\ \ \ \ \ \ \ \ \ \ x_0x_1x_4x_5x_6x_7x_9,\ x_1x_2x_4x_5x_6x_7x_9,\ x_1x_2x_3x_4x_5x_6x_9,\]\[\ \ \ \ \ \ \ \ \ \ x_1x_2x_4x_5x_6x_8x_9,\ x_0x_1x_4x_5x_6x_8x_9,\ x_0x_1x_3x_4x_5x_8x_9,\]\[\ \ \ \ \ \ \ \ \ \ x_0x_1x_2x_3x_4x_8x_9,\ x_0x_1x_2x_3x_4x_6x_9,\ x_0x_1x_2x_3x_6x_7x_9,\]\[x_0x_1x_2x_3x_5x_6x_7,\ x_1x_2x_3x_5x_6x_7x_8).\ \ \ \ \ \ \ \ \ \ \ \ 
\]

The Stanley--Reisner ideal for the Alexander dual of $\Delta_2$ is:
\[I_{\Delta_2^{\vee}}=(x_2 x_3 x_5 x_6x_7 x_8 x_9,\ x_0x_3x_5x_6x_7x_8x_9,\ x_0x_1x_3x_6x_7x_8x_9,\]
\[\ \ \ \ \ \ \ \ \ \ x_0x_3x_4x_6x_7x_8x_9,\ x_2x_3x_4x_6x_7x_8x_9,\ x_1x_2x_3x_4x_7x_8x_9,\]
\[x_2x_3x_4x_5x_7x_8x_9,\ x_0x_3x_4x_5x_7x_8x_9).\ \ \ \ \ \ \ \ \ \ \ \ 
\]

The ideals $I_{\Delta_1^{\vee}}$ and $I_{\Delta_2^{\vee}}$ have linear quotients by Theorem \ref{linquot}. Since
	\[I_{\Delta_1^{\vee}}\cap I_{\Delta_2^{\vee}}=I_{(\Delta_1\cap\Delta_2)^{\vee}}=(x_1x_2 x_3 x_5 x_6x_7 x_8 x_9,\ x_0x_1x_2 x_3 x_6 x_7x_8 x_9,\]
	\[\ \ \ \ \ \ \ \ \ \ \ \ \ \ \ \ \ \ \ \ \ \ \ \ \ \ \ \ \ \ \ \ \ \ \ \ \ x_0x_1x_2 x_3 x_4x_7 x_8 x_9,\ x_0x_1x_3 x_4 x_5x_7 x_8 x_9),
\]
the ideal $I_{\Delta_1^{\vee}}\cap I_{\Delta_2^{\vee}}$ has linear quotients by Theorem \ref{linquot}. Hence $I_{\Delta^{\vee}}=I_{\Delta_1^{\vee}}+ I_{\Delta_2^{\vee}}$ is a square-free constructible ideal and $I_{\Delta^{\vee}}$ does not have linear quotients since $\Delta$ is not shellable.
\end{Example1}

Next, we consider a monomial ideal which has linear resolution and it is not constructible.

\begin{Example1}\label{CMnonconst}\rm
We consider now the Dunce Hat. It is known that the Dunce Hat is Cohen--Macaulay but it is not constructible (see M. Hachimori \cite{Ha}).
	
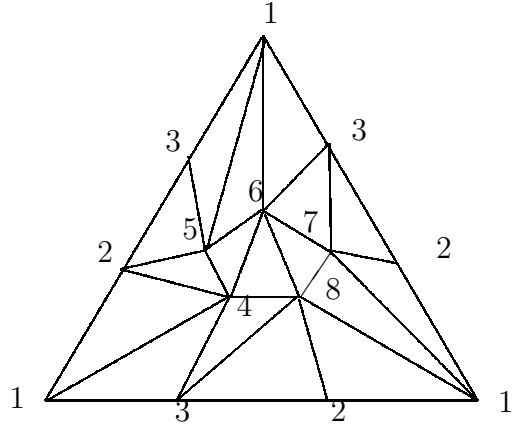
	\[
\]
\unitlength 1mm 
\linethickness{0.4pt}
\ifx\plotpoint\undefined\newsavebox{\plotpoint}\fi 
\begin{picture}(102,63.5)(0,0)
\multiput(70.5,58.5)(-.0337045721,-.0565650645){853}{\line(0,-1){.0565650645}}
\put(41.75,10.25){\line(1,0){57.5}}
\multiput(99.25,10.25)(-.0337278107,.0573964497){845}{\line(0,1){.0573964497}}
\put(70.75,58.75){\line(0,-1){23}}
\multiput(70.75,35.75)(-.0472561,-.03353659){164}{\line(-1,0){.0472561}}
\multiput(63,30.25)(.03350515,-.06443299){97}{\line(0,-1){.06443299}}
\put(66.25,24){\line(1,0){9.5}}
\put(75.75,24){\line(2,3){4}}
\multiput(79.75,30)(-.05487805,.03353659){164}{\line(-1,0){.05487805}}
\multiput(70.75,35.5)(-.03373016,-.08928571){126}{\line(0,-1){.08928571}}
\multiput(66.5,24.25)(-.059036145,-.03373494){415}{\line(-1,0){.059036145}}
\multiput(71,35)(.03358209,-.08022388){134}{\line(0,-1){.08022388}}
\multiput(75.5,24.25)(-.039156627,-.03373494){415}{\line(-1,0){.039156627}}
\multiput(59.25,10.25)(.03365385,.06610577){208}{\line(0,1){.06610577}}
\multiput(66.25,24)(-.12946429,.03348214){112}{\line(-1,0){.12946429}}
\multiput(51.75,27.75)(.15,.0333333){75}{\line(1,0){.15}}
\multiput(63,30.25)(-.0335821,.1865672){67}{\line(0,1){.1865672}}
\multiput(71,58.25)(-.03369565,-.12065217){230}{\line(0,-1){.12065217}}
\multiput(75.5,23.75)(.03348214,-.11830357){112}{\line(0,-1){.11830357}}
\multiput(75.75,24)(.057598039,-.03370098){408}{\line(1,0){.057598039}}
\multiput(99.25,10.25)(-.033703072,.034129693){586}{\line(0,1){.034129693}}
\multiput(79.5,30.25)(.1730769,-.0336538){52}{\line(1,0){.1730769}}
\multiput(79.75,29.75)(-.03125,1.84375){8}{\line(0,1){1.84375}}
\multiput(79.5,44.5)(-.033730159,-.034722222){252}{\line(0,-1){.034722222}}
\put(70.75,60.5){$1$}
\put(57.75,43.5){$3$}
\put(48.75,28.75){$2$}
\put(37,9.25){$1$}
\put(59,7.5){$3$}
\put(79.75,7.5){$2$}
\put(102,8.75){$1$}
\put(93.75,29.25){$2$}
\put(82.5,45){$3$}
\put(68.75,36.75){$6$}
\put(60,31.75){$5$}
\put(67.25,21.5){$4$}
\put(79,23.75){$8$}
\put(76,32.75){$7$}
\end{picture}

The Stanley--Reisner ideal for the Alexander dual of $\Delta$ is:
	\[
\]
	\[I_{\Delta^{\vee}}=(x_3x_5x_6x_7x_8,\ x_3x_4x_5x_6x_8,\ x_3x_4x_5x_6x_7,\ x_2x_5x_6x_7x_8,\]	
	\[\ \ \ \ \ \ \ \ \ \ \ x_2x_4x_6x_7x_8,\ x_2x_4x_5x_7x_8,\ x_2x_3x_4x_7x_8,\ x_2x_3x_4x_5x_6,\] 
	\[\ \ \ \ \ \ \ \ \ \ \ x_1x_4x_6x_7x_8,\ x_1x_4x_5x_6x_8,\ x_1x_4x_5x_6x_7,\ x_1x_3x_6x_7x_8,\]
\[\ \ \ \ \ \ \ \ \ \ \ x_1x_2x_5x_6x_7,\ x_1x_2x_4x_5x_8,\ x_1x_2x_3x_7x_8,\ x_1x_2x_3x_5x_7,\]\[x_1x_2x_3x_4x_5)\ \ \ \ \ \ \ \ \ \ \ \ \ \ \ \ \ \ \ \ \ \ \ \ \ \ \ \ \ \ \ \ \ \ \ \ \ \ \ \ 
\]
and it is not a constructible ideal but it has a linear resolution.
\end{Example1}

The last example is no longer square-free. It is an example of a constructible ideal which does not have linear quotients. We argue this using the fact that constructibility is "preserved" during the polarization process.

\begin{Example1}\label{nonsqfree}\rm Let $I\in k[x_1,\ldots,x_8]$ be the monomial ideal
\[I=(x_1x_2x_5x_6x_7x_8,\ x_2x_3x_5x_6x_7x_8,\ x_2^2x_3x_5x_6x_7,\ x_2^2x_3x_4x_6x_7,\ x_1x_2^2x_3x_6x_7,\]	
	\[\ \ \ \ \ \ x_2x_3x_4x_5x_7x_8,\ x_2^2x_3x_4x_7x_8,\ x_1x_2x_3x_4x_7x_8,\ x_1^2x_3x_4x_7x_8,\ x_1^2x_3x_4x_5x_8,\]
  \[\ \ \ \ \ \ \ x_1x_3x_4x_6x_7x_8,\ x_1x_4x_5x_6x_7x_8,\  x_1^2x_4x_5x_6x_8,\ x_1^2x_2x_4x_5x_8,\ x_1x_2^2x_5x_6x_8,\ \]
  \[\ \ \ \ \ \ \ \ \  x_1x_2^2x_3x_6x_8,\ x_1^2x_2^2x_3x_6,\ x_1^2x_2^2x_5x_6,\ x_1^2x_2x_5x_6x_7,\ x_1^2x_2x_4x_5x_7,\ x_1^2x_2^2x_4x_5)\]
Then $I=I_1+I_2$, where
\[I_1=(x_1x_2x_5x_6x_7x_8,\ x_2x_3x_5x_6x_7x_8,\ x_2^2x_3x_5x_6x_7,\ x_2^2x_3x_4x_6x_7,\ x_1x_2^2x_3x_6x_7,\]
\[ \ \ \ \ \ \ \ \ x_2x_3x_4x_5x_7x_8,\ x_2^2x_3x_4x_7x_8,\ x_1x_2x_3x_4x_7x_8,\ x_1^2x_3x_4x_7x_8,\ x_1^2x_3x_4x_5x_8,\]
  \[\ \ \ x_1x_3x_4x_6x_7x_8,\ x_1x_4x_5x_6x_7x_8,\ x_1^2x_4x_5x_6x_8,\ x_1^2x_2x_4x_5x_8)\ \ \ \ \ \ \ \ \ \ \ \ \] 
  and
  \[I_2=(x_1x_2^2x_5x_6x_8,\ x_1x_2^2x_3x_6x_8,\ x_1^2x_2^2x_3x_6,\ x_1^2x_2^2x_5x_6,\ x_1^2x_2x_5x_6x_7,\ x_1^2x_2x_4x_5x_7,\]
  \[ x_1^2x_2^2x_4x_5)\ \ \ \ \ \ \ \ \  \ \ \ \ \ \ \ \ \ \ \  \ \ \ \ \ \ \ \ \ \ \  \ \ \ \ \ \ \ \ \ \ \ \ \ \ \ \ \ \ \ \  \ \ \ \ \ \ \ \ \ \ \ \ \ \ \ \ \ \ \ \ \ \ \  \]
  with
	\[I_1\cap I_2=(x_1x_2^2x_5x_6x_7x_8,\ x_1^2x_2x_5x_6x_7x_8,\ x_1^2x_2x_4x_5x_7x_8,\ x_1^2x_2^2x_4x_5x_8,\ x_1x_2^2x_3x_6x_7x_8,\]
	\[ x_1^2x_2^2x_3x_6x_7)\ \ \ \ \ \ \ \ \  \ \ \ \ \ \ \ \ \ \ \  \ \ \ \ \ \ \ \ \ \ \  \ \ \ \ \ \ \ \ \ \ \ \ \ \ \ \ \ \ \ \  \ \ \ \ \ \ \ \ \ \ \ \ \ \ 
\]  
Since $I_1,\ I_2$ are monomial ideals with linear quotients generated in degree $6$, and $\ \ \ \ I_1\cap I_2$ is a monomial ideal generated in degree $7$ with linear quotients, $I$ is a constructible ideal.

Let $\Delta=\Delta_1\cup\Delta_2$ be the simplicial complex, presented by G.M. Ziegler, with $10$ vertices and $21$ facets (see G.M. Ziegler \cite{Zi}):
	\[\Delta_1=\langle \{1,2,3,4\},\ \{1,2,4,9\},\ \{1,4,8,9\},\ \{1,5,8,9\},\ \{1,4,5,8\},\ \{1,2,6,9\},\ \ \ \ \ \ \ \]
	\[\ \ \ \ \ \ \ \ \{1,5,6,9\},\ \{1,2,5,6\},\ \{2,5,6,10\},\ \{2,6,7,10\},\ \{1,2,5,10\},\ \{1,2,3,10\},
\]
	\[\{2,3,7,10\},\ \{2,3,6,7\}\rangle\ \ \ \ \ \ \ \ \  \ \ \ \ \ \ \ \ \ \ \  \ \ \ \ \ \ \ \ \ \ \  \ \ \ \ \ \ \ \ \ \ \ \ \ \ \ \ \ \ \ \  \ \ \ \ \ \ \ 
\]
	\[\Delta_2=\langle \{1,3,4,7\},\ \{1,4,5,7\},\ \{4,5,7,8\},\ \{3,4,7,8\}, \ \{2,3,4,8\},\ \{2,3,6,8\},\ \ \ \ \ \ \ 
\]
	\[ \{3,6,7,8\}\rangle.\ \ \ \ \ \ \ \ \  \ \ \ \ \ \ \ \ \ \ \  \ \ \ \ \ \ \ \ \ \ \  \ \ \ \ \ \ \ \ \ \ \ \ \ \ \ \ \ \ \ \  \ \ \ \ \ \ \ \ \ \ \ \ \ \ \ \ \ \ \ \ \ \ \ 
\]
This simplicial complex is constructible but non-shellable (see G.M. Ziegler \cite{Zi}).

The polarization of $I$ in the polynomial ring $k[x_{1},x_{1,1},x_2,x_{2,1},x_3,\ldots,x_8]$, with $x_{1,1}=x_9$ and $x_{2,1}=x_{10}$, is the Stanley--Reisner ideal of the Alexander dual associated to the above simplicial complex. Then the ideal $I$ does not have linear quotients, by Proposition \ref{pol}.  
\end{Example1}

%
%
%
%
%
%
\chapter{Subword complexes in Coxeter groups}\label{sub}
A. Knutson and E. Miller \cite{KnMi} proved that subword complexes in Coxeter groups are vertex-decomposable. Since any vertex-decomposable simplicial complex is she- llable (L.J. Billera and J.S. Provan \cite{BiPr}), subword complexes in Coxeter groups are shellable. 

In this chapter, we get the following facts which will be published in our paper \cite{O2}:
\begin{itemize}
	\item We prove directly that subword complexes in Coxeter groups are shellable, by using the Alexander duality.
	\item As a consequence, we get a shelling on the facets of the subword complex. 
	\item We study the Stanley--Reisner ideal of the Alexander dual for a special class of subword complexes. 
	\item For this class, we prove that the Stanley--Reisner ring is a complete intersection ring.
\end{itemize}

\section{Subword complexes in Coxeter groups and Alexander duality}

Let $(W,S)$ be a Coxeter system, $Q=(\sigma_1,\ldots,\sigma_n)$ be a word in $W$, with $\sigma_i\in S$ for all $1\leq i\leq n$, and $\pi$ be an element in $W$. Let $k[x_1,\ldots,x_n]$ be the polynomial ring in $n$ variables over a field $k$, where $n$ is the size of the word $Q$ and $\Delta(Q,\pi)$ be the subword complex. 
We aim to determine a shelling order on the facets of $\Delta(Q,\pi)$. For this purpose, we consider the Stanley--Reisner ideal of the Alexander dual associated to $\Delta(Q,\pi)$. 

Let $\Delta$ be a simplicial complex on the vertex set $[n]$. We recall that
	\[I_{\Delta^{\vee}}=(\mathbf{x}_{[n]\setminus F}\ \mid\ F\in\mathcal{F}(\Delta)),
\]
where we denote by $\mathbf{x}_F$ the monomial $\prod\limits_{i\in F}x_i$ and by $\mathcal{F}(\Delta)$ the set of all the facets of $\Delta$.

Since $\Delta=\Delta(Q,\pi)$ is shellable (Theorem \ref{vdshell}), by Theorem \ref{linquot}, the Stanley--Reisner ideal of the Alexander dual of $\Delta$, $I_{\Delta^{\vee}}$, has linear quotients. In order to obtain a shelling on the facets of $\Delta$, we have to define an order on the monomials from the minimal monomial set of generators of $I_{\Delta^{\vee}}$ such that $I_{\Delta^{\vee}}$ has linear quotients with respect to this order.

Henceforth, for $P=(\sigma_{i_1},\ldots, \sigma_{i_m})$ a subword of $Q$, $m\leq n$, we shall denote by $\mathbf{x}_{P}$ the monomial $x_{i_1}\cdots x_{i_{m}}$ in $k[x_1,\ldots,x_n]$. For a monomial ideal $I$, we will denote by $G(I)$ the minimal monomial generating set of $I$.

 In the special context of the subword complexes, the Stanley--Reisner ideal of the Alexander dual is
	\[I_{\Delta^{\vee}}=(\mathbf{x}_P\ \mid\ P\subseteq Q,\ P\ \mbox{represents}\ \pi).
\]

\begin{Theorem1}\label{sclinquot} Let $\Delta$ be the subword complex $\Delta(Q,\pi)$. Then the Stanley--Reisner ideal of the Alexander dual, $I_{\Delta^{\vee}}$, has linear quotients with respect to the lexicographical order of the minimal monomial generators.
\end{Theorem1}

\begin{proof} Let $G(I_{\Delta^{\vee}})=\{\mathbf{x}_{P_1},\ldots,\mathbf{x}_{P_{r}}\}$. We assume that $\mathbf{x}_{P_1}>_{lex}\ldots>_{lex}\mathbf{x}_{P_{r}}$. Note that $P_i\subseteq Q$ and $P_i$ represents $\pi$ for all $1\leq i\leq r$.

We have to prove that $I_{\Delta^{\vee}}$ has linear quotients with respect to the sequence of monomials $\mathbf{x}_{P_1},\ldots,\mathbf{x}_{P_r}$, that is, for all $i\geq2$ and for all $j<i$ there exists an integer $l\in[n]$ and an integer $k$, $1\leq k<i,$ such that $\mathbf{x}_{P_k}/\gcd(\mathbf{x}_{P_i},\mathbf{x}_{P_k})=x_l$ and $x_l$ divides $\mathbf{x}_{P_j}/\gcd(\mathbf{x}_{P_i},\mathbf{x}_{P_j})$. 

Let us fix $i\geq 2$ and $j<i$. Since $P_i,\ P_j$ represent $\pi$, they are subwords of $Q$ of size $\ell(\pi)$. We assume that $P_i=(\sigma_{i_1},\ldots,\sigma_{i_{\ell(\pi)}})$ and $P_j=(\sigma_{j_1},\ldots,\sigma_{j_{\ell(\pi)}})$. Let $l\in [n]$ be an integer such that $i_t=j_t$ for all $1\leq t\leq l-1$ and $j_l<i_l$. Such an integer exists since $j<i$, that is $\mathbf{x}_{P_j}>_{lex}\mathbf{x}_{P_i}$.

Let $T$ be the subword of $Q$ of size $\ell(\pi)+1$ obtained from $P_i$ by adding $\sigma_{j_l}$ between $\sigma_{i_{l-1}}$ and $\sigma_{i_l}$, that is $$T=(\sigma_{i_1},\ldots,\sigma_{i_{l-1}},\sigma_{j_l},\sigma_{i_l},\ldots,\sigma_{i_{\ell(\pi)}}).$$ Since $P_i$ is a subword of $T$ that represents $\pi$, we have $\delta(T)\succeq\pi$ by Lemma \ref{pi}(a).

Let us assume that $\delta(T)\succ\pi$. Hence, by Lemma \ref{m-1}(c), $T$ represents an element $\tau\in W$, $\tau\succ\pi$, such that $\ell(\tau)=\ell(\pi)+1$. Since $P_i,\ P_j$ represent $\pi$, we have that
	\[\pi=\sigma_{i_1} \cdots \sigma_{i_{l-1}}\sigma_{i_l}\cdots \sigma_{i_{\ell(\pi)}}=\sigma_{j_1} \cdots \sigma_{j_{l-1}}\sigma_{j_l}\cdots \sigma_{j_{\ell(\pi)}}
\]
are reduced expressions for $\pi$ and, by the choice of $l$, we obtain
	\[\sigma_{i_l}\cdots \sigma_{i_{\ell(\pi)}}=\sigma_{j_l}\cdots \sigma_{j_{\ell(\pi)}},
\]
that is
\begin{eqnarray}\sigma_{j_l}\sigma_{i_l}\cdots \sigma_{i_{\ell(\pi)}}=\sigma_{j_{l+1}}\cdots \sigma_{j_{\ell(\pi)}}.\label{**}
\end{eqnarray}
Now, since $T$ represents $\tau$,
	\[\tau=\sigma_{i_1}\cdots \sigma_{i_{l-1}} \sigma_{j_l}\sigma_{i_l}\cdots \sigma_{i_{\ell(\pi)}}
\]
is a reduced expression for $\tau$. On the other hand, using the equality (\ref{**}), we have that
	\[\tau=\sigma_{i_1}\cdots \sigma_{i_{l-1}}\sigma_{j_l}\sigma_{i_l}\cdots \sigma_{i_{\ell(\pi)}}=\sigma_{i_1} \cdots \sigma_{i_{l-1}}\sigma_{j_{l+1}}\cdots \sigma_{j_{\ell(\pi)}}.
\]
This is impossible since we obtained that $\tau$ can be written as a product of $\ell(\pi)-1$ simple reflections which contradicts the fact that $\ell(\tau)=\ell(\pi)+1$.

Hence, the Demazure product of the word $T$ is $\pi$. By Lemma \ref{m-1}(b), there exists a unique $\sigma_{i_t}\neq \sigma_{j_l}$ such that both $T\setminus \sigma_{j_l}=P_i$ and $T\setminus \sigma_{i_t}$ represent $\pi$. Let us denote $P'=T\setminus \sigma_{i_t}$. We will show that $\mathbf{x}_{P'}>_{lex}\mathbf{x}_{P_i}$ which will end the proof. We note that $\mathbf{x}_{P'}=x_{j_l}\mathbf{x}_P/x_{i_t}$. Also $\mathbf{x}_{P'}\neq \mathbf{x}_P$ since $j_l\neq i_t$.

Assume by contradiction that $\mathbf{x}_{P'}<_{lex}\mathbf{x}_{P_i}$, that is $i_t<j_l$. Since both $P_i$ and $P'$ represent $\pi$, we have that
	\[\pi=\sigma_{i_1}\cdots \sigma_{i_t}\cdots \sigma_{i_{l-1}}\sigma_{i_l}\cdots \sigma_{i_{\ell(\pi)}}=\sigma_{i_1}\cdots \sigma_{i_{t-1}}\sigma_{i_{t+1}}\cdots \sigma_{i_{l-1}}\sigma_{j_l}\sigma_{i_l}\cdots \sigma_{i_{\ell(\pi)}}
\]
are two reduced expressions for $\pi$. The above equality can be written as
	\begin{eqnarray}\sigma_{i_t}\sigma_{i_{t+1}}\cdots \sigma_{i_{l-1}}=\sigma_{i_{t+1}}\cdots \sigma_{i_{l-1}}\sigma_{j_l}.\label{x}
\end{eqnarray}
On the other hand, $P_j$ also represents $\pi$, that is $\pi=\sigma_{j_1}\cdots \sigma_{j_l}\cdots \sigma_{j_{\ell(\pi)}}$. Now, since for all $1\leq k<l$, $i_k=j_k$, using (\ref{x}) we obtain that
	\[\pi=\sigma_{j_1}\cdots \sigma_{j_{l-1}}\sigma_{j_l}\cdots \sigma_{j_{\ell(\pi)}}=\sigma_{i_1}\cdots \sigma_{i_{l-1}}\sigma_{j_l}\cdots \sigma_{j_{\ell(\pi)}}=\]\[=\sigma_{i_1}\cdots \sigma_{i_{t-1}}\sigma_{i_{t}}\cdots \sigma_{i_{l-1}}\sigma_{j_l}\cdots \sigma_{j_{\ell(\pi)}}=\sigma_{i_1}\cdots \sigma_{i_{t-1}}\sigma_{i_{t+1}}\cdots \sigma_{i_{l-1}}\sigma_{j_l}\sigma_{j_l}\sigma_{j_{l+1}}\cdots \sigma_{j_{\ell(\pi)}}=\]\[=\sigma_{i_1}\cdots \sigma_{i_{t-1}}\sigma_{i_{t+1}}\cdots \sigma_{i_{l-1}}\sigma_{j_{l+1}}\cdots \sigma_{j_{\ell(\pi)}}.
\]
We obtained an expression for $\pi$ with $\ell(\pi)-2$ simple reflections, which is impossible.

Hence, we must have $\mathbf{x}_{P'}>_{lex}\mathbf{x}_{P_i}$. Thus, there exists a monomial $\mathbf{x}_{P'}$ and an integer $j_l\in[n]$ such that $\mathbf{x}_{P'}/\gcd(\mathbf{x}_{P'},\mathbf{x}_{P_i})=x_{j_l}$ and $x_{j_l}$ divides $\mathbf{x}_{P_{j}}/\gcd(\mathbf{x}_{P_i},\mathbf{x}_{P_j})$ which ends our proof.
\end{proof}

\begin{Example1}\label{ex}\rm$ \ $ Let $(S_4,S)$ be the Coxeter system and $Q$ be the following word of size $8$, $$Q=(s_1,\ s_2,\ s_1,\ s_3,\ s_1,\ s_2,\ s_3,\ s_1).$$ Let $\pi=(1,2,4)$ be an element in $S_4$, $\ell(\pi)=4$. The set of all the reduced expressions of $\pi$ is $$\{s_1s_2s_3s_2,\ s_1s_3s_2s_3,\ s_3s_1s_2s_3\}.$$ 
Let us denote  $$Q=(\sigma_1,\ \sigma_2,\ \sigma_3,\ \sigma_4,\ \sigma_5,\ \sigma_6,\ \sigma_7,\ \sigma_8).$$ The set of all the subwords of $Q$ that represent $\pi$ is $$\{(\sigma_1,\sigma_2,\sigma_4,\sigma_6),\ (\sigma_1,\sigma_4,\sigma_6,\sigma_7),\ (\sigma_3,\sigma_4,\sigma_6,\sigma_7),\ (\sigma_4,\sigma_5,\sigma_6,\sigma_7)\}.$$ 

Let $k[x_1,\ldots,x_8]$ be the polynomial ring over a field $k$. The Stanley--Reisner ideal of the Alexander dual of $\Delta$ is the square-free monomial ideal whose minimal monomial set of generators is $$G(I_{\Delta^{\vee}})=\{x_1x_2x_4x_6,\ x_1x_4x_6x_7,\ x_3x_4x_6x_7,\ x_4x_5x_6x_7\}.$$ We denote $\mathbf{x}_{P_1}=x_1x_2x_4x_6,\ \mathbf{x}_{P_2}=x_1x_4x_6x_7,\ \mathbf{x}_{P_3}=x_3x_4x_6x_7,\ \mathbf{x}_{P_4}=x_4x_5x_6x_7$ and we have $\mathbf{x}_{P_1}>_{lex}\ldots>_{lex}\mathbf{x}_{P_4}$. Since $$(\mathbf{x}_{P_1})\colon \mathbf{x}_{P_2}=(x_2),\ (\mathbf{x}_{P_1},\mathbf{x}_{P_2})\colon \mathbf{x}_{P_3}=(x_1)\ \mbox{and}\ (\mathbf{x}_{P_1},\mathbf{x}_{P_2},\mathbf{x}_{P_3})\colon \mathbf{x}_{P_4}=(x_1,x_3),$$ $I_{\Delta^{\vee}}$ has linear quotients with respect to this order of the monomials from $G(I_{\Delta^{\vee}})$.
\end{Example1}

By using the Theorem \ref{linquot}, we can state the following corollary:

\begin{Corollary1}\label{shell} Let $\Delta$ be the subword complex $\Delta(Q,\pi)$ and let $G(I_{\Delta^{\vee}})=\{\mathbf{x}_{P_1},\ldots,$ $\mathbf{x}_{P_r}\}$, where $\mathbf{x}_{P_1}>_{lex}\ldots>_{lex}\mathbf{x}_{P_r}$. Then $Q\setminus P_{1},\ldots, Q\setminus P_r$ is a shelling order on the facets of $\Delta$.
\end{Corollary1}
\begin{proof} Since $G(I_{\Delta^{\vee}})=\{\mathbf{x}_{P_1},\ldots,\mathbf{x}_{P_r}\}$, it follows that $Q\setminus P_1,\ldots,Q\setminus P_r$ are the facets of $\Delta$. Since $(Q\setminus P_i)\setminus(Q\setminus P_j)=P_j\setminus P_i$, the assertion follows from Corollary \ref{sqfreelinquot}. 
\end{proof}

\begin{Example1}\label{ex1}\rm For the subword complex from Example \ref{ex}, a shelling on the facets of $\Delta$ is $F_1=\{\sigma_3,\sigma_5,\sigma_7,\sigma_8\},\ F_2=\{\sigma_2,\sigma_3,\sigma_5,\sigma_8\},\ F_3=\{\sigma_1,\sigma_2,\sigma_5,\sigma_8\},\ F_4=\{\sigma_1,\sigma_2,\sigma_3,\sigma_8\}$.
\end{Example1}

\begin{Remark1}\rm The shelling from Corollary \ref{shell} for the subword complex $\Delta(Q,\pi)$ coincides with the shelling inductively constructed by vertex-decomposing the subword complex $\Delta(Q,\pi)$ (see, for instance, A. Bj\"orner and M.L. Wachs \cite[Theorem 11.3]{BjWa}).
\end{Remark1}

\begin{Example1}\rm We study the same subword complex as in Example \ref{ex}. We obtained that $F_1=\{\sigma_3,\sigma_5,\sigma_7,\sigma_8\},\ F_2=\{\sigma_2,\sigma_3,\sigma_5,\sigma_8\},\ F_3=\{\sigma_1,\sigma_2,\sigma_5,\sigma_8\},\ F_4=\{\sigma_1,\sigma_2,\sigma_3,\sigma_8\}$, is a shelling on the facets of $\Delta$ (see Example \ref{ex1}). We shall prove that the same shelling is obtained inductively by vertex-decomposing $\Delta$. 

Let $Q'=Q\setminus \sigma_1$. Since $\ell(\sigma_1\pi)<\ell(\pi)$, by the proof of Theorem \ref{vdshell}, one obtains that $$\lk(\sigma_1,\Delta)=\Delta(Q',\pi)=\langle\{\sigma_2,\sigma_5,\sigma_8\},\ \{\sigma_2,\sigma_3,\sigma_8\}\rangle$$ and $$\del(\sigma_1,\Delta)=\del(Q',\sigma_1\pi)=\langle\{\sigma_3,\sigma_5,\sigma_7,\sigma_8\},\ \{\sigma_2,\sigma_3,\sigma_5,\sigma_8\}\rangle.$$ We denote $\Delta_1=\langle\{\sigma_2,\sigma_5,\sigma_8\},\ \{\sigma_2,\sigma_3,\sigma_8\}\rangle$ and $\Delta_2=\langle\{\sigma_3,\sigma_5,\sigma_7,\sigma_8\},\ \{\sigma_2,\sigma_3,\sigma_5,\sigma_8\}\rangle$. We apply the same procedure to $\Delta_1$. Let $Q''=Q'\setminus\sigma_2$. Since $\ell(\sigma_2\pi)>\ell(\pi)$, we have $$\lk(\sigma_2,\Delta_1)=\del(\sigma_2,\Delta_1)=\Delta(Q'',\pi)=\langle\{\sigma_5,\sigma_8\},\ \{\sigma_3,\sigma_8\}\rangle.$$ 
Let us denote this simplicial complex by $\Delta_1'$ and $Q'''=Q''\setminus \sigma_2$. We have $\ell(\sigma_3\pi)<\ell(\pi)$. Hence, $$\lk(\sigma_3,\Delta_1')=\Delta(Q''',\pi)=\langle\{\sigma_8\}\rangle$$ and $$\del(\sigma_3,\Delta_1')=\Delta(Q''',\sigma_3\pi)=\Delta(Q''',s_2s_3s_2)=\langle\{\sigma_5,\sigma_8\}\rangle.$$

For the simplicial complex $\Delta_2$, since $\ell(\sigma_2\sigma_1\pi)<\ell(\sigma_1\pi)$, one has that $$\lk(\sigma_2,\Delta_2)=\Delta(Q'',\sigma_1\pi)=\Delta(Q'',s_2s_3s_2)=\langle\{\sigma_3,\sigma_5,\sigma_8\}\rangle$$ and $$\del(\sigma_2,\Delta_2)=\Delta(Q'',\sigma_2\sigma_1\pi)=\Delta(Q'',s_3s_2)=\langle\{\sigma_3,\sigma_5,\sigma_7,\sigma_8\}\rangle.$$

Hence, we get the following shelling on the facets of the subword complex $\Delta$
	\[\{\sigma_3,\sigma_5,\sigma_7,\sigma_8\},\ \{\sigma_2,\sigma_3,\sigma_5,\sigma_8\},\ \{\sigma_1,\sigma_2,\sigma_5,\sigma_8\},\ \{\sigma_1,\sigma_2,\sigma_3,\sigma_8\}
\]
which is the same shelling as the one obtained in the Example \ref{ex}. 
\end{Example1}

Note that the subword complexes are not shifted simplicial complexes, as one can see in the following example.

\begin{Example1}\rm Let $(S_4,S)$ be the Coxeter system and $Q$ be the word of size $6$, $$Q=(s_1,\ s_3,\ s_3,\ s_1,\ s_2,\ s_3).$$ Let, as before, $\pi=(1,2,4)\in S_4$ with $\ell(\pi)=4$. The set of all the reduced expressions of $\pi$ is $$\{s_1s_2s_3s_2,\ s_1s_3s_2s_3,\ s_3s_1s_2s_3\}.$$ We denote  $$Q=(\sigma_1,\ \sigma_2,\ \sigma_3,\ \sigma_4,\ \sigma_5,\ \sigma_6).$$ The set of all the subwords of $Q$ that represent $\pi$ is $$\{(\sigma_1,\sigma_2,\sigma_5,\sigma_6),\ (\sigma_1,\sigma_3,\sigma_5,\sigma_6),\ (\sigma_2,\sigma_4,\sigma_5,\sigma_6),\ (\sigma_3,\sigma_4,\sigma_5,\sigma_6)\}.$$ Hence the subword complex $\Delta=\Delta(Q,\pi)$ is the simplicial complex with the facets $$\mathcal{F}(\Delta)=\{\{\sigma_3,\sigma_4\},\ \{\sigma_2,\sigma_4\},\ \{\sigma_1,\sigma_3\},\ \{\sigma_1,\sigma_2\}\}.$$ If we consider a label of the vertices such that $\sigma_1<\sigma_3$, looking to the facet $\{\sigma_3,\sigma_4\}$ and replacing $\sigma_3$ by $\sigma_1$, we obtain the set $\{\sigma_1,\sigma_4\}$ which is not a face in $\Delta$. If we order the vertices such that $\sigma_3<\sigma_1$, looking to the facet $\{\sigma_1,\sigma_2\}$ and replacing $\sigma_1$ by $\sigma_3$, we obtain the set $\{\sigma_2,\sigma_3\}$ which is not a face in $\Delta$. Hence $\Delta$ is not a shifted simplicial complex.
\end{Example1}

The projective dimension of the Stanley--Reisner ring can be easily determined.
\begin{Proposition1} Let $\Delta$ be the subword complex $\Delta(Q,\pi)$ and let $n$ be the size of $Q$. Then $$\projdim(k[\Delta])=\ell(\pi).$$
\end{Proposition1} 

\begin{proof} By Theorem \ref{vdshell}, $\Delta$ is shellable. Since any subword $P\subseteq Q$ that represents $\pi$ is of size $\ell(\pi)$ and $\Delta$ is pure, we have that $\dim(\Delta)=n-\ell(\pi)-1$, so $\dim(k[\Delta])=n-\ell(\pi)$. Therefore, $\projdim(k[\Delta])=\ell(\pi)$.
\end{proof}

One may note that the projective dimension of the Stanley--Reisner ring of a subword complex does not depend on the word $Q$.

Next we determine all the elements of $\set(\mathbf{x}_P)$, where $\mathbf{x}_P\in G(I_{\Delta^{\vee}})$ where the monomials are ordered decreasing in the lexicographical order. 
We will need the following lemma. 
\begin{Lemma1}\label{snmid} Let $I$ be a square-free monomial ideal with $G(I)=\{w_1,\ldots,w_r\}$ and $w_1>_{lex}\ldots>_{lex}w_r$ such that $I$ has linear quotients with respect to this order of the generators. Then
	\[\set(w_i)\subseteq [\max(w_i)]\setminus \supp(w_i),
\]
for all $1\leq i\leq r$, where $[\max(w_i)]=\{1,2,\ldots,\max(w_i)\}$.
 \end{Lemma1}
 \begin{proof} Let $i\geq2$ and $k\in\set(w_i)$. Then $x_kw_i\in(w_1,\ldots,w_{i-1})$ that is there exist a variable $x_t$, $t\neq k$ and $j<i$ such that $w_ix_k=w_jx_t$. If $k\in\supp(w_i)$, since $k\neq t$, we get that $x_k^2\mid w_j$, contradiction. Thus $k\notin\supp(w_i)$. Since $t\neq k$, we have $x_t\mid w_i$. Then 
	\[w_j=\frac{x_kw_i}{x_t}>_{lex} w_i,
\]
which implies $k<t\leq\max(w_i)$.
 \end{proof}

In general, the inclusion is strict, as one can see in the following example.
\begin{Example1}\rm We consider the same subword complex as in Example \ref{ex}. We note that $\set(\mathbf{x}_{P_2})=\{2\}$, $\set(\mathbf{x}_{P_3})=\{1\}$ and $\set(\mathbf{x}_{P_4})=\{1,3\}$. We have $$[\max(\mathbf{x}_{P_4})]\setminus\supp\{\mathbf{x}_{P_4}\}=\{1,2,3,4,5,6,7\}\setminus\{4,5,6,7\}=\{1,2,3\},$$ and $\set(\mathbf{x}_{P_4})\subsetneq\{1,2,3\}$.
\end{Example1}

Therefore we will denote by $P_j\setminus P_i$ the word obtained from $P_j$ by omitting the simple reflections that appear both in $P_i$ and $P_j$. We also use $\min(P_j\setminus P_i)$ for $\min(\mathbf{x}_{P_j}/\gcd(\mathbf{x}_{P_j},\mathbf{x}_{P_i}))$. 

\begin{Proposition1}\label{set} Let $\Delta$ be the subword complex $\Delta(Q,\pi)$ and let $G(I_{\Delta^{\vee}})=\{\mathbf{x}_{P_1},\ldots,$ $\mathbf{x}_{P_r}\}$, where $\mathbf{x}_{P_1}>_{lex}\ldots>_{lex}\mathbf{x}_{P_r}$. For any $1\leq i\leq r$, we have $$\set(\mathbf{x}_{P_i})=\{\min(P_j\setminus P_i)\ \mid\ 1\leq j<i\}.$$
\end{Proposition1}

\begin{proof} We will use $I$ instead of $I_{\Delta^{\vee}} $ to simplify the notation. 

Let $s\in\set(\mathbf{x}_{P_i})$. We have to prove that there exists a monomial $\mathbf{x}_{P_j}$ such that $s=\min(P_j\setminus P_i)$.

Since $s\in\set(\mathbf{x}_{P_i})$, $x_s\mathbf{x}_{P_i}\in (\mathbf{x}_{P_1},\ldots,\mathbf{x}_{P_{i-1}})$. Hence, there exist $j<i$ and a variable $x_t$ such that 
	\[x_s\mathbf{x}_{P_i}=x_t\mathbf{x}_{P_j}.
\]
One may note that $s\neq t$, since $\mathbf{x}_{P_i}\neq\mathbf{x}_{P_j}$. By Lemma \ref{snmid}, we have that $x_s\nmid\mathbf{x}_{P_i}$. So we have that
\[\mathbf{x}_{P_j}=\frac{x_s\mathbf{x}_{P_i}}{x_t}.
\] Hence $|\{P_j\setminus P_i\}|=1$ and $s=\min(P_j\setminus P_i)$.

For the other inclusion, let $j<i$ and $s=\min(P_j\setminus P_i)$. By the proof of Theorem \ref{linquot}, we have that there exists a monomial $\mathbf{x}_{P}\in G(I)$, $\mathbf{x}_{P}\neq\mathbf{x}_{P_i}$, $\mathbf{x}_P>_{lex} \mathbf{x}_{P_i}$ and there exists a variable $t\in [n]$, $t\neq s$ such that $x_s\mathbf{x}_{P_i}=x_t\mathbf{x}_{P}\in(\mathbf{x}_{P_1},\ldots,\mathbf{x}_{P_{i-1}})$, which ends the proof. 
\end{proof}

We note that, for an arbitrary square-free monomial ideal which has linear quotients with respect to the lexicographical order of its minimal monomial generators $w_1,\ldots,w_r$, the equality $\set(w_i)=\{\min(\supp(w_j)\setminus\supp(w_i))\ :\ 1\leq j<i\}$ might be not true.

\begin{Example1}\rm Let $I=(x_1x_2x_3,\ x_2x_3x_4,\ x_2x_4x_5)$ be a square-free monomial ideal in the polynomial ring $k[x_1,\ldots,x_5]$. Denote $w_1=x_1x_2x_3,\ w_2= x_2x_3x_4,\ w_3=x_2x_4x_5$. One may note that $w_1>_{lex}w_2>_{lex}w_3$ and $I$ has linear quotients with respect to this order of generators. We have that $\set(w_2)=\{1\}$ and $\set(w_3)=\{3\}$. If we denote $F_i=\supp(w_i)$, $1\leq i\leq3$, we have that $\min(F_1\setminus F_3)=\{1\}$ and $\{1\}\notin\set(w_3)$.
\end{Example1}
 
Let $\Delta$ be the subword complex $\Delta(Q,\pi)$ and $G(I_{\Delta^{\vee}})=\{\mathbf{x}_{P_1},\ldots,\mathbf{x}_{P_r}\}$ be the minimal monomial set of generators for $I_{\Delta^{\vee}}$. For a monomial $\mathbf{x}_{P_i}$ from $G(I_{\Delta^{\vee}})$, we denote $d_i=|\set(\mathbf{x}_{P_i})|$. We note that $d_i\leq i-1$. 

\begin{Example1}\rm Let $(S_4,S)$ be the Coxeter system and $Q$ be the word of size $7$, $$Q=(s_1,\ s_1,\ s_1,\ s_3,\ s_1,\ s_2,\ s_3)$$and $\pi=(1,2,4)\in S_4$ with $\ell(\pi)=4$. Note that $s_1s_2s_3s_2,\ s_1s_3s_2s_3,\ s_3s_1s_2s_3$ are all the reduced expressions of $\pi$. 
We denote  $$Q=(\sigma_1,\ \sigma_2,\ \sigma_3,\ \sigma_4,\ \sigma_5,\ \sigma_6,\ \sigma_7).$$ The set of all the subwords of $Q$ that represent $\pi$ is $$\{(\sigma_1,\sigma_4,\sigma_6,\sigma_7),\ (\sigma_2,\sigma_4,\sigma_6,\sigma_7),\ (\sigma_3,\sigma_4,\sigma_6,\sigma_7),\ (\sigma_4,\sigma_5,\sigma_6,\sigma_7)\}.$$ 

Let $k[x_1,\ldots,x_7]$ be the polynomial ring over a field $k$. The Stanley--Reisner ideal of the Alexander dual of $\Delta$ has the minimal monomial set of generators $$G(I_{\Delta^{\vee}})=\{x_1x_4x_6x_7,\ x_2x_4x_6x_7,\ x_3x_4x_6x_7,\ x_4x_5x_6x_7\}.$$ Then $$\mathbf{x}_{P_1}=x_1x_4x_6x_7>_{lex} \mathbf{x}_{P_2}=x_2x_4x_6x_7>_{lex} \mathbf{x}_{P_3}=x_3x_4x_6x_7>_{lex} \mathbf{x}_{P_4}=x_4x_5x_6x_7.$$ Since $$(\mathbf{x}_{P_1})\colon \mathbf{x}_{P_2}=(x_1),\ (\mathbf{x}_{P_1},\mathbf{x}_{P_2})\colon \mathbf{x}_{P_3}=(x_1,x_2),\ (\mathbf{x}_{P_1},\mathbf{x}_{P_2},\mathbf{x}_{P_3})\colon \mathbf{x}_{P_4}=(x_1,x_2,x_3)$$ we have that $d_1=0,\ d_2=1,\ d_3=2,\ d_4=3$.
\end{Example1}
 
Next, we esablish an upper bound for the projective dimension of the Stanley--Reisner ideal of the Alexander dual.
\begin{Theorem1}\label{pd} Let $Q=(\sigma_1,\ldots,\sigma_n)$ be a word in $W$, $\pi$ be an element in $W$ and $\Delta$ be the subword complex $\Delta(Q,\pi)$. Then
	\[\projdim(I_{\Delta^{\vee}})\leq n-\ell(\pi).
\]
\end{Theorem1}
\begin{proof} Let $G(I_{\Delta^{\vee}})=\{\mathbf{x}_{P_1},\ldots,\mathbf{x}_{P_r}\}$, where $\mathbf{x}_{P_1}>_{lex}\ldots>_{lex}\mathbf{x}_{P_r}$ and $P_i\subseteq Q$ represents $\pi$ for all $1\leq i\leq r$. 

Since $I_{\Delta^{\vee}}$ has linear quotients with respect to the sequence $\mathbf{x}_{P_1},\ldots,\mathbf{x}_{P_r}$, we have $\projdim(I_{\Delta^{\vee}})=\max\{d_1,\ldots,d_r\}$. Let us assume by contradiction that $$\projdim(I_{\Delta^{\vee}})> n-\ell(\pi).$$ Hence, there exists $1\leq k\leq r$ such that $$\projdim(I_{\Delta^{\vee}})=d_k> n-\ell(\pi).$$
By Lemma \ref{snmid}, we have  
	\[\set(\mathbf{x}_{P_k})\cap\supp(\mathbf{x}_{P_k})=\emptyset.
\]
Since $\mathbf{x}_{P_k}$ is a square-free monomial, $|\supp(\mathbf{x}_{P_k})|=\ell(\pi)$.

We have that 
	\[|\set(\mathbf{x}_{P_k})|+|\supp(\mathbf{x}_{P_k})|>n-\ell(\pi)+\ell(\pi)=n
\]
that is 
	\[|\set(\mathbf{x}_{P_k})\cup\supp(\mathbf{x}_{P_k})|>n
\]
which is a contradiction.
\end{proof}

\begin{Remark1}\rm Let $Q=(\sigma_1,\ldots,\sigma_n)$ be a word in $W$, $\pi$ an element in $W$ and $\Delta$ the subword complex $\Delta(Q,\pi)$. Let $G(I_{\Delta^{\vee}})=\{\mathbf{x}_{P_1},\ldots,\mathbf{x}_{P_r}\}$ be the minimal monomial set of generators for $I_{\Delta^{\vee}}$ with $\mathbf{x}_{P_1}>_{lex}\ldots>_{lex}\mathbf{x}_{P_r}$. By Theorem \ref{pd}, we have that, if there exists $1\leq i\leq r$ such that $d_i=i-1$, then $i\leq n-\ell(\pi)+1$.
\end{Remark1}

As a consequence of Theorem \ref{pd}, we get an upper bound for the Castelnuovo--Mumford regularity of the Stanley--Reiner ideal of a subword complex.

\begin{Corollary1}Let $Q=(\sigma_1,\ldots,\sigma_n)$ be a word in $W$, $\pi\in W$ be an element and $\Delta$ the subword complex $\Delta(Q,\pi)$. Then
	\[\reg(I_{\Delta})\leq n-\ell(\pi)+1.
\]
\end{Corollary1}
\begin{proof} By Theorem \ref{Terai}, we have that 
	\[\reg(I_{\Delta})=\projdim k[\Delta^{\vee}].
\]
On the other hand, by Theorem \ref{pd}
	\[\projdim k[\Delta^{\vee}]=\projdim(I_{\Delta^{\vee}})+1\leq n-\ell(\pi)+1
\]
which ends our proof.
\end{proof}
The following results will be used in the next section. 

\begin{Lemma1}\label{min} Let $u,v,w$ be monomials of the same degree in $k[x_1,\ldots,x_n]$. Assume that $u,v>_{lex}w$ and $\min(u/\gcd(u,w))\neq \min(v/\gcd(v,w))$. Then $u>_{lex} v$ if and only if $\min(u/\gcd(u,w))$ $<\min(v/\gcd(v,w))$.
\end{Lemma1}
\begin{proof} In the following, for a monomial $m=x_1^{\alpha_1}\cdots x_n^{\alpha_n}$, we denote by $\nu_i(m)$ the exponent of the variable $x_i$ in $m$, that is $\nu_i(m)=\alpha_i$, $i=1,\ldots,n$. Since $u>_{lex}w$ there exists an integer $l'$ such that for all $i<l'$, $\nu_i(u)=\nu_i(w)$ and $\nu_{l'}(u)>\nu_{l'}(w)$. Similar, since $v>_{lex}w$ there exists an integer $l''$ such that for all $i<l''$, $\nu_i(v)=\nu_i(w)$ and $\nu_{l''}(v)>\nu_{l''}(w)$. Therefore, $l'=\min(u/\gcd(u,w))$ and $l''=\min(v/\gcd(v,w))$ and, by the hypothesis, $l'\neq l''$.

"$\Leftarrow$" Since $l'<l''$, we have that $\nu_{l'}(u)>\nu_{l'}(w)=\nu_{l'}(v)$ and for all $i<l'$ $\nu_{i}(u)=\nu_{i}(w)=\nu_{i}(v)$. Thus, $u>_{lex} v$.

"$\Rightarrow$" Since $u>_{lex}v$, there exists an integer $l\in [n]$ such that, for all $i<l$, $\nu_i(u)=\nu_i(v)$ and $\nu_{l}(u)>\nu_{l}(v)$.  Assume that $l<\min(l',l'')$. Then, one has $\nu_l(w)=\nu_l(v)<\nu_{l}(u)=\nu_l(w)$, which is impossible. Let us suppose that $\l>\min(l',l'')$. If we consider that $l'<l''$, we have $l>l'$ and we get $\nu_{l'}(w)<\nu_{l'}(u)=\nu_{l'}(v)=\nu_{l'}(w)$, contradiction. Similar, if $l''<l'$, we have $l>l''$ and we get $\nu_{l''}(w)<\nu_{l''}(v)=\nu_{l''}(u)=\nu_{l''}(w)$, contradiction.  

Therefore, we must have $l=\min(l',l'')$. Let us assume that $l''<l'$. Hence $l=l''$ and we get that $\nu_{l''}(u)=\nu_{l''}(w)<\nu_{l''}(v)<\nu_{l''}(u)$ which is impossible. Thus, we must have $l'<l''$, that is $\min(u/\gcd(u,w))<\min(v/\gcd(v,w))$.
\end{proof}

\begin{Lemma1}\label{i-1} Let $\Delta$ be the subword complex $\Delta(Q,\pi)$ and let $G(I_{\Delta^{\vee}})=\{\mathbf{x}_{P_1},$ $\ldots,$ $\mathbf{x}_{P_r}\}$ be the minimal monomial set of generators for $I_{\Delta^{\vee}}$, with $\mathbf{x}_{P_1}>_{lex}\ $ $\ldots>_{lex} \mathbf{x}_{P_r}$. Assume that there exists $2\leq i\leq r$ such that $d_i=i-1$. Then, for all $1\leq j<i$, $d_j=j-1$. 
\end{Lemma1}
\begin{proof} Since $d_i=i-1$, we have that $\min(P_j\setminus P_i)\neq\min(P_k\setminus P_i)$, for all $1\leq j,k<i$, $j\neq k$. Hence, by Lemma \ref{min}, 
	\[\min(P_1\setminus P_i)<\ldots<\min(P_{i-1}\setminus P_i)
\]

Let us fix $j<i$ and assume that $P_j=(\sigma_{j_1},\ldots,\sigma_{j_{\ell(\pi)}})$ and $P_i=(\sigma_{i_1},\ldots,\sigma_{i_{\ell(\pi)}})$. Since $\mathbf{x}_{P_j}>_{lex}\mathbf{x}_{P_i}$, we have that  
\begin{eqnarray}\mbox{for all}\ t<\min(P_j\setminus P_i)\ i_t=j_t\ \mbox{and}\ j_{\min(P_j\setminus P_i)}<i_{\min(P_j\setminus P_i)}.\label{p1}
\end{eqnarray}

Let $1\leq k<j$. We prove that $\min(P_k\setminus P_i)=\min(P_k\setminus P_j)$. This will imply that, for all $1\leq k,s<j$ with $k\neq s$, $\min(P_k\setminus P_j)\neq\min(P_s\setminus P_j)$ and hence $d_j=j-1$.

Since $k<j$, by Lemma \ref{min}, \begin{eqnarray}\min(P_k\setminus P_i)<\min(P_j\setminus P_i).\label{p2}
\end{eqnarray}

On the other hand, since $k<i$, we have that $\mathbf{x}_{P_k}>_{lex}\mathbf{x}_{P_i}$ that is for all $t<\min(P_k\setminus P_i)$ $i_t=k_t$ and $k_{\min(P_k\setminus P_i)}<i_{\min(P_k\setminus P_i)}$. By (\ref{p1})and (\ref{p2}) we have that  $k_{\min(P_k\setminus P_i)}<j_{\min(P_k\setminus P_i)}$ and for all $t<\min(P_k\setminus P_i)$ $k_t=j_t$, which means that $\min(P_k\setminus P_i)=\min(P_k\setminus P_j)$.
\end{proof}

\begin{Lemma1}\label{unic l} Let $\Delta$ be the subword complex $\Delta(Q,\pi)$ and let $G(I_{\Delta^{\vee}})=\{\mathbf{x}_{P_1},$ $\ldots,$ $\mathbf{x}_{P_r}\}$ be the minimal monomial set of generators for $I_{\Delta^{\vee}}$, with $\mathbf{x}_{P_1}>_{lex}\ldots>_{lex}\mathbf{x}_{P_r}$. Let $2\leq i\leq r$. Then $d_i=i-1$ if and only if there exists a unique $l\in\supp(\mathbf{x}_{P_i})$ such that $\mathbf{x}_{P_j}=x_{\min(P_j\setminus P_i)}\mathbf{x}_{P_i}/x_l$, for all $1\leq j<i$. 
\end{Lemma1}
\begin{proof} "$\Rightarrow$" Let us assume that $d_i=i-1$. Then $\min(P_j\setminus P_i)\neq\min(P_k\setminus P_i)$ for all $1\leq j,k<i$ with $j\neq k$. 

By Proposition \ref{set}, we have that $$(\mathbf{x}_{P_1},\ldots,\mathbf{x}_{P_{i-1}}):\mathbf{x}_{P_i}=(x_{\min(P_j\setminus P_i)}\ :\ 1\leq j<i).$$
Let $j<i$ and $1\leq k<j<i$. Let us assume by contradiction that there exist $i_t,i_{t'}\in\supp(\mathbf{x}_{P_i})$, $i_t\neq i_{t'}$ such that
	\[\mathbf{x}_{P_j}=x_{\min(P_j\setminus P_i)}\frac{\mathbf{x}_{P_i}}{x_{i_t}},\ \mathbf{x}_{P_k}=x_{\min(P_k\setminus P_i)}\frac{\mathbf{x}_{P_i}}{x_{i_{t'}}}.
\]

By the proof of Lemma \ref{i-1}, we have that $\min(P_k\setminus P_i)=\min(P_k\setminus P_j)$. Since $j<i$, by Lemma \ref{i-1}, $d_j=j-1$ and there exists $j_{t''}\in\supp(\mathbf{x}_{P_j})$ such that
	\[\mathbf{x}_{P_k}=x_{\min(P_k\setminus P_j)}\frac{\mathbf{x}_{P_j}}{x_{j_{t''}}},
\]
Replacing $\mathbf{x}_{P_k}$, we have that $\mathbf{x}_{P_j}x_{i_{t'}}=\mathbf{x}_{P_i}x_{j_{t''}}$. We note that $x_{i_{t'}}\neq x_{j_{t''}}$ since $\mathbf{x}_{P_j}\neq\mathbf{x}_{P_i}$. Hence, we have that
	\[x_{i_{t'}}x_{\min(P_j\setminus P_i)}\frac{\mathbf{x}_{P_i}}{x_{i_{t}}}=\mathbf{x}_{P_i}x_{j_{t''}}
\]
that is
	\[x_{i_{t'}}x_{\min(P_j\setminus P_i)}=x_{j_{t''}}x_{i_t}.
\]
Since $x_{i_{t'}}\mid x_{j_{t''}}x_{i_t}$ and $x_{i_{t'}}\neq x_{j_{t''}}$, we have that $x_{i_{t'}}= x_{i_{t}}$ contradiction with our assumption.

"$\Leftarrow$" If there exists a unique $l\in\supp(\mathbf{x}_{P_i})$ such that $\mathbf{x}_{P_j}=x_{\min(P_j\setminus P_i)}\mathbf{x}_{P_i}/x_l$ for all $1\leq j<i$ we have that $\min(P_j\setminus P_i)\in\set(\mathbf{x}_{P_i})$ for all $1\leq j\leq i-1$. Since $\mathbf{x}_{P_j}\neq\mathbf{x}_{P_k}$ for all $1\leq j,k\leq i-1$, $j\neq k$, we get that $d_i=i-1$.
\end{proof}
\section{A special class of subword complexes}

In this section we consider only subword complexes $\Delta=\Delta(Q,\pi)$ such that the minimal monomial generating system of $I_{\Delta^{\vee}}$ has $r\leq n-\ell(\pi)+1$ elements, where $n$ is the size of $Q$, and for which $d_r=r-1$. In the following proposition, we construct classes of such subword complexes.

\begin{Proposition1} Let $\pi\in W$ be an element and $\sigma_1\cdots\sigma_{\ell(\pi)}$ a reduced expression for $\pi$. Let $1\leq i\leq \ell(\pi)$ be a fixed integer and let $$Q=(\sigma_1,\sigma_2,\ldots,\sigma_{i-1},\sigma_i,\sigma_i,\ldots,\sigma_i,\sigma_{i+1},\ldots,\sigma_{\ell(\pi)})$$ be a word of size $n$ in $W$. Then the minimal monomial generating system of $I_{\Delta^{\vee}}$ has $N=n-\ell(\pi)+1$ elements and $d_N=N-1$.
\end{Proposition1}
\begin{proof} Since $\sigma_1\cdots\sigma_{\ell(\pi)}$ is a reduced expression for $\pi$, any subword of $Q$ that represents $\pi$ is a copy of this reduced expression. So, $I_{\Delta^{\vee}}$ has the minimal monomial generating system $G(I_{\Delta^{\vee}})=\{x_1\cdots x_{i-1}x_jx_{i+N}\cdots x_{n}|\ i\leq j\leq N+i-1\}$, where $N= n-\ell(\pi)+1$. Hence, $|G(I_{\Delta^{\vee}})|=N$. One may note that $d_N=N-1$.
\end{proof}

Not all the words $Q$ such that, for the subword complex $\Delta=\Delta(Q,\pi)$, the minimal monomial generating system of $I_{\Delta^{\vee}}$ has $r\leq n-\ell(\pi)+1$ elements, where $n$ is the size of $Q$, and for which $d_r=r-1$ have the same form as in the above proposition.

\begin{Example1}\rm Let $(S_4,S)$ be the Coxeter system, $Q$ the following word of size $9$
	\[Q=(s_2,\ s_3,\ s_2,\ s_3,\ s_1,\ s_3,\ s_2,\ s_3,\ s_2)
\]
and $\pi=(14)(23)\in S_4$, $\ell(\pi)=6$. The set of all the reduced expressions of $\pi$ is
	\[\{s_2s_3s_2s_1s_2s_3,\ s_3s_2s_3s_1s_2s_3,\ s_3s_2s_1s_3s_2s_3,\ s_3s_2s_1s_2s_3s_2\}.
\]Let us denote $$Q=(\sigma_1,\sigma_2,\sigma_3,\sigma_4,\sigma_5,\sigma_6,\sigma_7,\sigma_8,\sigma_9)$$ All the subwords of $Q$ that represent $\pi$ are
	\[(\sigma_1,\sigma_2,\sigma_3,\sigma_5,\sigma_7,\sigma_8), (\sigma_2,\sigma_3,\sigma_4,\sigma_5,\sigma_7,\sigma_8), (\sigma_2,\sigma_3,\sigma_5,\sigma_6,\sigma_7,\sigma_8),(\sigma_2,\sigma_3,\sigma_5,\sigma_7,\sigma_8,\sigma_9), 
	\]
	hence
	\[G(I_{\Delta^{\vee}})=\{\mathbf{x}_{P_1}=x_1x_2x_3x_5x_7x_8,\ \mathbf{x}_{P_2}=x_2x_3x_4x_5x_7x_8,\ \mathbf{x}_{P_3}=x_2x_3x_5x_6x_7x_8,\]
	\[ \mathbf{x}_{P_4}=x_2x_3x_5x_7x_8x_9\}.
\]One may note that $\mathbf{x}_{P_1}>_{lex}\mathbf{x}_{P_2}>_{lex}\mathbf{x}_{P_3}>_{lex}\mathbf{x}_{P_4}$ and $I_{\Delta^{\vee}}$ has linear quotients with respect to the sequence $\mathbf{x}_{P_1},\mathbf{x}_{P_2},\mathbf{x}_{P_3},\mathbf{x}_{P_4}$. We have that $d_2=1,\ d_3=2$ and $d_4=3$.
\end{Example1}

For this special class of subword complexes, the Stanley--Reisner ideal of the Alexander dual is of a certain form.

\begin{Lemma1}\label{idelta}Let $\Delta$ be the subword complex $\Delta(Q,\pi)$ and suppose that the size of $Q$ is $n$. Assume that $G(I_{\Delta^{\vee}})=\{\mathbf{x}_{P_1},\ldots,\mathbf{x}_{P_r}\}$ with $\mathbf{x}_{P_1}>_{lex}\ldots>_{lex}\mathbf{x}_{P_r}$, $r\leq n-\ell(\pi)+1$ and $d_r=r-1$. Then there exists a unique $l\in\supp(\mathbf{x}_{P_r})$ such that
\[I_{\Delta^{\vee}}=\frac{\mathbf{x}_{P_r}}{x_l}(x_{\min(P_1\setminus P_r)},\ldots,x_{\min(P_{r-1}\setminus P_r)},x_l).
\]
\end{Lemma1}
\begin{proof} Since $d_r=r-1$, by Lemma \ref{unic l} we have that there exists a unique $l$ such that $\mathbf{x}_{P_j}=x_{\min(P_j\setminus P_r)}\mathbf{x}_{P_r}/x_l$ for all $1\leq j\leq r-1$. Hence,
	\[I_{\Delta^{\vee}}=\frac{\mathbf{x}_{P_r}}{x_l}(x_{\min(P_1\setminus P_r)},\ldots,x_{\min(P_{r-1}\setminus P_r)},x_l).
\]
\end{proof}

As a consequence, for this class one can compute the height of the Stanley--Reisner ideal of the Alexander dual.
\begin{Corollary1}\label{ht} Let $\Delta$ be the subword complex $\Delta(Q,\pi)$ and suppose that the size of $Q$ is $n$. Assume that $G(I_{\Delta^{\vee}})=\{\mathbf{x}_{P_1},\ldots,\mathbf{x}_{P_r}\}$ with $\mathbf{x}_{P_1}>_{lex}\ldots>_{lex}\mathbf{x}_{P_r}$, $r\leq n-\ell(\pi)+1$ and $d_r=r-1$. Then $\height(I_{\Delta^{\vee}})=1$. 
\end{Corollary1}

\begin{proof} By Lemma \ref{idelta}, 
\[I_{\Delta^{\vee}}=\frac{\mathbf{x}_{P_r}}{x_l}(x_{\min(P_1\setminus P_r)},\ldots,x_{\min(P_{r-1}\setminus P_r)},x_l).
\]
One may note that for any $t\in\supp(\mathbf{x}_{P_r}/x_l)$, the ideal $(x_t)$ is a minimal prime ideal of $I_{\Delta^{\vee}}$. Hence $\height(I_{\Delta^{\vee}})=1$.
\end{proof}

Henceforth, we denote $R=k[x_1,\ldots,x_n]$. We determine the minimal graded free resolution for the Stanley--Reisner ring of the Alexander dual.
\begin{Theorem1}\label{Koszul} Let $\Delta$ be the subword complex $\Delta(Q,\pi)$ and let $n$ be the size of $Q$. Assume that $G(I_{\Delta^{\vee}})=\{\mathbf{x}_{P_1},\ldots,\mathbf{x}_{P_r}\}$ with $\mathbf{x}_{P_1}>_{lex}\ldots>_{lex}\mathbf{x}_{P_r}$, $r\leq n-\ell(\pi)+1$ and $d_r=r-1$. Then there exists a unique integer $l\in[n]$ such that the Koszul complex associated to the sequence $$x_{\min(P_1\setminus P_r)},\ldots,x_{\min(P_{r-1}\setminus P_r)},x_l$$ is isomorphic to the minimal graded free resolution of $k[\Delta^{\vee}]$.
\end{Theorem1}

\begin{proof} By Lemma \ref{idelta} there exists a unique $l\in[n]$ such that 
	\[I_{\Delta^{\vee}}=\frac{\mathbf{x}_{P_r}}{x_l}(x_{\min(P_1\setminus P_r)},\ldots,x_{\min(P_{r-1}\setminus P_r)},x_l).
\] Therefore the multiplication by $\mathbf{x}_{P_r}/x_l$ defines an isomorphism of $R$--modules be- tween $(x_{\min(P_1\setminus P_r)},\ldots,x_{\min(P_{r-1}\setminus P_r)},x_l)$ and $I_{\Delta^{\vee}}$. Since $$x_{\min(P_1\setminus P_r)},\ldots,\ x_{\min(P_{r-1}\setminus P_r)},\ x_l$$ is a regular sequence, the assertion follows by Theorem \ref{Koszulres}.
\end{proof}

The following result is a simple consequence of Theorem \ref{Koszul}.
\begin{Corollary1}\label{Betti}Let $\Delta$ be the subword complex $\Delta(Q,\pi)$ and suppose that the size of $Q$ is $n$. Assume that $G(I_{\Delta^{\vee}})=\{\mathbf{x}_{P_1},\ldots,\mathbf{x}_{P_r}\}$ with $\mathbf{x}_{P_1}>_{lex}\ldots>_{lex}\mathbf{x}_{P_r}$, $r\leq n-\ell(\pi)+1$ and $d_r=r-1$. Then
	\[\beta_i(I_{\Delta^\vee})=\left(\twoline{r}{i+1} \right), 
\]
for all $i$.
\end{Corollary1}

Now we may get the Hilbert numerator of the Hilbert series.
\begin{Corollary1}\label{hilbd} Let $\Delta$ be the subword complex $\Delta(Q,\pi)$ and let $n$ be the size of $Q$. Assume that $G(I_{\Delta^{\vee}})=\{\mathbf{x}_{P_1},\ldots,\mathbf{x}_{P_r}\}$ with $\mathbf{x}_{P_1}>_{lex}\ldots>_{lex}\mathbf{x}_{P_r}$, $r\leq n-\ell(\pi)+1$ and $d_r=r-1$. Then the Hilbert numerator of the Hilbert series is
	\[\mathcal{K}_{I_{\Delta^{\vee}}}(t)=\sum_{i=0}^{r-1}(-1)^i\left(\twoline{r}{i+1} \right)t^{i+\ell(\pi)}.
\]
\end{Corollary1}
\begin{proof} Since $I_{\Delta^{\vee}}$ has a $\ell(\pi)$--linear resolution, $\beta_{ij}(I_{\Delta^{\vee}})=0$ for all $j\neq i+\ell(\pi)$. Hence $\beta_i(I_{\Delta^{\vee}})=\beta_{i,i+\ell(\pi)}(I_{\Delta^{\vee}})$. Since $\projdim(I_{\Delta^{\vee}})=r-1$,
\[\mathcal{K}_{I_{\Delta^{\vee}}}(t)=\sum_{i=0}^{r-1}(-1)^i\beta_i(I_{\Delta^{\vee}})t^{i+\ell(\pi)}.
\]
Hence, using Corollary \ref{Betti}, we get that the Hilbert numerator is 
\[\mathcal{K}_{I_{\Delta^{\vee}}}(t)=\sum_{i=0}^{r-1}(-1)^i\left(\twoline{r}{i+1} \right)t^{i+\ell(\pi)}.
\]
\end{proof}

\begin{Corollary1}\label{numw} Let $Q$ be a word in $W$ of size $n$ that contains $\pi$, $\Delta$ the subword complex  $\Delta(Q,\pi)$ and $G(I_{\Delta^{\vee}})=\{\mathbf{x}_{P_1},\ldots,\mathbf{x}_{P_r}\}$ with $\mathbf{x}_{P_1}>_{lex}\ldots>_{lex}\mathbf{x}_{P_r}$, $r\leq n-\ell(\pi)+1$ and $d_r=r-1$. Then there are $\left(\twoline{r}{j+1} \right)$ subwords $P$ of $Q$ such that $\delta(P)=\pi$ and $|P|=j+\ell(\pi)$ for $0\leq j\leq r-1$.
\end{Corollary1}
\begin{proof} By Lemma \ref{hilb}, in the fine grading, the Hilbert numerator of $I_{\Delta^\vee}$ is
\begin{eqnarray*}
  \mathcal{K}_{I_{\Delta^\vee}}( t_1,\ldots,t_n) &=& \sum_\twoline{P \subseteq Q}
  {\delta(P)=\pi} (-1)^{|P|-\ell(\pi)} \mathbf{t}^P.
\end{eqnarray*}
where $\mathbf{t}^P=\prod\limits_{\sigma_i\in P}t_i$.
In $\mathbb{Z}$--grading, we have 
\[
  \mathcal{K}_{I_{\Delta^\vee}}(t) = \sum_\twoline{P \subseteq Q}
  {\delta(P)=\pi} (-1)^{|P|-\ell(\pi)} t^{|P|}=\sum_\twoline{P \subseteq Q}
  {\delta(P)=\pi}\sum^n_\twoline{j=\ell(\pi)}
{|P|=j} (-1)^{j-\ell(\pi)} t^j=\]
  \[=\sum_\twoline{P \subseteq Q}
  {\delta(P)=\pi}\sum^{n-\ell(\pi)}_\twoline{j=0}
  {|P|=j+\ell(\pi)} (-1)^{j} t^{j+\ell(\pi)}=\sum_{j=0}^{n-\ell(\pi)}(-1)^j m_{j+\ell(\pi)} t^{j+\ell(\pi)},\]
where we denoted by $m_{j+\ell(\pi)}$ the number of the subwords $P$ of $Q$ such that $\delta(P)=\pi$ and $|P|=j+\ell(\pi)$. Comparing with the formula from Corollary \ref{hilbd}, we obtain that $m_{j+\ell(\pi)}=\left(\twoline{r}{j+1} \right)$ for all $0\leq j\leq r-1$, and $m_{j+\ell(\pi)}=0$ for all $r\leq j\leq n-\ell(\pi)$.
\end{proof}

Next, we characterize all the subword complexes from this special class which are simplicial spheres.
\begin{Corollary1} Let $Q$ be a word in $W$ of size $n$ that contains $\pi$, $\Delta$ the subword complex $\Delta(Q,\pi)$ and $G(I_{\Delta^{\vee}})=\{\mathbf{x}_{P_1},\ldots,\mathbf{x}_{P_r}\}$ with $\mathbf{x}_{P_1}>_{lex}\ldots>_{lex}\mathbf{x}_{P_r}$, $r\leq n-\ell(\pi)+1$ and $d_r=r-1$. Then $\Delta$ is a simplicial sphere if and only if $r=n-\ell(\pi)+1$.
\end{Corollary1}

\begin{proof} By Theorem \ref{sphere}, $\Delta$ is a simplicial sphere if $\delta(Q)=\pi$. Hence, in the Hilbert numerator, the coefficient of $t^{|Q|}$ must be non-zero. By Corollary \ref{numw}, the coefficient of $t^n$ is $m_n=\left(\twoline{r}{n-\ell(\pi)+1}\right)$. Hence $m_n\neq 0$ if and only if $r=n-\ell(\pi)+1$.
\end{proof}

\begin{Proposition1}\label{complint} Let $\Delta$ be the subword complex $\Delta(Q,\pi)$ and let $n$ be the size of $Q$. Assume that $G(I_{\Delta^{\vee}})=\{\mathbf{x}_{P_1},\ldots,\mathbf{x}_{P_r}\}$ with $\mathbf{x}_{P_1}>_{lex}\ldots>_{lex}\mathbf{x}_{P_r}$, $r \leq n-\ell(\pi)+1$ and $d_r=r-1$. Then $k[\Delta]$ is a complete intersection ring.
\end{Proposition1}
\begin{proof} By Lemma \ref{idelta}, we have that there exists a unique integer $l$ such that \[I_{\Delta^{\vee}}=\frac{\mathbf{x}_{P_r}}{x_l}(x_{\min(P_1\setminus P_r)},\ldots,x_{\min(P_{r-1}\setminus P_r)},x_l).
\]
Hence
\[I_{\Delta^{\vee}}=\left(\bigcap_{k\in\supp\left(\mathbf{x}_{P_r}/x_l\right)}(x_k)\right)\cap\left(x_{\min(P_1\setminus P_r)},\ldots, x_{\min(P_{r-1}\setminus P_r)},x_l\right)
\]
and
	\[I_{\Delta}=\left(x_{\min(P_1\setminus P_r)}\cdots x_{\min(P_{r-1}\setminus P_r)}x_l\right)+\left(x_k:k\in\supp\left(\mathbf{x}_{P_r}/x_l\right)\right).
\]
Since $\supp(\mathbf{x}_{P_r}/x_l)\cap\left\{\min(P_1\setminus P_r),\ldots, \min(P_{r-1}\setminus P_r),l\right\}=\emptyset$, $I_{\Delta}$ is a complete intersection ideal. 
\end{proof}

As a consequence, we get the minimal graded free resolution for the Stanley--Reisner ring of a subword complex from this class.
\begin{Corollary1} Let $\Delta$ be the subword complex $\Delta(Q,\pi)$ and suppose that the size of $Q$ is $n$. Assume that $G(I_{\Delta^{\vee}})=\{\mathbf{x}_{P_1},\ldots,\mathbf{x}_{P_r}\}$ with $\mathbf{x}_{P_1}>_{lex}\ldots>_{lex}\mathbf{x}_{P_r}$, $r\leq n-\ell(\pi)+1$ and $d_r=r-1$. Then there exists a unique integer $l$ such that the Koszul complex associated to the sequence $x_{\min(P_1\setminus P_r)}\ldots x_{\min(P_{r-1}\setminus P_r)}x_l,\ x_i\ :i\in\supp(\mathbf{x}_{P_r}/x_l)$ is a minimal graded free resolution of $k[{\Delta}]$.
\end{Corollary1}
\begin{proof} By Proposition \ref{complint}, we have that $$G(I_{\Delta})=\{x_{\min(P_1\setminus P_r)}\ldots x_{\min(P_{r-1}\setminus P_r)}\}\cup\{x_l,\ x_i\ :i\in\supp(\mathbf{x}_{P_r}/x_l)\}.$$ The statement follows by Theorem \ref{Koszulres}.
\end{proof}

We can describe all the subword complexes from this special class such that the Stanley--Reisner ring of the Alexander dual is Cohen--Macaulay.
\begin{Proposition1} Let $\Delta$ be the subword complex $\Delta(Q,\pi)$ and let $n$ be the size of $Q$. Assume that $G(I_{\Delta^{\vee}})=\{\mathbf{x}_{P_1},\ldots,\mathbf{x}_{P_r}\}$ with $\mathbf{x}_{P_1}>_{lex}\ldots>_{lex}\mathbf{x}_{P_r}$, $r\leq n-\ell(\pi)+1$ and $d_r=r-1$. Then $k[\Delta^{\vee}]$ is Cohen--Macaulay if and only if $I_{\Delta^{\vee}}$ is a principal monomial ideal.
\end{Proposition1}
\begin{proof} By Eagon--Reiner theorem \cite{EaRe}, $k[\Delta^{\vee}]$ is Cohen--Macaulay if and only if $I_{\Delta}$ has a linear resolution. In particular, $I_{\Delta}$ is generated in one degree. The statement follows by Proposition \ref{complint}. 
\end{proof}
\addcontentsline{toc}{chapter}{Ideas for future}
\chapter*{Ideas for future}
Alexander duality has become an important tool in the study of the square-free monomial ideals, due to the remarkable result of J.A. Eagon and V.Reiner \cite{EaRe} and extensions by N. Terai \cite{T}, relating data of the resolution of the Stanley--Reisner ring of a simplicial complex to that of its Alexander dual. The Eagon--Reiner theorem says that the ideal $I_{\Delta}$ associated to a simplicial complex $\Delta$ has a linear resolution if and only if the Alexander dual of $\Delta$ is a Cohen--Macaulay complex. In the same spirit, it was proved that $\Delta$ is shellable if and only if the Stanley--Reisner ideal associated to the Alexander dual of $\Delta$ has linear quotients \cite{HeHiZh} and that $\Delta$ is constructible if and only if the Stanley--Reisner ideal associated to the Alexander dual of $\Delta$ is constructible \cite{O1}.
 
Along these ideas, we want to find classes of examples of ideals which have linear resolutions but does not have linear quotients. We are also looking for classes of simplicial complexes which are constructible but they are not shellable and simplicial complexes which are Cohen--Macaulay but they are not shellable or constructible. 
 
As it was seen in the work of R. Stanley, the Hilbert function is a bridge between combinatorics and commutative algebra. Lexsegment ideals, as well as square-free lexsegment ideals, play the key role in the discussion. Based on the formula of S. Eliahou and M. Kervaire \cite{EK} and basic techniques on generic initial ideals, in 1993, A. Bigatti \cite{B} and H. Hulett \cite{Hul} independently obtained the result that among all the graded ideals with a given Hilbert function, the lexsegment ideal possesses the maximal graded Betti numbers, provided the base field is of characteristic zero. In arbitrary characteristic this has been shown by K. Pardue. This theorem has as square-free analogue a theorem due to A. Aramova, J. Herzog, and T. Hibi \cite{AHH}. We studied general lexsegment ideals and we proved that all the lexsegment ideals with a linear resolution are in fact ideals with linear quotients. We also computed some invariants such as the dimension and the depth and we characterized all the lexsegment ideals which are Cohen--Macaulay \cite{EOS}.
 
We want to carry out a similar study of arbitrary square-free lexsegment ideals with regard to their resolutions and other  invariants like  depth, dimension, Cohen--Macaulay-ness, etc. We also want to classify all the lexsegment ideals and square-free lexsegment ideals which are Gotzmann. 

\backmatter
\fancyhead{}

\end{document}